\documentclass[11pt]{article}
\usepackage{amsmath,bm,amssymb,amsthm,amsfonts,afterpage,hyperref,
xcolor,graphicx,soul,multirow,multicol,siunitx}
\usepackage{mathtools}

\usepackage[normalem]{ulem} 
\usepackage{cancel} 
\usepackage[textsize=tiny,color=blue!50!white]{todonotes} 
\hypersetup{colorlinks,linkcolor=blue,citecolor=blue,urlcolor=blue}
\definecolor{darkgreen}{rgb}{0.01, 0.75, 0.24}
\DeclareMathAlphabet\mathbfcal{OMS}{cmsy}{b}{n}
\newtheorem{lemma}{Lemma}[section]

\newtheorem{assumption}[lemma]{Assumption}
\newtheorem{theorem}[lemma]{Theorem}

\newtheorem{remark}[lemma]{Remark}
\newtheorem{proposition}[lemma]{Proposition}
\newtheorem{alg}[lemma]{Algorithm}

\newcommand\curl{{\rm \bf curl}}
\renewcommand\div{\rm div}
\newcommand\B{\mathcal{B}}
\newcommand\I{\mathcal{I}}
\newcommand\bI{\bm{I}}
\newcommand\dif{{\bm K}}
\newcommand\T{\mathcal{T}}
\newcommand\V{\mathcal{V}}
\newcommand{\ver}{{\bm{a}}}
\newcommand\tb{{\bm b}}
\newcommand\omaj{\omega_{j}^\ver}
\newcommand\omaJ{\omega_{J}^\ver}
\newcommand\omaH{\omega_{H}^{\ver}}

\newcommand\N{\bm{\mathcal{N}}\hspace{-0.05em}}
\renewcommand\H{\bm{H}}
\newcommand\bV{\bm{V}}

\newcommand\bX{\bm{X}}
\newcommand\bVaj{\bV_{j}^\ver}
\newcommand\Vash{\bV_{h}^\ver}
\newcommand\bVajz{\bV_{j}^{\ver,0}}
\newcommand\Vaz{\bV_{h}^{\ver,0}}

\renewcommand\u{\bm{u}}
\newcommand\bv{\bm{v}}
\newcommand\bw{\bm{w}}
\newcommand\bz{\bm{z}}
\newcommand\bx{\bm{x}}
\newcommand\br{\bm{r}}
\newcommand\n{\bm{n}}

\newcommand\bxi{\bm{\xi}}

\newcommand\bzeta{\bm{\zeta}}
\newcommand\bphi{\bm{\phi}}
\newcommand\brho{\bm{\rho}}

\newcommand\balpha{\bm{\alpha}}
\newcommand\bdelta{\bm{\delta}}
\newcommand\rhoa{\brho_{h,\ver}^i}
\newcommand\bvja{\bv_{j,\ver}}
\newcommand\va{\bv_{h,\ver}}
\newcommand\etalg{\eta^i_{\mathrm{alg}}}
\newcommand\refvs[3]{\stackrel{\mathmakebox[\widthof{#3}]{#1}}{#2}}
\newcommand\eq{:=}
\newcommand\qe{=:}
\newcommand\RT{\bm{\mathcal{R\hspace{-0.1em}T\hspace{-0.05em}}}}
\newcommand\BDM{\bm{\mathcal{B\hspace{-0.05em}D\hspace{-0.1em}M}}}
\newcommand\PP{{\mathcal P}}

\title{
A-posteriori-steered $p$-robust multigrid and domain decomposition 
methods with optimal step-sizes for mixed finite element discretizations 
of elliptic problems\thanks{This project has received funding from 
Austrian Science Fund (FWF) project
\href{https://www.fwf.ac.at/en/research-radar/10.55776/F65}{10.55776/F65} 
(SFB F65 ``Taming complexity in PDE systems''),
the scientific mobility grant Christiana H\"orbiger 2021, 
Czech Academy of Sciences RVO 67985840, 
the Grant Agency of the Czech Republic grant no. 23-06159S, 
Inria Paris Visiting Professorship, NSF grants DMS 2111129 and DMS 241068. 
Jan Papež is a member of Nečas Center for Mathematical Modeling.
}
}
\author{
Ani Mira{\c c}i\footnotemark[5]
\and
Jan Pape{\v z}\footnotemark[6]
\and
Martin Vohral{\'i}k\footnotemark[2] \footnotemark[3]
\and
Ivan Yotov\footnotemark[4]
}
\date{\today}

\begin{document}

\maketitle
\renewcommand{\thefootnote}{\fnsymbol{footnote}}
\footnotetext[5]{TU Wien, Institute of Analysis and Scientific Computing, 
Wiedner Hauptstr. 8-10/E101/4, 1040 Vienna, Austria}
\footnotetext[6]{Institute of Mathematics, Czech Academy of Sciences, 
\v{Z}itn\'{a}~25, 115\,67 Prague, Czech Republic}
\footnotetext[2]{Inria, 48 rue Barrault, 75647 Paris, France}
\footnotetext[3]{CERMICS, Ecole nationale des ponts et chauss\'ees, IP 
Paris, 77455 Marne-la-Vall\'ee, France}
\footnotetext[4]{Department of Mathematics, University of Pittsburgh, 
Pittsburgh, PA 15260 USA}
\renewcommand{\thefootnote}{\arabic{footnote}}

\begin{abstract}
In this work, we develop algebraic solvers
for linear systems arising from the discretization of second-order elliptic 
partial differential equations by saddle-point mixed finite element 
methods of arbitrary polynomial degree $p \ge 0$ on possibly highly graded 
simplicial meshes. We present a multigrid and a two-level domain decomposition 
approach in two and three space dimensions, which are steered by their 
respective a~posteriori estimators of the algebraic error.
First, we extend the results of [Miraçi, Papež, and Vohral{\'i}k, 
\textit{SIAM J. Sci. Comput.} 43 (2021), S117–S145] to the mixed finite element 
setting.
Extending the multigrid procedure itself is rather natural. 
To obtain analogous theoretical results, however, a $p$-robust multilevel 
stable decomposition of the velocity space is needed.
In two space dimensions, we can treat the velocity space as the curl of a 
stream-function Lagrange space, for which the previous results apply. 
In three space 
dimensions, we design a novel stable decomposition by combining a 
one-level high-order local stable decomposition of 
[Falk and Winther, \textit{Found.
Comput. Math.} (2025), DOI~10.1007/s10208-025-09700-2] and a multilevel 
lowest-order stable decomposition of [Hiptmair, Wu, and Zheng, 
\textit{Numer. Math. Theory Methods Appl.} 5 (2012), 297–332]. 
This allows us to prove that our multigrid solver contracts the 
algebraic error at each iteration $p$-robustly and, simultaneously, 
that the associated a~posteriori 
estimator is $p$-robustly efficient.
Next, we use this multilevel methodology to define a two-level domain 
decomposition method where the subdomains consist of overlapping patches of 
coarse-level elements sharing a common coarse-level vertex. We again establish a 
$p$-robust contraction of the solver and $p$-robust efficiency of the 
a~posteriori 
estimator. Numerical results presented both for the multigrid approach and the 
domain decomposition method confirm the theoretical findings.
\end{abstract}


\tableofcontents

\section{Introduction}

In many physical problems studying fluid flows, the main focus is to obtain an 
accurate representation of the velocity variable. While different discretization 
methods can be used to approximate the fluid velocity, the mixed finite element 
method, see e.g., Boffi, Brezzi, 
and Fortin~\cite{Boffi_Brezzi_Fortin_book} and the 
references therein, has been one of the most attractive approaches because of 
the accuracy, robustness, and instantaneous local mass conservation it provides.
In order to benefit from these advantages, suitable algebraic iterative solvers 
should also be considered. One difficulty is that some common formulations lead 
to a saddle-point form with an indefinite linear system, see e.g., Benzi, Golub, 
and Liesen~\cite{Ben_Gol_Lie_05} or Brenner~\cite{Brenner_MFE_Solvers_09} and 
the references therein.

Amidst a large class of algebraic iterative solvers, multilevel methods such as 
multigrid methods and domain decomposition methods with a coarse-grid 
solve have 
proven to be efficient, accurate, and robust in many different applications.
For a domain decomposition approach in the setting of mixed discretizations, 
we refer to, e.g.,
the works of Glowinski and Wheeler~\cite{Glow_Whe_MFE_DD_88}, 
Cowsar, Mandel, and Wheeler~\cite{Cow_Man_Whee_BDD_MFE_95}, and more recently 
to, e.g., 
Jayadharan, Khattatov, and Yotov~\cite{JayadharanKhattatovYotov21}. 
Another possible approach is that of multigrid methods, see, e.g.,
Brenner~\cite{Bren_MG_MFE_92} for an optimal multigrid in a lowest-order 
Raviart--Thomas setting, Wheeler and Yotov~\cite{Wheeler_Yotov_00} for a 
non-matching grids setting, 
Sch\"oberl and Zulehner~\cite{Schoeberl_Zulehner_03}, 
Takacs and Zulehner~\cite{TakacsZulehner13_MGall}, or Brenner, Oh, and 
Sung~\cite{BrennerOhSung18_MGDarcy} for all-at-once multigrid methods, where 
all the 
unknowns are treated simultaneously and the error analysis is similar to that 
of non-conforming methods since mesh-dependent inner products are used. 
An abstract framework for multigrid convergence in Raviart--Thomas spaces 
is developed in
Arnold, Falk, and Winther~\cite{Arn_Falk_Wint_MG_H_div_H_curl_00}.
In some cases, it is possible to rewrite the problem such that, if one 
first constructs a suitable initial approximation of the velocity which 
satisfies a divergence constraint, then only a symmetric and positive 
definite divergence-free problem remains to be solved. 
This approach was followed in 
Chavent et al.~\cite{Chav_Coh_Jaf_Dup_Rib_MFE2D_84} and then later in
Ewing and Wang~\cite{Ewing_Wang_92,Ewing_Wang_94},
Mathew~\cite{Mathew_93},
Hiptmair and Hoppe~\cite{Hipt_Hopp_MFE3D_99}, or Cai et al.~\cite{Cai_etal} 
to develop multilevel solvers for nested spaces. We adopt the same setting here.

In this work, we present two solvers for mixed finite element discretizations: 
a multigrid solver and a domain decomposition solver. Our analysis is however, 
unified, since our domain decomposition solver fits within our multigrid 
setting. One iteration of our multigrid solver
consists in a V-cycle with zero pre- and solely one post-smoothing step with 
additive Schwarz (block-Jacobi) as a smoother and 
\emph{optimal levelwise step-sizes} given by line search at the error 
correction stage, see e.g., Heinrichs~\cite{Heinrichs_88_line_relaxation_MG}.
This, as in Miraçi, Papež, and Vohralík~\cite{Mir_Pap_Voh_lambda} for 
conforming finite 
elements, leads to the following {\em Pythagorean formula} for the decrease of 
the algebraic error in each step:
\begin{align}\label{intro_formula}
	\big\|\dif^{-1/2}(\u_J - \u_J^{i+1})\big\|^2
	= \big\|\dif^{-1/2}(\u_J - \u_J^i)\big\|^2 
	- \underbrace{\sum_{j=0}^J 
	\big(\lambda_j^i\big\|\dif^{-1/2}\brho_j^i\big\|\big)^2}_{\big(
		\etalg \big)^2},
\end{align}
where $i$ is the solver iteration counter, $\dif$ is the diffusion tensor, 
$j \in \{0, \ldots, J \}$ is the level counter, $\u_J $ is the (unknown) exact 
algebraic solution, $\u_J^i$ denotes the current iterate, $\u_J^{i+1}$ is the 
iterate that is being computed via $\brho_j^i$, the levelwise smoothing 
corrections, and $\lambda_j^i$ are the levelwise optimal step-sizes.
In particular, formula~\eqref{intro_formula} gives a \emph{computable 
levelwise decrease} of the (square of the) algebraic error by the factors
$ \big(\lambda_j^i\big\|\dif^{-1/2}\brho_j^i\big\|\big)^2 $.
This naturally defines a built-in (no additional construction is required) 
a~posteriori estimator $\etalg$, representing a guaranteed lower bound of the 
algebraic error $\big\|\dif^{-1/2}(\u_J - \u_J^{i})\big\|$ on step $i$
\begin{align}
	\etalg 
	\le \big\|\dif^{-1/2}(\u_{J} - \u_{J}^{i})\big\|.
\end{align}
In this 
sense, we refer to the present algebraic solver as \emph{a-posteriori steered}. 
Moreover, one can use the block-smoothing structure of the solver and the 
definitions of $\lambda_j^i$ to further rewrite the estimator on the algebraic 
error as
\begin{align}\label{intro_locformula}
	\big( \etalg\big)^2 
	=  \big\|  \dif^{-1/2} \brho_0^i \big\| ^2 
	+ \sum_{j=1}^J \lambda_j^i \sum_{\ver \in \V_j} 
	\big\|\dif^{-1/2}\brho_{j,\ver}^i\big\|_{\omaj}^2,
\end{align}
which localizes our estimate on the algebraic error with respect to mesh levels 
as well as with respect to vertex patches $\omaj$ on each level.
This multigrid methodology is then used to define a two-level 
\emph{domain decomposition} method, where the subdomains are the overlapping 
patches of the coarsest level, a fine level solve is applied on each of them, 
and a coarse solver ensures uniform contraction. In particular, our methods 
do not need any additional smoothing steps or any damping/relaxation parameters 
which might require tuning. Details form the content of 
Theorem~\ref{thm_error_contr}.

Further main results read as follows. First, we prove that the introduced 
multigrid solver as well as the domain decomposition solver \emph{contract} 
the algebraic error at each iteration, i.e.,
\begin{align}
	\big\|\dif^{-1/2}(\u_{J} - \u_{J}^{i+1})\big\| 
	\le \alpha \big\|\dif^{-1/2}(\u_{J} - \u_{J}^{i})\big\|, 
	\qquad 0 < \alpha < 1,
\end{align}
see Theorem~\ref{thm_converg} for details. Second, we show that the 
associated a~posteriori estimators are \emph{efficient} in that
\begin{align}
	\etalg 
	\ge \beta \big\|\dif^{-1/2}(\u_{J} - \u_{J}^{i})\big\|, 
	\qquad 0 < \beta < 1,
\end{align}
see Theorem~\ref{thm_upper_bound}. 
In fact, proving the two above results is equivalent, owing to the 
connection between solvers and estimators as described by the Pythagorean 
formula~\eqref{intro_formula}. Importantly, the above results 
hold independently of the polynomial degree $p$, i.e., \emph{$p$-robustly}.

A crucial ingredient needed for our analysis is a polynomial-degree-robust 
\emph{multilevel stable splitting}. In two space dimensions, we can adapt 
the one on the discrete stream-function Lagrange spaces from
Miraçi, Papež, and Vohralík~\cite{Mir_Pap_Voh_19}, itself obtained by 
combining the $p$-robust one-level stable splitting achieved in 
Sch\"oberl et al.~\cite{SchMelPechZag_08} and a multilevel piecewise affine 
stable splitting from
Xu, Chen, and Nochetto~\cite{Xu_Chen_Noch_opt_MG_loc_ref_09}. 
In three space dimensions, the recent result by Falk and 
Winther~\cite{Falk_Winther_25} yields a suitable $p$-robust decomposition over 
one level that we can combine with the lowest-order multilevel splitting of 
Hiptmair, Wu, and Zheng~\cite{Hipt_Wu_Zheng_cvg_adpt_MG_12}. Alternatively, 
a not necessarily $p$-robust decomposition follows by the one-level patch-wise 
decomposition of a divergence-free velocity owing to 
Chaumont-Frelet and Vohral{\'i}k~\cite{Chaum_Voh_Maxwell_equil_23}.
We emphasize that the analysis for both the multigrid solver and the domain 
decomposition method rely on the same arguments and that the results hold for 
simplicial meshes that are quasi-uniform as well as possibly highly 
graded, i.e., 
locally refined in specific regions of the domain. 
Indeed, even though we consider initial meshes to be quasi-uniform, the use of 
local refinement via newest-vertex bisection can produce highly graded meshes 
which are included in our analysis. In such a case, the resulting fine mesh, 
see, e.g., Figure~\ref{fig_DD}, is not quasi-uniform, meaning that elements 
which are far away from each other may have very different sizes, but it is 
shape-regular since the minimal angles do not degrade.

This work is organized as follows. In Section~\ref{sec_model_pb} we present 
the model problem and its mixed finite element discretization. 
The multilevel setting and assumptions used in our theory are collected 
in Section~\ref{sec_ML_sett}.
In Section~\ref{sec_mg} we present the a-posteriori-steered multigrid solver 
with its associated a~posteriori estimator of the algebraic error and in 
Section~\ref{sec_AS} we similarly present the domain decomposition method with 
the associated a~posteriori estimator. 
Our main results are summarized in Section~\ref{sec_main_res} and the 
numerical experiments are presented in Section~\ref{sec_num}. 
In Section~\ref{sec_stable_decompositions}, we state the crucial 
$p$-stable multilevel decomposition for divergence-free Raviart--Thomas 
piecewise polynomials. In Sections~\ref{sec_MLD_2D} and~\ref{sec_MLD_3D}, we 
then respectively prove it in two and three space dimensions, passing through 
intermediate results which can be of independent interest. 
Finally, the proof of the central Theorem~\ref{thm_upper_bound} 
forms the content of 
Section~\ref{sec_proofs}. Some concluding remarks are given in 
Section~\ref{sec_conclusions}.

\section{Model elliptic problem and its mixed finite element discretization}
\label{sec_model_pb}
Let $\Omega \subset {\mathbb R}^d$, $d = 2, 3$,
be a polygonal or a polyhedral domain, i.e., open, bounded, connected set, with Lipschitz boundary. We additionally assume that $\Omega$ is simply connected
(no holes for $d=2$, no tunnels for $d=3$). This guarantees exactness of the
de~Rham complex with homogeneous boundary conditions at the level of
$\H_{0}({\div}; \Omega) $: every divergence-free
$\boldsymbol{v}\in \H_{0}({\div}; \Omega) $ is the curl of a potential
with vanishing tangential trace, i.e., $\boldsymbol{v}=\nabla \times \psi$
with $\psi\in H^1_0(\Omega)$ for $d=2$ and
$\boldsymbol{v}=\nabla \times \boldsymbol{\psi}$ with
$\boldsymbol{\psi}\in \H_{0}(\curl; \Omega) $ for $d=3$,
see, e.g., Amrouche, Bernardi, Dauge, and Girault~\cite{Amr_Ber_Dau_Gir_98_potentials3d}
and Arnold~\cite[Section~4.5.5]{Arnold_FEEC_18}. Note the subtlety that, without boundary
conditions, the analogous statement requires instead $\partial\Omega$ connected (no holes for $d=2$, no cavities for $d=3$),
cf. Girault and Raviart~\cite[Corollary 2.4, Theorem 2.9, 
Remark 3.10]{Gir_Rav_NS_86} or Cantarella, DeTurck, and 
Gluck~\cite{Cant_DeTurck_Gluck_calc_top_3D_02}. This exactness holds at the discrete level as well. Hence, under the same assumption, every
divergence-free $\boldsymbol{v}_h\in  {\RT}_p(\mathcal{T}_h)\cap
\H_{0}({\div}; \Omega)$ satisfies $\boldsymbol{v}_h= \nabla \times \psi_h$
with $\psi_h\in\mathcal{P}_{p+1}(\mathcal{T}_h)\cap H^1_0(\Omega)$ for $d=2$,
and $\boldsymbol{v}_h=\nabla \times \boldsymbol{\psi}_h$ with
$\boldsymbol{\psi}_h\in{\N}_p(\mathcal{T}_h)\cap
\H_{0}(\curl; \Omega) $ for $d=3$.

\subsection{Model elliptic problem}

Consider an elliptic partial differential equation in a mixed form modeling, 
for example, the single phase flow in porous media:
\begin{equation}\label{elliptic}
	\u = - \dif\nabla \gamma , 
	\quad \nabla {\cdot} 
	\u = f \quad \text{in } \Omega, 
	\quad \u{\cdot}\n = 0 \quad \text{on } \partial \Omega,
\end{equation}
where $\gamma$ is the fluid pressure, $\u$ is the fluid Darcy velocity, 
$f$ is the source term such that 
$\int_\Omega f(\bx) \mathrm{d}\bx = 0$ in $\Omega$, 
$\n$ is the outward normal vector on $\partial \Omega$, and $\dif$
is a bounded, symmetric, and uniformly positive definite tensor 
representing the material permeability divided by the fluid viscosity. 
More precisely, we assume that there exist 
$0 < \Lambda_{\min} \le \Lambda_{\max}$ such that all eigenvalues of 
$ \dif$ belong to the interval $[\Lambda_{\min} , \Lambda_{\max}]$. 
The homogeneous Neumann boundary condition in~\eqref{elliptic} 
is assumed merely for simplicity of exposition; inhomogeneous and 
mixed boundary conditions can be considered as well.

\subsection{Function spaces and weak formulations}

Let $({\cdot},{\cdot})_S$ and $\|{\cdot}\|_S$, $S \subset {\mathbb R}^d$, be the
$L^2(S)$ inner product and norm, respectively, where we omit the
subscript if $S = \Omega$.
Define the spaces
\begin{subequations}\label{eq_spaces}
	\begin{align}
	\bV & \eq \H_{0}({\div}; \Omega) 
	:=  \{ \bv \in \H ({\div}; \Omega), \bv{\cdot}\n =0 
	\text{ on } \partial \Omega 
	\}, 
	\\
	W   & \eq L^2_{0}(\Omega) 
	:= \{ w \in L^2(\Omega), (w, 1) = 0 \text{ in } \Omega \} .
	\end{align}
\end{subequations}
Here, $\bv {\cdot} \n = 0$ on $\partial \Omega$ means that $(\nabla {\cdot} \bv, \phi) + (\bv, \nabla \phi) = 0$ for all $\phi \in H^1(\Omega)$.
The mixed weak formulation of~\eqref{elliptic}, see e.g., 
Boffi, Brezzi, and Fortin~\cite{Boffi_Brezzi_Fortin_book}, is: 
find $\u \in \bV$ and $\gamma \in W$
such that
\begin{subequations}\label{weak}
	\begin{alignat}{2}
	(\dif^{-1}\u,\bv ) - (\gamma,\nabla {\cdot} \bv) 
	& = 0 \qquad \quad & & \forall \, \bv \in \bV, \label{weak_1} \\
	(\nabla {\cdot} \u,w) & = (f,w) & & \forall \, w \in W. \label{weak_2}
	\end{alignat}
\end{subequations}
The saddle-point structure of~\eqref{weak} may cause (after discretization) severe troubles to iterative algebraic solvers. Let us thus immediately remind that
problem~\eqref{weak} can be written equivalently via the 
dual weak formulation: find $\u \in \bV^f$ such that
\begin{equation}\label{weak_div_free}	
	(\dif^{-1}\u,\bv ) = 0 \quad \forall \, \bv \in \bV^0,
\end{equation}
where, for $g \in L^2_{0}(\Omega)$, 
\begin{equation} \label{eq:Vz}
	\bV^g \eq \{\bv \in \bV: \nabla {\cdot} \bv = g\}.
\end{equation}

\subsection{Mixed finite element discretization} \label{sec_MFE}

In order to discretize the model problem~\eqref{weak}, we use a 
shape-regular mesh $\T_{J}$, partitioning $\Omega$ into $d$-simplices 
(triangles or tetrahedra). This will be the finest mesh from a hierarchy described in Section~\ref{sec_ML_hier} below.
Then, we fix an integer $p \ge 0$ which denotes the polynomial degree used 
in our mixed finite element (MFE) spaces
$\bV_J \times W_J \subset \bV \times W$.
For our setting, we shall work with the Raviart--Thomas (${\RT_p}$) spaces, 
see Raviart and Thomas~\cite{Ra_Tho_MFE_77} for two space dimensions, 
N\'ed\'elec~\cite{Ned_mix_R_3_80} for three space dimensions, or 
Boffi, Brezzi, and Fortin~\cite[Section 2.3.1]{Boffi_Brezzi_Fortin_book}.
We introduce the space $\PP_{p}(K)$ of scalar-valued polynomials of 
degree $p$ on an element $K \in \T_{J}$ and denote by 
${\RT}_{p} (K) := [\PP_{p}(K)]^d + \PP_{p}(K)\bx$ 
the Raviart--Thomas(--N\'ed\'elec) space on $K \in \T_{J}$.
Define 
${\RT}_{p}(\T_{J}) 
:= \{ \bv_{J} \in {\bm{L}}^2(\Omega),  \bv_{J}|_K \in {\RT}_{p} (K) 
\ \forall K \in \T_{J} \}$ 
and $\PP_p(\T_{J}) := \{ w_{J} \in L^2(\Omega), \  w_{J}|_K \in \PP_{p}(K) \ \forall K \in \T_{J} \}$ the broken (elementwise) spaces on the mesh $\T_J$ and
\begin{align}
	\bV_{J} & 
	:= \{ \bv_{J} \in \bV,  \bv_{J}|_K \in {\RT}_{p} (K) 
	\ \forall K \in \T_{J} \} = {\RT}_{p}(\T_{J}) \cap \H_{0}({\div}; \Omega), \label{V_h} \\
	W_{J}   & 
	:= \{ w_{J} \in W, \  w_{J}|_K \in \PP_{p}(K) 
	\ \forall K \in \T_{J}   \} = \PP_p(\T_{J}) \cap L^2_{0}(\Omega) . \label{W_h}
\end{align}
Note that $\bV_{J}$ is a discrete conforming subspace of $\bV = \H_{0}({\div}; \Omega)$. 
In particular, having discrete functions that belong to $\H({\div}; \Omega)$ 
means that the piecewise ${\RT_p}$ functions satisfy normal continuity 
across mesh faces.

We search for $\u_{J} \in \bV_{J}$ and $\gamma_{J} \in W_{J}$ such that
\begin{subequations}\label{mfe}
	\begin{alignat}{2}
	(\dif^{-1}\u_{J},\bv_{J}) - (\gamma_{J},\nabla {\cdot} \bv_{J}) 
	& = 0
	\qquad \qquad & & \forall \, \bv_{J} \in \bV_{J}, \label{mfe_1} 
	\\
	(\nabla {\cdot} \u_{J},w_{J}) & = (f,w_{J}) & & \forall \, w_{J} \in W_{J}. 
	\label{mfe_2}
	\end{alignat}
\end{subequations}
As above, denoting 
$\bV_{J}^g \eq \{\bv_{J} \in \bV_{J}: 
(\nabla {\cdot} \bv_{J},w_{J}) = (g,w_{J})
\ \forall \, w_{J} \in W_{J}\}$, the method \eqref{mfe} can again
be written equivalently as: find $\u_{J} \in \bV_{J}^f$ such that
\begin{equation}\label{mfe_div_free}
	(\dif^{-1}\u_{J},\bv_{J}) = 0 \quad \forall \, \bv_{J} \in \bV_{J}^0.
\end{equation}

\begin{remark}[Other choices of discrete spaces]
Consider the Brezzi--Douglas--Marini spaces, see 
Brezzi, Douglas, and Marini~\cite{BDM_86}, 
and denote the elementwise divergence-free spaces
$\BDM^0_{p} (K) 
:= \{ \bv  \in {\BDM}_{p} (K), \ \nabla {\cdot} \bv=0 \}$,
for all $K \in \T_{J}$. 
Let similarly
$\RT^0_{p} (K) := \{ \bv \in {\RT}_{p} (K), \ \nabla {\cdot} \bv=0 \}$.
From Boffi, Brezzi, and 
Fortin~\cite[Corollary~2.3.1]{Boffi_Brezzi_Fortin_book}, one has 
${\RT}^0_{p} (K) = \BDM^0_{p} (K)$ for all $K \in \T_{J}$. 
Since, after an initial step of constructing a velocity that satisfies 
the divergence constraint, we only work with divergence-free functions, 
the results of this work hold also for this choice of mixed finite element 
spaces.
\end{remark}

\section{Hierarchy of meshes and spaces}\label{sec_ML_sett}
We introduce here the assumptions on the hierarchy of meshes and the 
associated hierarchy of spaces used in this manuscript.

\subsection{Hierarchy of meshes}\label{sec_ML_hier}

To define our solvers, we consider a hierarchy of nested matching simplicial 
meshes of $\Omega$,
$\{\T_j \}_{ 0 \leq j \leq J}$, $J \ge 1$,
where $\T_j$ is a refinement of $\T_{j-1}, \ 1 \leq j \leq J$, and where $\T_J$ 
has been considered in Section~\ref{sec_MFE}. For any element 
$ K \in \T_j $, we denote $h_K :=  {\rm diam} (K) $, and we also use 
$h_j := \underset{K \in \T_j}{\max} \, h_K$.
Let
$\kappa_{\T} >0 $ denote the shape-regularity parameter, i.e., the smallest 
constant such that
\begin{align}
	\max_{K \in \T_j} \cfrac{h_K}{\rho_K}\le \kappa_{\T} \quad
	\forall 0 \leq j \leq J, \label{meshreg}
\end{align}
where $\rho_K$ denotes the diameter of the largest ball inscribed in $K$.

We will work in one of the two settings corresponding to the two assumptions 
below. In the first setting, we assume:

\begin{assumption}[Mesh quasi-uniformity and refinement strength]
\label{assumption_refinement_quasiuniformity} 	
There exists a fixed positive real number $0< C_{\rm qu}  \le  1$ 
such that for all ${ j \in  \{0, \ldots, J\}}$ and for all $K \in \T_{j}$, 
there holds
\begin{align}\label{quasiuniformity}
	C_{\rm qu} h_j \le h_{K} \le h_j.
\end{align}
There further exists a fixed positive real number 
$0< C_{\rm ref} \le 1$ such that for all $ j \in \{1, \ldots, J\}$,
for all $K \in \T_{j-1}$, and all $K^* \in \T_{j}$
such that $K^* \subset K $, there holds
\begin{align}
C_{\rm ref} h_K \le h_{K^*} \le h_K.\label{refinement}
\end{align}
\end{assumption}

In the second setting, we assume that our hierarchy is generated from a 
quasi-uniform coarse mesh $\T_{0}$ by a series of bisections, e.g., 
newest vertex bisection,
see, e.g., Stevenson~\cite{stevenson2008} when $d=2$ and admissible 
coarse mesh, 
Karkulik, Pavlicek, and Praetorius~\cite{kpp2013} 
for general coarse mesh when $d=2$,
and the recent work Diening, Gehring, and Storn~\cite{dgs2023} 
for general coarse mesh when $d \ge 2$. 
We illustrate in Figure
\ref{fig_graded_grids} for $d=2$ how
refining one edge of $\T_{j-1}$, $ j  \in  \{1, \ldots, J\}$, 
leads to a new finer mesh $ \T_j$. 
Let us now denote by $ \V_j$ the set of vertices in $\T_j$ and by 
$\B_j \subset \V_j$ the set consisting of the new vertex obtained 
after the bisection together with its two neighbors on the refinement edge.
Moreover, let $h_{\B_j}$ be the maximal diameter of elements 
having a vertex 
in the set $ \B_j$, for $ j  \in  \{1, \ldots, J\}$. We assume:

\begin{figure}[!ht]
	\phantom{}\hspace{7cm}
	$\T_j$ obtained by bisection of $\T_{j-1}$\\
	\phantom{}\hspace{7cm}
	neighboring vertices ${\bf{b}}_{j_1}$, ${\bf{b}}_{j_3}$ on 
	the refinement edge\\
	\phantom{}\hspace{7cm}
	new vertex after refinement ${\bf{b}}_{j_2}$\\
	\phantom{}\hspace{7cm}
	$ \B_j = \{ {\bf{b}}_{j_1}, {\bf{b}}_{j_2}, {\bf{b}}_{j_3}\}$

	\vspace{-2cm}
	\hspace{2cm}
	\includegraphics[width=0.2\textwidth]{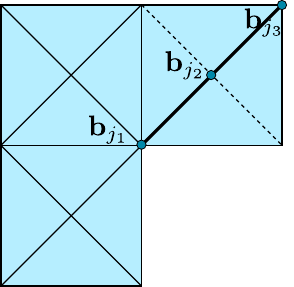}
\caption{
Illustration of the set $\B_j$; the refinement $\T_{j}$ 
(dotted lines) of the mesh $\T_{j-1}$ (full lines).
}
\label{fig_graded_grids}
\end{figure}

\begin{assumption}[Coarsest mesh quasi-uniformity and local refinement 
strength of bisection-generated meshes]\label{graded_grids}
The coarsest mesh $\T_0$ is a quasi-uniform mesh in the sense 
of~\eqref{quasiuniformity}, with parameter $0< C^0_{\rm qu} \le 1$.
The (possibly highly graded) conforming mesh $\T_J$ is generated from 
$\T_0$ by a series of bisections. There exists a fixed positive real number 
$0< C_{\rm loc,qu} \le 1$ such that for all $ j  \in  \{1, \ldots, J\}$, 
there holds
\begin{align}\label{loc_quasiuniformity}
	C_{\rm loc,qu} h_{\B_j}  \le h_K  \le h_{\B_j} 
	\quad \forall K  \in  \T_{j} 
	\text{ such that a vertex of $K$ belongs to $\B_j$}.
\end{align}
\end{assumption}

\subsection{Hierarchy of spaces}\label{sec_ML_hier_spaces}

To define our solvers, we also consider a hierarchy of nested 
mixed finite element spaces associated to the nested meshes $\{ \T_j\}_{j=0}^J$.
First, fix a sequence of non-decreasing polynomial degrees
\begin{equation}\label{eq_pol_degs_j}
0 \le p_0 \le p_1 \le \ldots \le p_J = p.
\end{equation}
Then, for $0 \le j \le J$, define the levelwise mixed finite element space
\begin{align}
	\bV_j & := \{ \bv_j \in \bV, \  \bv_j|_K \in {\RT}_{p_j} (K) 
	\ \forall K \in \T_j \}.
\label{V_j}
\end{align}
We also define, for $0 \le j \le J$,
the levelwise divergence-free discrete spaces
\begin{align}
	\bV^0_j & := \{ \bv_j \in \bV_j, \ \nabla {\cdot} \bv_j =0 \}.
\label{V_j_0}
\end{align}
Any sequence satisfying~\eqref{eq_pol_degs_j} is admissible. In practice, $0 = p_0 = p_1 = \ldots = p_{J-1} \le p_J = p$ gives the smallest-possible spaces $\bV_j$, lowest-order for all levels except the last one.

\section{A-posteriori-steered multigrid solver}\label{sec_mg}

We can now develop a 
multigrid solver for iterative approximation of the discrete 
problem~\eqref{mfe} using the introduced multilevel setting.

\subsection{Setting for patchwise smoothing}\label{sec_ml_solver_setting}

The solver we develop involves solving in each iteration: 
1) a global coarse-grid problem (global residual solve, lowest-order if $p_0=0$); 
2) local problems on patches of elements sharing a vertex on all other 
grids $\T_j$ of the hierarchy (block-Jacobi smoothing, actually explicit smoothing if $p_j=0$ and $d=2$).
We begin with the definition of the patches.  Let $\V_j$ be the set of
vertices of the mesh $\T_j$ and let $ \V_K$ be the set of vertices of
an element $K$ of $\T_j$. Given a vertex $\ver \in \V_j$, 
$ j \in \{ 1, \ldots , J\}$,
we denote the patch associated to $\ver$ by
\begin{align}
	\T_{j}^\ver & :=  \{ K \in \T_j, \ver \in \V_K \}. \label{patch_tau}
\end{align}
Denote the corresponding open patch subdomain by $\omaj 
:= \text{interior} (
\cup_{K \in \mathcal T_j^a} K)$. Define the local 
MFE spaces on $\omaj$ associated with $\T_j$ as
\begin{align}
	\bVaj & := \{ \bv_j \in \bV_j |_{\omaj},
	\ \bv_j {\cdot}\n = 0 \text{ on } \partial \omaj \}. \label{V_ja}
\end{align}
Finally, define
\begin{align}
	\bVajz \eq \{\bvja \in \bVaj,\ \nabla {\cdot} \bvja = 0\}. \label{V_ja_0}
\end{align}

\begin{remark}[Choice of patch subdomains]
Other types of patches can also be considered. For example, in
\emph{\cite{Mir_Pap_Voh_19}}, larger patches, obtained by combining
all elements in the coarser mesh $\T_{j-1}$ that share a vertex in
$\V_{j-1}$, are also studied. The trade-off is that, 
though the local problems are larger in size, there are fewer larger patches. The theoretical results also apply in this case.
For simplicity, we limit the presentation here to the above 
standard vertex patches.
\end{remark}

\subsection{Multigrid solver}\label{sec_ml_solver}

We now proceed with the definition of the iterative solver.

\begin{alg}[A-posteriori-steered multigrid solver]\label{Definition_solver}
~
\begin{enumerate}
\item Initialize $\u_{J}^0 \in \bV_{J}^f$. 
Thus, $\u_{J}^0$ has the requested divergence given by $f$, 
and all the subsequent corrections will be looked for as divergence-free. 
Let $i := 0$ and define a tolerance $\tau > 0$.

\item Perform the following steps (a)--(d): \label{it_2}

\begin{enumerate}
\item Solve the coarse-grid problem: find the global correction 
$\brho_0^i \in \bV_0^0$ as the solution of the global residual problem
\begin{equation}\label{rho_0}
(\dif^{-1}\brho_0^i,\bv_0) =
- (\dif^{-1}\u_{J}^i,\bv_0)  \quad  \forall \,  \bv_0 \in \bV_0^0.
\end{equation}
Define
$\lambda_0^i := 1$ and the coarse-grid update
\begin{equation}\label{u_0_i}
\widetilde \u_0^i := \u_{J}^i + \lambda_0^i\brho_0^i \in \bV_{J}^f.
\end{equation}

\item For higher levels $ 1 \le j \le J$:\\

Compute the local corrections $\brho_{j,\ver}^i  \in \bVajz$ 
as solutions of the patch residual problems,
for all level $j$ vertices $\ver \in \V_j$,
\begin{equation}\label{patch_problem}
	(\dif^{-1}\brho_{j,\ver}^i,\bvja)_{\omaj} 
	= -(\dif^{-1}\widetilde \u_{j-1}^i,\bvja)_{\omaj}
	\quad \forall \bvja \in \bVajz.
\end{equation}

Define the levelwise correction $\brho_j^i \in \bV_j^0$ by
\begin{equation}\label{rho_j}
	\brho_j^i := \sum_{\ver \in \V_j} \brho_{j,\ver}^i.
\end{equation}

If $\brho_j^i \ne 0$, define the 
step size by the line search
\begin{equation}\label{step_size}
	\lambda_j^i 
	:= 
	-\frac{(\dif^{-1}\widetilde \u_{j-1}^i,
	\brho_j^i)}{\big\|\dif^{-1/2}\brho_j^i\big\|^2};
\end{equation}
otherwise set $\lambda_j^i := 1$. Define the level $j$ update
\begin{equation}\label{update_j}
	\widetilde \u_j^i := \widetilde \u_{j-1}^i 
	+ \lambda_j^i \brho_j^i \in \bV_{J}^f.
\end{equation}

\item \label{it_2c} Set the next iterate
$  \u_{J}^{i+1} := \widetilde \u_{J}^{i}$.
Define the a~posteriori estimator of the algebraic error
\begin{equation}\label{def_eta}
	\etalg 
	:=  \left( \sum_{j=0}^J 
	\big(\lambda_j^i \big\|\dif^{-1/2}\brho_j^i\big\|\big)^2 \right)^{1/2} .
\end{equation}

\item If $\etalg \le \tau $, then stop the solver. 
Otherwise set $i := i+1$ and go to step~\ref{it_2}.

\end{enumerate}
\end{enumerate}
\end{alg}

\subsection{Remarks and basic properties}\label{sec_ml_solver_remarks}

\begin{remark}[Initialization and its cost]
The first step of Algorithm~\ref{Definition_solver} is to construct a 
function $\u_{J}^0 \in \bV_{J}$ with a prescribed divergence, $(\nabla {\cdot} \u_{J}^0,w_{J}) = (f,w_{J}) \ \forall \, w_{J} \in W_{J}$. This is a generic question studied in, e.g., Alonso Rodr\'{\i}guez~et al. \cite{AlonsRodr_Cam_San_div_free_18}, see also the references therein. In the multigrid setting, Algorithm~3 of Papež and Vohralík~\cite{PV22} finds a suitable $\u_{J}^0$ and consists in a coarse grid solve and one multigrid-type iteration.
\end{remark}

\begin{remark}[Cost per iteration]
Algorithm~\ref{Definition_solver} is a multigrid V(0,1)-cycle on the hierarchy of the divergence-free spaces $\bV^0_j$ from~\eqref{V_j_0}. This means that there is zero pre-smoothing and merely one post-smoothing step. The coarse-grid problem~\eqref{rho_0} is of the lowest-order if $p_0=0$ in~\eqref{eq_pol_degs_j}. The patch problems~\eqref{patch_problem} are then of size $1 \times 1$ if $d=2$ and $p_j=0$, $1 \le j \le J-1$, since the dimension of the spaces $\bVajz$ from~\eqref{V_ja_0} is then one (they are spanned by the rotated gradient of the hat functions). In contrast, the patch problems~\eqref{patch_problem} grow quickly in size with the polynomial degree $p_J$, which, recall, is always equal to $p$ from~\eqref{eq_pol_degs_j}. Nevertheless, they can always be solved in parallel, since our approach is of additive Schwarz type. A crucial ingredient both in theory and in practice is then the line search~\eqref{step_size}, which, though not parallel, merely consists of two integrations/scalar products.
\end{remark}

\begin{remark}[A-posteriori-steered solver]
Note from step~\ref{it_2c} of Algorithm~\ref{Definition_solver}
that both the new solver iterate $\u_{J}^{i+1}$ and the a~posteriori 
estimator $\etalg$ of the algebraic error are constructed from the 
levelwise algebraic residual liftings $\brho_j^i$ and the step size 
parameters $\lambda_j^i$. In this sense, the solver has a built-in 
algebraic error estimator and that is why we call it \emph{a-posteriori-steered}.
\end{remark}

\begin{remark}[Compact formulas]
The new iterate can be written in the compact form
\begin{equation}\label{compact_update}
	\u_{J}^{i+1} 
	= \u_{J}^i + \sum_{j=0}^J \lambda_j^i \brho_j^i  
	\stackrel{\eqref{rho_j}}= 
	\u_{J}^i + \brho_0^i 
	+\sum_{j=1}^J \lambda_j^i \sum_{\ver \in \V_j} \brho_{j,\ver}^i.
\end{equation}
It is also easy to see that the local updates satisfy for 
$j \in \{1, \ldots, J \}$
\begin{align}
	(\dif^{-1}  \brho_{j,\ver}^i,  \bv_{j,\ver})_{\omaj} 
	=
	& - (\dif^{-1}\u_{J}^i,  \bv_{j,\ver})_{\omaj}
	- \sum_{m=0}^{j-1} \lambda_m^i(\dif^{-1} \brho_{m}^i, \bv_{j,\ver})_{\omaj}    
	\quad
	\forall \, \bvja   \in  \bVajz.       
	\label{rho_j_compact}
\end{align}
\end{remark}

Now, we explain through the following lemma that the choice of the
step-sizes~\eqref{step_size} leads to the best possible decrease of the 
algebraic 
error along the direction given by $\brho_j^i$, as also seen and used in, e.g., 
Heinrichs~\cite{Heinrichs_88_line_relaxation_MG} 
and~\cite{Mir_Pap_Voh_lambda,Mir_Pap_Voh_W}.

\begin{lemma}[Optimal step-sizes]\label{lem_step_size}
For $j \in \{1, \ldots, J \}$, the step size $\lambda_j^i$ defined in 
\eqref{step_size} satisfies
\begin{equation}
	\lambda_j^i = \underset{\lambda \in \mathbb{R}}{\mathrm{argmin}}
	\big\|\dif^{-1/2}(\u_{J} 
    - (\widetilde \u_{j-1}^i + \lambda \brho_j^i))\big\|.
\end{equation}
\end{lemma}
\begin{proof}
The result follows from determining the minimum of the quadratic function
\begin{align*}
	\big\|\dif^{-1/2}(\u_{J} - (\widetilde \u_{j-1}^i  
    + \lambda\brho_j^i))\big\|^2
	& =
	\big\|\dif^{-1/2}(\u_{J} 
        - \widetilde \u_{j-1}^i)\big\|^2  
	-  2\lambda(\dif^{-1}(\u_{J}  
        - \widetilde \u_{j-1}^i),\brho_j^i)\\
	& \qquad \qquad \qquad \qquad  \quad \ \ \ 
	+ \lambda^2\big\|\dif^{-1/2}\brho_j^i\big\|^2 \\
	&
	\stackrel{\eqref{mfe_div_free}}=
	\big\|\dif^{-1/2}(\u_{J}   
        - \widetilde \u_{j-1}^i)\big\|^2  
	+  2\lambda(\dif^{-1} 
        \widetilde \u_{j-1}^i,\brho_j^i)
	\\
	& \qquad \qquad \qquad \qquad  \quad \ \ \
	+ \lambda^2\big\|\dif^{-1/2}\brho_j^i\big\|^2. 
\end{align*}
Here, in the second equality, we have decisively used that the levelwise 
correction/algebraic residual lifting $\brho_j^i$ is conforming and 
divergence-free, i.e., it belongs to the space $\bV_{J}^0$, which allows 
to eliminate the unknown exact solution $\u_{J}$ from the second term in 
the development via~\eqref{mfe_div_free}.
The above expression implies that
%
\begin{equation*}
	\lambda_{\min} 
	=  
	-\frac{(\dif^{-1} 
	\widetilde \u_{j-1}^i,\brho_j^i)}{\big\|\dif^{-1/2}\brho_j^i\big\|^2}.
\end{equation*}
\end{proof}

\begin{lemma}[Norm of the levelwise corrections as sum of norms of the 
local corrections] \label{lem_cor}
For $\brho_j^i$ given by \eqref{patch_problem}--\eqref{rho_j}, 
$j \in \{ 1, \hdots , J \}$, we have
\begin{align}
	\sum_{\ver \in \V_j} \big\|\dif^{-1/2} \brho_{j,\ver}^i\big\|^2_{\omaj}
	\stackrel{\eqref{patch_problem}} =  
	- \sum_{\ver \in \V_j} (\dif^{-1}
    \widetilde \u_{j-1}^i,\brho_{j,\ver}^i)_{\omaj} 
	\stackrel{\eqref{rho_j}} =
	-(\dif^{-1}\widetilde \u_{j-1}^i,\brho_j^i)
	\stackrel{\substack{\eqref{step_size}}}= 
	\lambda_j^i \big\|\dif^{-1/2}\brho_j^i\big\|^2. \label{brhoaj_lam}
\end{align}
\end{lemma}

\section{A-posteriori-steered domain decomposition solver}\label{sec_AS}
\sectionmark{An overlapping additive Schwarz method}

In this section, we present how to adapt the multigrid methodology
developed in Section~\ref{sec_mg} to a domain decomposition setting.

\subsection{Setting for subdomains and a coarse grid}\label{sec_DD_setting}

We consider a hierarchy of \emph{two} nested matching meshes of $\Omega$, 
denoted by
$\T_{H}$ and $\T_{h}$, where the mesh $\T_{h}$ is obtained from $\T_{H}$ 
by a \emph{sequence} of refinements. More precisely, in the setting of 
Section~\ref{sec_ML_sett} with the multilevel mesh hierarchy $\{\T_j \}_{
0 \leq j \leq J}$, we take $J \ge 1$, $\T_{H} \eq \T_0$, and $\T_{h} \eq \T_J$. 
Only the levels $H$ and $h$ will be used in the algorithm for the 
domain decomposition, whereas all the levels $0\le j \le J$ will be used in 
the theoretical analysis. If $\T_H$ is obtained from $\T_h$ by 
a sequence of $J$ uniform refinements, then we obtain $J=\log_2 (H/h)$.

\begin{figure}
	\center
	\includegraphics[width=0.25\textwidth]{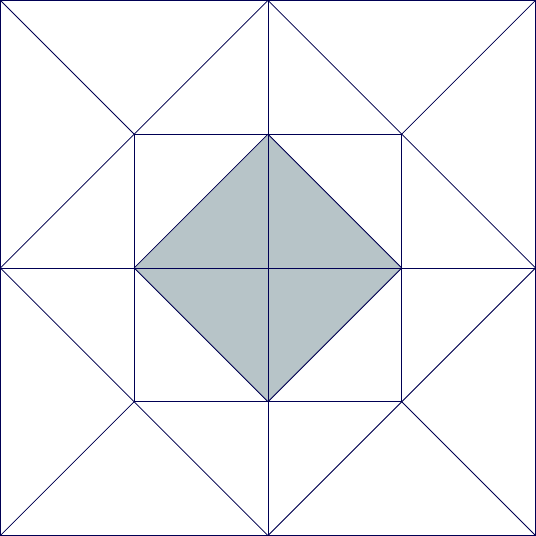}
	\ \ \ \ \
	\includegraphics[width=0.25\textwidth]{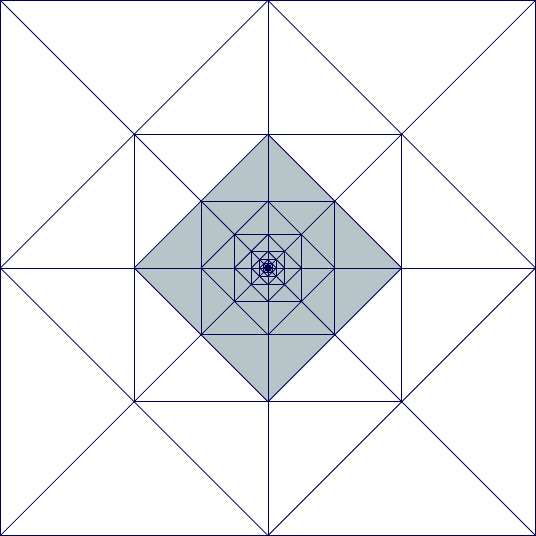}
\caption{Patch in the two-level overlapping additive Schwarz
method: coarse grid $\T_{H}$ (left), fine
grid $\T_{h}$ (right). The highlighted patch consists of four coarse elements 
of $\T_{H}$ which share a vertex and form a subdomain. 
The subdomains (coarse patches) are discretized with
the fine grid~$\T_{h}$.}
\label{fig_DD}
\end{figure}

Denote by $\V_{H}$ the set of vertices of $\T_{H}$ and $\V_K$ the set of 
vertices of an element $K$ of $\T_{H}$.
For each coarse vertex $\ver \in \V_{H}$, the associated patch of 
coarse elements sharing $\ver$ is 
$\T_{H}^\ver \ := \{ K \in \T_{H} , \ver \in \V_K \} $.
Next, we denote the open patch subdomain
corresponding to $\T_{H}^\ver$ by $\omaH$. 
These will be used as \emph{subdomains} 
for the overlapping domain decomposition method. 
Figure~\ref{fig_DD} gives an illustration.
Define the global and local MFE space as
\begin{equation}\label{eq_spaces_H_h}
	\bV_H \eq \bV_0, \quad \bV_h \eq \bV_J, 
	\quad \bV_H^0 \eq \bV^0_0, \quad \bV_h^0 \eq \bV^0_J
\end{equation}
and
\begin{align}
	\Vash & := \{ \bv_{h} \in \bV_{h} |_{\omaH},
	\ \bv_{h} {\cdot}\n = 0 \text{ on } \partial \omaH \}. \label{V_ja_DD}
\end{align}
Remark that the latter spaces are restrictions of the fine-mesh MFE space 
$\bV_{J} $ on the subdomains $\omaH$ with homogeneous Neumann boundary 
conditions on the whole boundary of the subdomain $\omaH$.
Finally, define their divergence-free subspaces as
\begin{align}
	\Vaz \eq \{\va \in \Vash,\ \nabla {\cdot} \va = 0\}. \label{V_ja_0_DD}
\end{align}

\subsection{Domain decomposition solver: overlapping additive Schwarz 
with a coarse-grid solve}
\label{sec_2_level_solver}

Our domain decomposition solver is similar to the multilevel
Algorithm~\ref{Definition_solver} in the case of two levels only. It reads:

\begin{alg}[A-posteriori-steered additive Schwarz domain 
decomposition solver]\label{AS_solver}
~
\begin{enumerate}
\item Initialize $\u_{h}^0 \in \bV_{h}^f$
Let $i := 0$ and define a tolerance 
$\tau > 0$.

\item \label{it_2_DD} Perform the following steps (a)--(d):

\begin{enumerate}
\item Solve the coarse-grid problem: find the global correction 
$\brho_{H}^i \in \bV_{H}^0$ as the solution of the global residual problem
\begin{equation}\label{rho_0_DD}
(\dif^{-1}\brho_H^i,\bv_H) =
- (\dif^{-1}\u_{h}^i,\bv_H)  \quad  \forall \,  \bv_H \in \bV_H^0.
\end{equation}
Define the corse-grid update
\begin{equation}\label{u_0_i_DD}
\u_{H}^i := \u_{h}^i +  \brho_{H}^i \in \bV_{h}^f.
\end{equation}

\item Compute the local corrections $\rhoa \in \Vaz$ as solutions of the 
subdomain residual problems, for all coarse mesh vertices $\ver \in \V_{H}$,
\begin{equation}\label{patch_problem_1}
	(\dif^{-1}\rhoa,\va)_{\omaH} 
	= -(\dif^{-1}\u_{H}^i,\va)_{\omaH}
	\quad \forall \va \in \Vaz.
\end{equation}
Define the overall correction $\brho_{h}^i \in \bV_{h}^0$ by
\begin{equation}\label{rho_h_i_a}
\brho_{h}^i := \sum_{\ver \in \V_{H}} \rhoa.
\end{equation}
If $\brho_{h}^i \ne 0$, define the optimal step size by
\begin{equation}\label{step_size_AS}
	\lambda_{h}^i 
	:= 
	-\frac{(\dif^{-1}\u_{H}^i,\brho_{h}^i)}{\big\|\dif^{-1/2}
	\brho_{h}^i\big\|^2};
\end{equation}
otherwise set $\lambda_{h}^i := 1$.

\item Set the next iterate
\begin{equation}\label{update_h}
	\u_{h}^{i+1} := \u_{H}^{i} + \lambda_{h}^i \brho_{h}^i  \in \bV_{h}^f.
\end{equation}
Define the a~posteriori estimator of the algebraic error
\begin{equation}\label{def_eta_DD}
\etalg := \Big(\big\|\dif^{-1/2}\brho_{H}^i\big\|^2 +
\big(\lambda_{h}^i \big\|\dif^{-1/2}\brho_{h}^i\big\|\big)^2 \Big)^{1/2}.
\end{equation}

\item If $\etalg \le \tau $, then stop the solver. Otherwise set 
$i := i+1$ and go to step~\ref{it_2_DD}.

\end{enumerate}

\end{enumerate}
\end{alg}

\begin{remark}[Initialization and its cost]
The first step of Algorithm~\ref{AS_solver} is again to construct a 
fine-level Raviart--Thomas field with a prescribed divergence, $\u_{h}^0 \in \bV_{h}^f$. 
In the domain decomposition setting, Ewing and 
Wang~\cite[Theorem~3.1]{Ewing_Wang_92} propose an algorithm which consists in a coarse grid solve and of one step of a nonoverlapping domain decomposition iteration. Its more elaborate version is 
Construction~5.1 of Bastidas Olivares~et al. \cite{Bast_Beni_Voh_Yot_DD_MFE_25}. 
\end{remark}

Similarly to Lemma~\ref{lem_cor}, we also have here:

\begin{lemma}[Norm of the overall correction as sum of norms of the 
subdomain corrections]
For $\brho_{h}^i$ given by \eqref{patch_problem_1}--\eqref{rho_h_i_a}, we have
\begin{align}
	\sum_{\ver \in \V_{H}} \big\|\dif^{-1/2} \rhoa \big\|^2_{\omaH}
	\stackrel{\eqref{patch_problem_1}} = 
	- \sum_{\ver \in \V_{H}} (\dif^{-1}\u_{H}^i,\rhoa)_{\omaH} 
	\stackrel{\eqref{rho_h_i_a}}  =
	-(\dif^{-1}\u_{H}^i,\brho_{h}^i)
	\stackrel{\substack{\eqref{step_size_AS}}}= 
	\lambda_{h}^i \big\|\dif^{-1/2}\brho_{h}^i\big\|^2. \label{brhoaj_lam_DD}
\end{align}
\end{lemma}

\section{Main results}\label{sec_main_res}

We now present our main results for the multigrid solver of 
Algorithm~\ref{Definition_solver} and the domain decomposition method of 
Algorithm~\ref{AS_solver}. 
To unify the presentation, we suitably equivalently use the symbols 
$J$ and $h$ and similarly for $0$ and $H$.

\subsection{Error representation and localization on each solver step}

First, we present an important error reduction property of the solvers, 
following as in \cite[Theorem~4.7]{Mir_Pap_Voh_lambda}.
\begin{theorem}[Error representation and localization on each solver step]\label{thm_error_contr} 
There holds
\begin{align}
	\big\|\dif^{-1/2}(\u_{J} & - \u_{J}^{i+1})\big\|^2
	= \big\|\dif^{-1/2}(\u_{J} -  \u_{J}^{i})\big\|^2 -
	(\etalg)^2. \label{error_contr}
\end{align}
Moreover, the above formula can be localized patch-wise and levelwise by 
rewriting the a~posteriori estimator of the algebraic error for 
Algorithm~\ref{Definition_solver} as
\begin{equation}\label{loc_eta}
	\big( \etalg \big)^2 
	= \big\|\dif^{-1/2}\brho_0^i\big\|^2  
	+\sum_{j=1}^J \lambda_j^i \sum_{\ver \in \V_j} 
	\big\|\dif^{-1/2}\brho_{j,\ver}^i\big\|^2_{\omaj},
\end{equation}
and, for Algorithm~\ref{AS_solver}, as
\begin{equation}\label{loc_eta_DD}
	\big( \etalg \big)^2 
	=\big\|\dif^{-1/2}\brho_{H}^i\big\|^2 
	+ \lambda_{h}^i  \sum_{\ver \in \V_{H}} 
	\big\|\dif^{-1/2}\rhoa\big\|^2_{\omaH}.
\end{equation}

\end{theorem}

\begin{proof} We first present the proof in the case of 
Algorithm~\ref{AS_solver}. 
The error representation~\eqref{error_contr} is obtained by using the 
definition of optimal step-sizes and visiting the levels from 
fine to coarse
\begin{align*}
	\big\|\dif^{-1/2}(\u_{J} - \u_{h}^{i+1})\big\|^2
	& \stackrel{\eqref{update_h}}=  
	\big\|\dif^{-1/2}(\u_{J} - \u_{H}^{i})\big\|^2  -
	2 \lambda_{h}^i (\dif^{-1}(\u_{J} - \u_{H}^{i}), \brho_{h}^i)  
	+ \big(\lambda_{h}^i \big\|\dif^{-1/2}\brho_{h}^i\big\|\big)^2
	\\
	& \stackrel{\eqref{mfe_div_free}}=  
	\big\|\dif^{-1/2}(\u_{J} - \u_{H}^{i})\big\|^2  
	+2 \lambda_{h}^i (\dif^{-1}\u_{H}^{i}, \brho_{h}^i)  
	+ \big(\lambda_{h}^i \big\|\dif^{-1/2}\brho_{h}^i\big\|\big)^2
	\\
	& \stackrel{\eqref{step_size_AS}}=  
	\big\|\dif^{-1/2}(\u_{J} - \u_{H}^{i})\big\|^2  
	- \big(\lambda_{h}^i \big\|\dif^{-1/2}\brho_{h}^i\big\|\big)^2
	\\
	& \!\! \! \! \!            
	\stackrel{\substack{\eqref{u_0_i_DD},\eqref{rho_0_DD}\\
	\eqref{mfe_div_free}}}=  
        \! \! \! \!
	\big\|\dif^{-1/2}(\u_{J} - \u_{h}^{i})\big\|^2  
	-\big\|\dif^{-1/2}\brho_{H}^i\big\|^2  
	- \big(\lambda_{h}^i \big\|\dif^{-1/2}\brho_{h}^i\big\|\big)^2
	\\
	& \stackrel{\eqref{def_eta_DD}}=  
	\big\|\dif^{-1/2} (\u_{J} - \u_{h}^{i} )\big\|^2
	- \big( \etalg \big)^2 .
\end{align*}

We proceed similarly for Algorithm~\ref{Definition_solver}. 
Indeed, recall that the levelwise update~\eqref{update_j} gives
\begin{align}\label{update}
{\widetilde \u}_{J-1}^{i} = \u_{J}^{i} 
	+ \sum_{j=0}^{J-1} \lambda_j^i {\brho}_j^i,
    \end{align}
    and the step-size~\eqref{step_size} for the finest level gives
	\begin{align}\label{step}
        \Big(\boldsymbol{K}^{-1} 
		\boldsymbol{\widetilde  u}_{J-1}^{i},\boldsymbol{\rho}_J^i\Big) 
    = - \lambda_J^i \big\|\boldsymbol{K}^{-1/2} \boldsymbol{\rho}_{J}^i\big\|^2.
    \end{align} Thus it holds that
\begin{align*}
	\big\|\dif^{-1/2} & (\u_{J} - \u_{J}^{i+1})\big\|^2
	\stackrel{\eqref{compact_update}}=  
	\Big\|\dif^{-1/2} \big(\u_{J} - (\u_{J}^{i} 
	+ \sum_{j=0}^J \lambda_j^i \brho_j^i)\big)\Big\|^2\\
	& \stackrel{\eqref{mfe_div_free}}=  
	\Big\|\dif^{-1/2} \big(\u_{J} - (\u_{J}^{i} 
	+ \sum_{j=0}^{J-1} \lambda_j^i \brho_j^i)\big)\Big\|^2 
	+ 2 \lambda_J^i \Big(\dif^{-1} (\u_{J}^{i} 
	+ \sum_{j=0}^{J-1} \lambda_j^i \brho_j^i),\brho_J^i\Big)
	+ \big(\lambda_J^i \big\|\dif^{-1/2}\brho_{J}^i\big\|\big)^2\\
	& \stackrel{\eqref{update}}=  
	\Big\|\dif^{-1/2} \big(\u_{J} - (\u_{J}^{i} 
	+ \sum_{j=0}^{J-1} \lambda_j^i \brho_j^i)\big)\Big\|^2
	+  2 \lambda_J^i \big(\dif^{-1} {\widetilde \u}_{J-1}^{i},{\brho}_J^i\big)
	+ \big(\lambda_J^i \big\|\dif^{-1/2} {\brho}_{J}^i\big\|\big)^2\\
    & \stackrel{\eqref{step}}=  
	\Big\|\dif^{-1/2} \big(\u_{J} - (\u_{J}^{i} 
	+ \sum_{j=0}^{J-1} \lambda_j^i \brho_j^i)\big)\Big\|^2
	-  \big(\lambda_J^i\big\|\dif^{-1/2}\brho_{J}^i\big\|\big)^2\\
    & \,\, = \ldots 
	\\
	& \stackrel{\phantom{\eqref{step_size}}}=  
	\big\|\dif^{-1/2} (\u_{J} - \u_{J}^{i} )\big\|^2
	-  \sum_{j=0}^{J}  \big(\lambda_{j}^i\big\|\dif^{-1/2}\brho_j^i\big\|\big)^2
	\stackrel{\eqref{def_eta}}=  
	\big\|\dif^{-1/2} (\u_{J} - \u_{J}^{i} )\big\|^2
	- \big( \etalg \big)^2,
\end{align*}
where we have proceeded similarly as for the finest level $J$ on the coarser levels.
The localized writing of the a~posteriori estimator~\eqref{loc_eta} 
(resp.~\eqref{loc_eta_DD}) then follows then from its definition~\eqref{def_eta} 
(resp.~\eqref{def_eta_DD}) and patch-localization~\eqref{brhoaj_lam} 
(resp.~\eqref{brhoaj_lam_DD}).
\end{proof}

\subsection{Reliability and efficiency of the estimate on the algebraic error}

For the a~posteriori estimators we introduced, there holds:

\begin{theorem}[$p$-robust reliability and efficiency of the algebraic 
	error estimators]\label{thm_upper_bound}
Let either Assumption~\ref{assumption_refinement_quasiuniformity} 
or Assumption~\ref{graded_grids} hold. Let $\u_{J}^{i} \in \bV_{J}^f$ 
be arbitrary. Let $\etalg$ be constructed from $\u_{J}^{i}$ by~\eqref{def_eta} 
in Algorithm~\ref{Definition_solver} or by~\eqref{def_eta_DD} 
in Algorithm~\ref{AS_solver}.
Then, there holds
\begin{equation}\label{lower_bound}
	\big\|\dif^{-1/2}(\u_{J} - \u_{J}^{i})\big\| \ge \etalg
\end{equation}
and
\begin{equation}\label{upper_bound}
	\beta \big\|\dif^{-1/2}(\u_{J} - \u_{J}^{i})\big\| \le \etalg,
\end{equation}
where $0 < \beta \leq 1$ only depends on: 
\begin{itemize}
    \setlength\itemsep{0ex}
    \item the mesh shape regularity parameter $\kappa_{\T}$, 
    \item  the parameters $C_{\rm qu}$ and $C_{\rm ref}$ when 
Assumption~\ref{assumption_refinement_quasiuniformity} is satisfied or the 
parameters $C^0_{\rm qu}$ and $ C_{\rm loc,qu}$ when 
Assumption~\ref{graded_grids} is satisfied,
    \item the diffusion inhomogeneity or anisotropy ratio $\Lambda_{\max}/\Lambda_{\min}$, 
    \item the domain $\Omega$ if $d=3$, and 
    \item at most linearly on the number of mesh levels $J$ (i.e. on $\log_2 (H/h)$ if $\T_H$ is obtained from $\T_h$ by a sequence of 
$J$ uniform refinements).
\end{itemize}
\end{theorem}
The proof is given in Section~\ref{sec_proofs}. 
Note that the fact that the fully computable estimator $ \etalg$ is a 
\emph{guaranteed lower bound} on the unknown algebraic error, 
cf.~\eqref{lower_bound}, is an immediate consequence of~\eqref{error_contr}, 
as the left-hand-side term in~\eqref{error_contr} corresponding to the new algebraic error is nonnegative.

\subsection{Error contraction on each solver step}

Finally, for the two solvers, there holds:

\begin{theorem}[$p$-robust error contraction]\label{thm_converg}
Let either Assumption~\ref{assumption_refinement_quasiuniformity} or 
Assumption~\ref{graded_grids} hold. Let $\u_{J}^{i} \in \bV_{J}^f$ be arbitrary. 
Let $\u_{J}^{i+1} \in \bV_{J}^f$ be constructed from $\u_{J}^{i}$ on 
step~\ref{it_2} of Algorithm~\ref{Definition_solver} or on step~\ref{it_2_DD} of 
Algorithm~\ref{AS_solver}. Then, there holds
\begin{equation}\label{converg}
	\big\|\dif^{-1/2}(\u_{J} - \u_{J}^{i+1})\big\| 
	\le \alpha \big\|\dif^{-1/2}(\u_{J} - \u_{J}^{i})\big\|,
\end{equation}
where $0 \leq \alpha < 1$ is given by $\alpha = \sqrt{1- \beta^2}$
with $\beta$ the constant from \eqref{upper_bound}.
\end{theorem}

\begin{proof}
The (immediate) proof follows from the crucial error representation formula~\eqref{error_contr}, implying equivalence of~\eqref{upper_bound} 
and~\eqref{converg}, as in~\cite[Corollary~6.7]{Mir_Pap_Voh_lambda}.
We present it here for completeness. Starting from \eqref{converg}
with $0 \leq \alpha < 1$,
\begin{align*}
	&\big\|\dif^{-1/2}(\u_{J} - \u_{J}^{i+1})\big\|^2 
	\le \alpha^2 \big\|\dif^{-1/2}(\u_{J} - \u_{J}^{i})\big\|^2 \\
	& \stackrel{\eqref{error_contr}}{\Longleftrightarrow}
	\big\|\dif^{-1/2}(\u_{J} - \u_{J}^{i})\big\|^2 - (\etalg)^2
	\le \alpha^2 \big\|\dif^{-1/2}(\u_{J} - \u_{J}^{i})\big\|^2 \\
	& \stackrel{\phantom{\eqref{error_contr}}}{\Longleftrightarrow}
	(1 - \alpha^2)\big\|\dif^{-1/2}(\u_{J} - \u_{J}^{i})\big\|^2 \le (\etalg)^2.
\end{align*}
\end{proof}

\section{Numerical experiments}\label{sec_num}

In this section, we present numerical experiments for the multigrid solver of 
Algorithm~\ref{Definition_solver} and the domain decomposition method of 
Algorithm~\ref{AS_solver}.
Though the implementation aspects are not the main focus of the work, 
we point out that several options are possible. One can possibly use 
divergence-free basis functions; for some early contributions on this subject, 
see, e.g., Thomasset~\cite{Thomasset81_book}, Hecht~\cite{Hecht84}, or 
Scheichl~\cite{Scheichl_PhD_00,Schei_MFE_dec_03}. 
This approach may, however, be involved, especially in three space dimensions 
and for higher polynomial degrees $p$. We follow here another approach, 
which consists in using basis functions of the entire spaces $\bV_{J}$ 
from~\eqref{V_h}; these are typically available in finite element 
software packages. Then, to solve the problems~\eqref{rho_0} 
and~\eqref{patch_problem} (or, similarly, \eqref{rho_0_DD} 
and~\eqref{patch_problem_1}), one uses the fact that they are equivalent to 
coarse-grid/local saddle-point problems just as~\eqref{mfe_div_free} is 
equivalent to~\eqref{mfe}. Then, no divergence-free basis functions are needed. 
A special attention for the choice of the basis functions and quadrature formulas needs to be paid when high polynomial degrees $p$ are used, 
cf., e.g., Kirby~\cite{Kirby_FE_de_Rham_14}, 
Brubeck and Farrell~\cite{Brubeck_Farrell_hoFEM_21}, and the references therein.

The experiments are designed to highlight our crucial result stating the 
$p$-robustness (proven in Theorems~\ref{thm_upper_bound} and~\ref{thm_converg}). 
For this reason, we stop the iterations when the 
a~posteriori algebraic error estimator $\etalg$ is reduced by the 
factor $10^{5}$; $p$-robustness is indicated by the number of iterations 
independent of the polynomial degree $p$.
The experiments are presented for meshes arising from both uniform and 
local adaptive refinement. For the latter, we employ newest vertex bisection, 
see e.g. Mitchell~\cite{Mitchell_91}, Traxler~\cite{trax_97}, or 
Stevenson~\cite{stevenson2008}, using D\"orfler's bulk-chasing criterion, 
cf. D\"orfler~\cite{Dorfler_marking_96}, with marking parameter 
$\theta_{\textrm{mark}}$ (specified below) and the true discretization error 
(instead of a discretization error estimator) for the sake of reproducibility. 
This criterion reads as: find a set of marked elements 
$\mathcal{M}_J \subseteq \mathcal{T}_J$ of minimal cardinality that satisfies
\begin{align*}
	\theta_{\textrm{mark}} \,\big\|\dif^{-1/2}(\u- \u_J )\big\|^2 \le
	\sum_{K \in \mathcal{M}_J} \big\|\dif^{-1/2}(\u- \u_J) \big\|^2_K.
\end{align*}
Thus, in the local adaptive case, the obtained mesh hierarchies can be 
highly graded. Moreover, since we aim to compare the performance of the 
solver with respect to increasing polynomial degree without other parameters 
changing, the mesh hierarchy is pre-computed once for each problem using 
$p=1$ and re-used in all experiments.

\begin{figure}[ht]
	\center
	\includegraphics[width=0.3\textwidth]{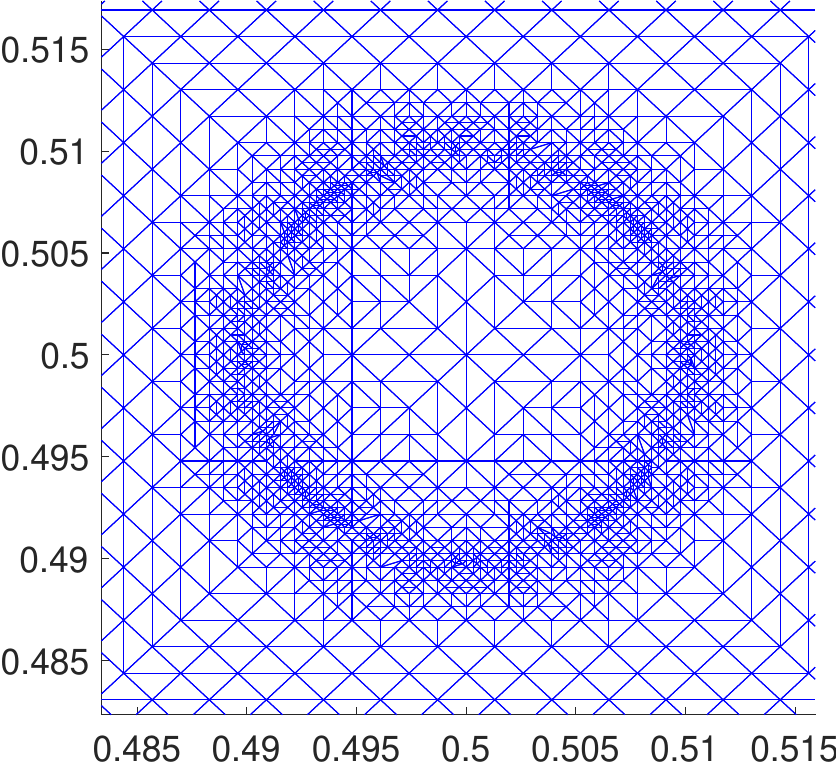}
	\includegraphics[width=0.32\textwidth]{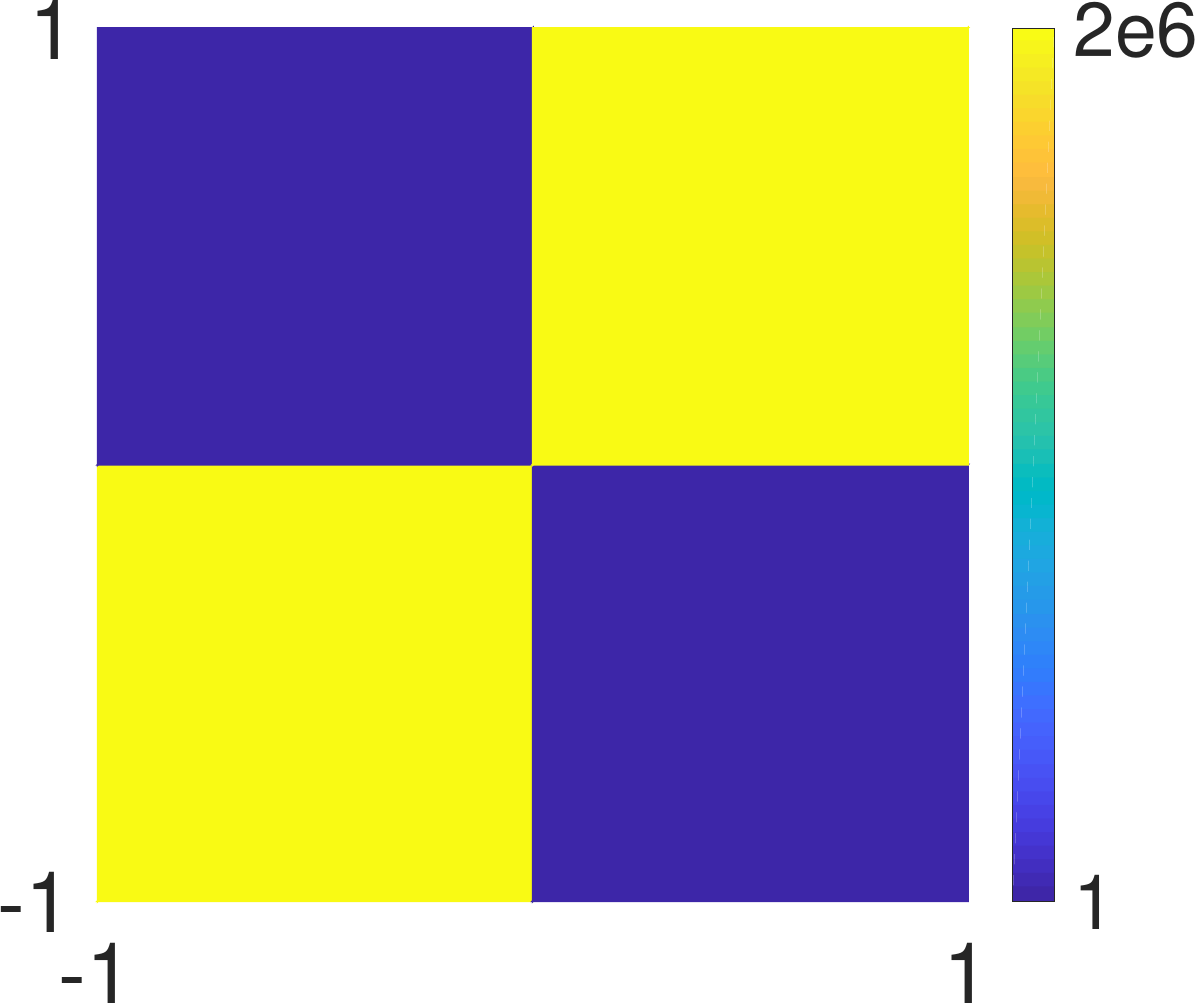}
	\includegraphics[width=0.28\textwidth]{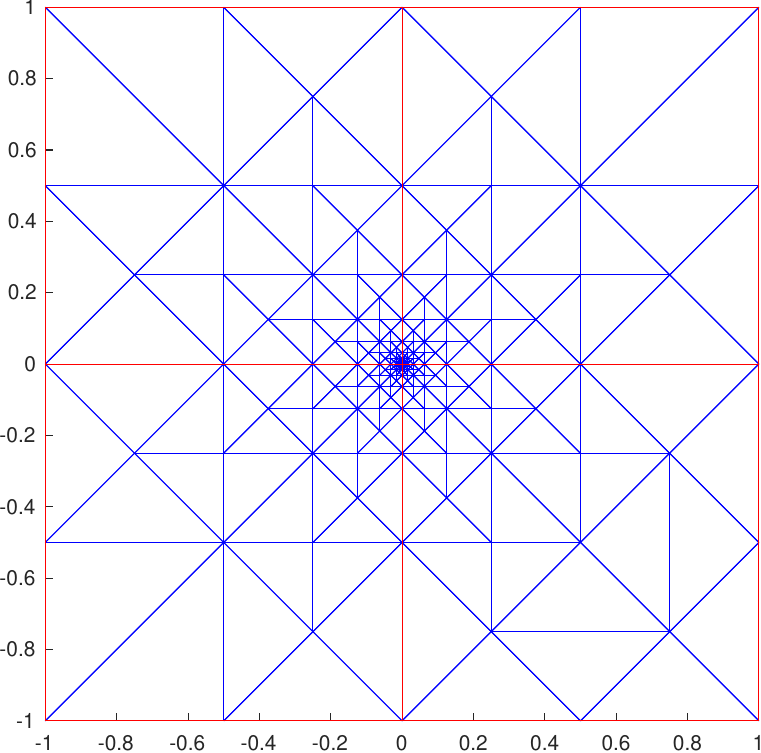}
\caption{Left: mesh (zoomed in) used for the well wavefront test 
case~\eqref{Wellwavefront}, obtained after $J=10$ local adaptive refinements and 
consisting of 9814 elements. Center: variations of the coefficient $c(x, y)$ 
for the piecewise constant diffusion tensor ${\dif = c(x,y)I}$ across the 
domain for the checkerboard test case~\eqref{Checkerboard}. Right: mesh used for the checkerboard test case~\eqref{Checkerboard}, obtained 
after $J=20$ local adaptive refinements and consisting of 1038 elements.}
\label{fig_meshes}
\end{figure}

We consider the following test cases for $d=2$ (we always take 
$\u = - \dif \nabla  \gamma$ and compute $f$ correspondingly):
\begin{itemize}
\item \textbf{Smooth test case and uniform mesh refinement.} 
This test is taken from Ewing and Wang~\cite{Ewing_Wang_94}:
\begin{align}\label{Smooth}
	\text{
	$\gamma (x,y) = {\rm cos}( \pi x ){\rm cos}( \pi y ),
	\  \Omega  =  (0,1)^2$.}
\end{align}
The initial (coarse) mesh size is taken to be $0.5$ and $\dif = I $. 
The refinement is uniform: each existing element is split into four congruent 
ones by joining the midpoints of its edges.
\item 	\textbf{Well wavefront test case and adaptive mesh refinement.}
This test case is taken from Mitchell~\cite{Mitchell_10}:
\begin{align}\label{Wellwavefront}
	\text{
	$\gamma(r) = {\rm tan}^{-1} (\alpha ( r - r_0)), \
	\Omega  =  (0,1)^2$,}
\end{align}
where $r = \sqrt{(x - x_c)^2 + (y - y_c)^2} $ and the parameters are 
$\alpha = \num{1000}$, $x_c = 0.5$, $y_c = 0.5$, $r_0 = 0.01$.
The initial (coarse) mesh size is taken to be 0.5, $\dif = I $, and the bulk-chasing 
parameter for the adaptive refinement is chosen to be 
$\theta_{\textrm{mark}} =0.7$. 
Figure~\ref{fig_meshes} (left) showcases the mesh obtained after 
$J=10$ refinements, zoomed in the vicinity of the ``well''.
\item \textbf{Checkerboard diffusion test case and adaptive mesh refinement.}
Consider a singular exact solution written in polar coordinates as
\begin{align}\label{Checkerboard}
	\text{
	$\gamma(r, \varphi ) = r^{\gamma} \mu (\varphi), \
	\Omega  =  (-1,1)^2$,}
\end{align}
where we follow Kellogg~\cite{Kellogg_75} to define $\mu (\varphi)$ and have a 
regularity parameter $\gamma = \num{0.0009}$ and diffusion contrast 
\num{2001405.429972}. We illustrate the diffusion coefficient across 
the domain in Figure~\ref{fig_meshes} (center). This problem has non-zero 
Neumann boundary conditions, and the implementation is done by re-adjusting 
the right-hand side first. The initial mesh size is taken to be 1 and the 
bulk-chasing parameter for the adaptive refinement is
chosen to be $\theta_{\textrm{mark}} =0.3$; 
see Figure~\ref{fig_meshes} (right) for the mesh 
in the graded hierarchy obtained after $J=20$ refinements. 
\end{itemize}

\begin{table}[tbhp]
\caption{Summary of problem sizes in mesh elements and number of degrees of 
freedom for all the test cases, multigrid and domain decomposition solvers, 
different polynomial degrees $p=1,6$, and mesh hierarchy of $J$ levels.}
\label{Tab_summary_data}
\resizebox{\textwidth}{!}{
\begin{tabular}{cc||cc|cc|cc}
&
&\multicolumn{2}{c|}{\multirow{2}{*}{\begin{tabular}[c]{@{}c@{}}
Smooth\\  \ \ \ \ \ \ \  $J=5$ \  \ \ \ \ \end{tabular}}}
&\multicolumn{2}{c|}{\multirow{2}{*}{\begin{tabular}[c]{@{}c@{}}
Wellwavefront\\  \ \ \ \ \ \ \  $J=12$ \ \ \ \ \ \end{tabular}}}
&\multicolumn{2}{c}{\multirow{2}{*}{\begin{tabular}[c]{@{}c@{}}
Checkerboard\\  \ \ \ \ \ \ \ $J=28$ \ \ \ \ \ \end{tabular}}}
\\
&
\multirow{2}{*}{}
\multirow{3}{*}{}
&\multicolumn{2}{c|}{\multirow{2}{*}{\begin{tabular}[c]{@{}c@{}} \\
$p=1$   \  \ \ \ \ $p=6$  \end{tabular}}}
&\multicolumn{2}{c|}{\multirow{2}{*}{\begin{tabular}[c]{@{}c@{}} \\
$p=1$    \  \ \ \ \  $p=6$  \end{tabular}}}
&\multicolumn{2}{c}{\multirow{2}{*}{\begin{tabular}[c]{@{}c@{}} \\
$p=1$  \  \ \ \ \  $p=6$ \end{tabular}}}
\\
&
&\multicolumn{2}{c|}{}
&\multicolumn{2}{c|}{}
&\multicolumn{2}{c}{}  \\
\hline
\multicolumn{1}{c}{}
&
\multicolumn{1}{c||}{\multirow{4}{*}{\begin{tabular}[l]{@{}l@{}}
\#DoF (mixed) \phantom{min/max     } \\ \#DoF (div-free)
\\ \#$\mathcal{T}_J$ \\ \#$\mathcal{T}_0$\end{tabular}}}
& \num{130816} & \num{1318016} & \num{191486} & \num{1927051} 
& \num{33195} & \num{334215}\\
&& \num{82176} &  \num{861056} & \num{119734} & \num{1256969} 
& \num{20794} & \num{218134}\\
&& \num{16384} &   \num{16384} &  \num{23940} & \num{23940}
&  \num{4153} &   \num{4153}\\
&&    16 &      16 &     16 & 16&    24 &     24\\
\hline
\multicolumn{1}{c}{\multirow{3}{*}{\begin{tabular}[c]{@{}c@{}}
MG\\ (per patch) \end{tabular}}}
&
\multicolumn{1}{c||}{\multirow{3}{*}{\begin{tabular}[l]{@{}l@{}}
min/max \#elements \\ min/max \#DoF (mixed) \\
min/max \#DoF (div-free) \end{tabular}}}
&   2/6  & 2/6     &  2/10  &   2/10  &  2/9  &   2/9 \\
&& 11/41 & 147/462 & 11/69  & 147/770 & 11/62 & 147/693 \\
&&  6/24 & 91/294  &  6/40  &  91/490 &  6/36 &  91/441  \\
\hline
\multicolumn{1}{c}{\multirow{3}{*}{\begin{tabular}[c]{@{}c@{}}
DD\\ (per subdomain) \end{tabular}}}
&
\multicolumn{1}{c||}{\multirow{3}{*}{\begin{tabular}[cl]{@{}l@{}}
min/max \#elements \\ min/max \#DoF (mixed) \\
min/max \#DoF (div-free) \end{tabular}}}
&   \num{2048}/\num{6144}   &   \num{2048}/\num{6144}   &   
	\num{97}/\num{22785}    &     \num{97}/\num{22785}  &  
	\num{20}/\num{3898}     &     \num{20}/\num{3898}  \\
&&  \num{16256}/\num{48960} & \num{164416}/\num{493920} &  
	\num{749}/\num{182203}  &   \num{7714}/\num{1833923}& 
	\num{146}/\num{31156}   &   \num{1561}/\num{313691} \\
&& 	\num{10112}/\num{30528} & \num{107072}/\num{321888} &  
	\num{458}/\num{113848}  &   \num{4998}/\num{1195943}&  
	\num{86}/\num{19462}    &   \num{1001}/\num{204547}  \\
\end{tabular}
}
\end{table}

As for the choice of the polynomial degrees per level, 
recalling~\eqref{eq_pol_degs_j}, any non-decreasing sequence is authorized, 
allowing in particular for the lowest degree $p=0$ in the coarse solve, 
which was tested variedly in~\cite{Mir_Pap_Voh_19}. Here, we rather opt 
for the same polynomial degree $p$ on every level, including the coarsest 
one, i.e., $p_j = p$ for $0 \le j \le J$. The different mesh settings and 
numbers of unknowns are summarized in Table~\ref{Tab_summary_data}.

In the smooth and well wavefront test cases, using the highest polynomial 
degree and highest number of mesh refinements leads to saddle-point algebraic 
systems from~\eqref{mfe} with $\sim 2 {\cdot} 10^6$ 
degrees of freedom. In the checkerboard test case, due to the point singularity, 
meshes are aggressively refined towards the origin. In this case, a number of 
$J=28$ refinements yields a finest mesh of overall \num{4153} triangles,
but the ratio of the largest to smallest triangle is of order $10^{10}$. 
After this point, numerical computations on further refined meshes start to
become numerically unstable in double precision if no specific treatment is considered.

\subsection{Efficiency of the a~posteriori error estimator of the 
algebraic error}

\begin{figure}[!ht]
\center
	\begin{tikzpicture}
	\node (img1)  {\includegraphics[scale=0.125]{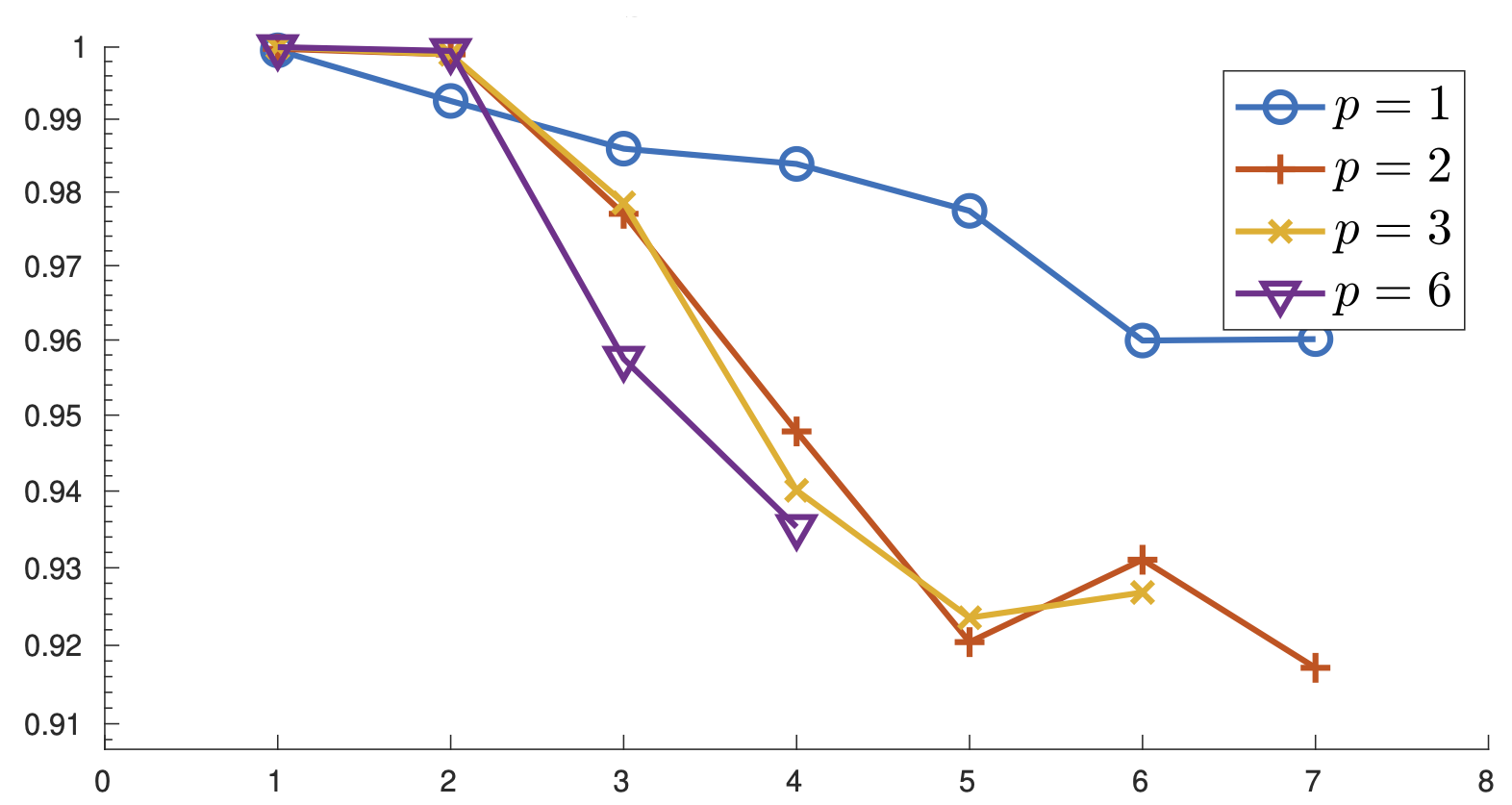}};
	\node[above=of img1, node distance=0cm,yshift=-1cm,font=\color{black}] 
	{MG solver, smooth solution, $J=5$};
	\node[below=of img1, node distance=0cm,yshift=1.2cm,font=\color{black}] 
	{\scriptsize iteration};
	\node[left=of img1, node distance=0cm, rotate=90, 
	anchor=center,xshift=0.2cm,yshift=-0.7cm,font=\color{black}] 
	{\scriptsize effectivity index};
	\end{tikzpicture}
\hfil
	\begin{tikzpicture}
	\node (img1)  {\includegraphics[scale=0.125]{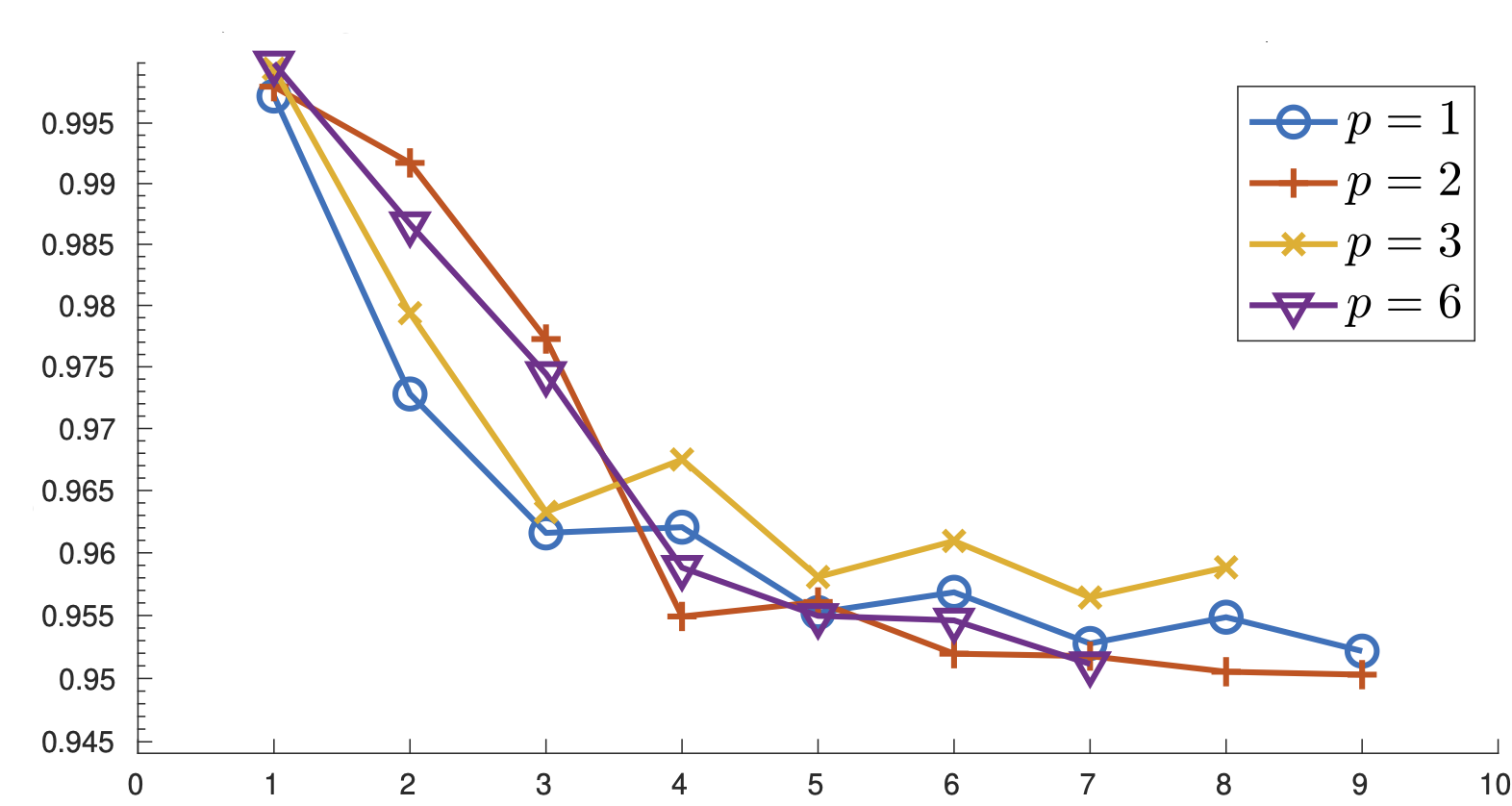}};
	\node[above=of img1, node distance=0cm, yshift=-1cm,font=\color{black}] 
	{DD solver, smooth solution, $J=5$};
	\node[below=of img1, node distance=0cm,yshift=1.2cm,font=\color{black}] 
	{\scriptsize iteration};
	\node[left=of img1, node distance=0cm, rotate=90, anchor=center, 
	xshift=0.2cm,yshift=-0.7cm,font=\color{black}] 
	{\scriptsize effectivity index};
	\end{tikzpicture}
\\
	\begin{tikzpicture}
	\node (img1)  {\includegraphics[scale=0.12]{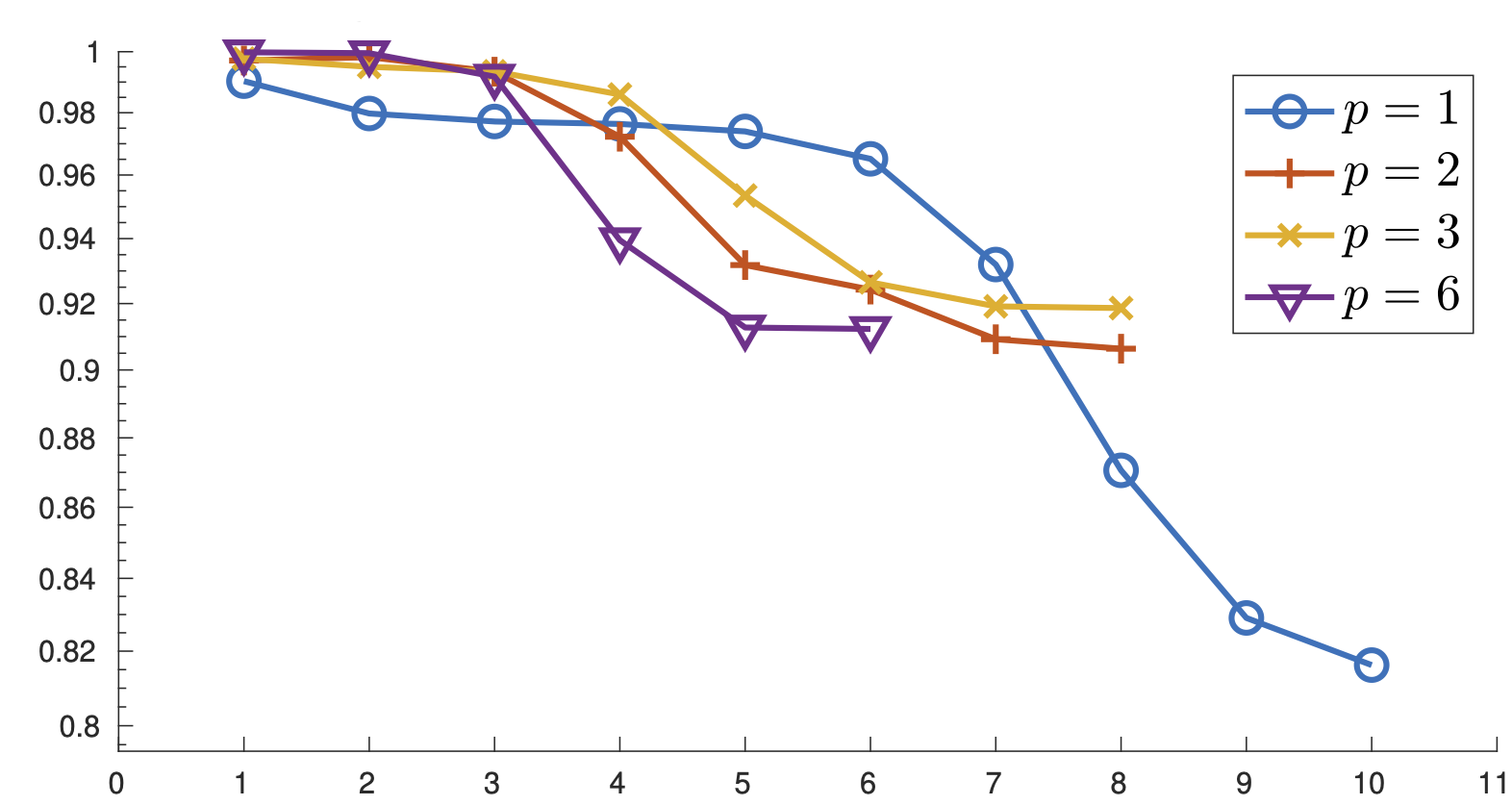}};
	\node[above=of img1, node distance=0cm,yshift=-1cm,font=\color{black}] 
	{MG solver, well wavefront solution, $J=12$};
	\node[below=of img1, node distance=0cm,yshift=1.2cm,font=\color{black}] 
	{\scriptsize iteration};
	\node[left=of img1, node distance=0cm, rotate=90, anchor=center, 
	xshift=0.2cm,yshift=-0.7cm,font=\color{black}] 
	{\scriptsize effectivity index};
	\end{tikzpicture}
\hfil
	\begin{tikzpicture}
	\node (img1)  {\includegraphics[scale=0.12]{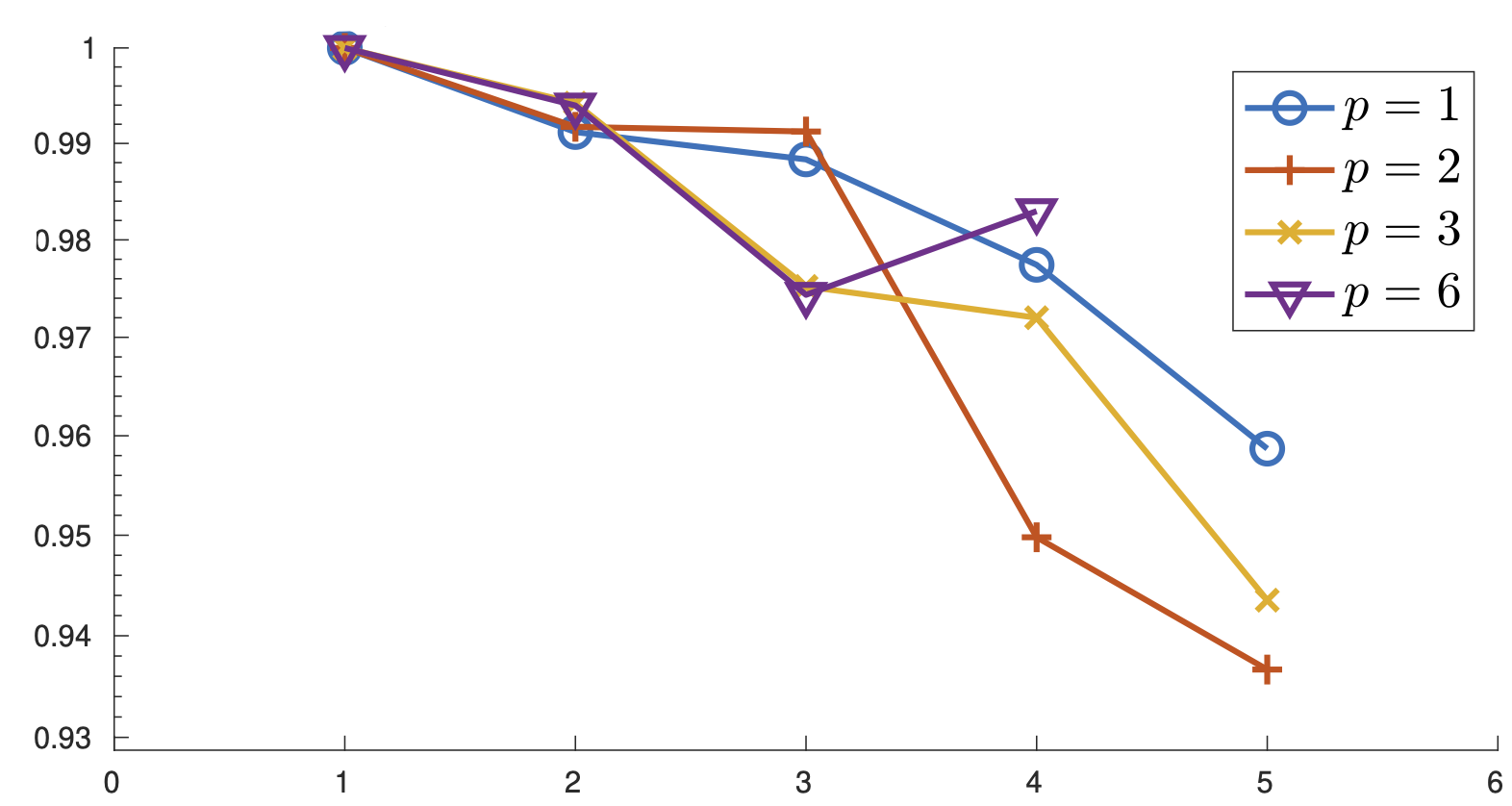}};
	\node[above=of img1, node distance=0cm, yshift=-1cm,font=\color{black}] 
	{DD solver, well wavefront solution, $J=12$};
	\node[below=of img1, node distance=0cm,yshift=1.2cm,font=\color{black}] 
	{\scriptsize iteration};
	\node[left=of img1, node distance=0cm, rotate=90, anchor=center, 
	xshift=0.2cm,yshift=-0.7cm,font=\color{black}] 
	{\scriptsize effectivity index};
	\end{tikzpicture}
\\
	\begin{tikzpicture}
	\node (img1)  {\includegraphics[scale=0.12]{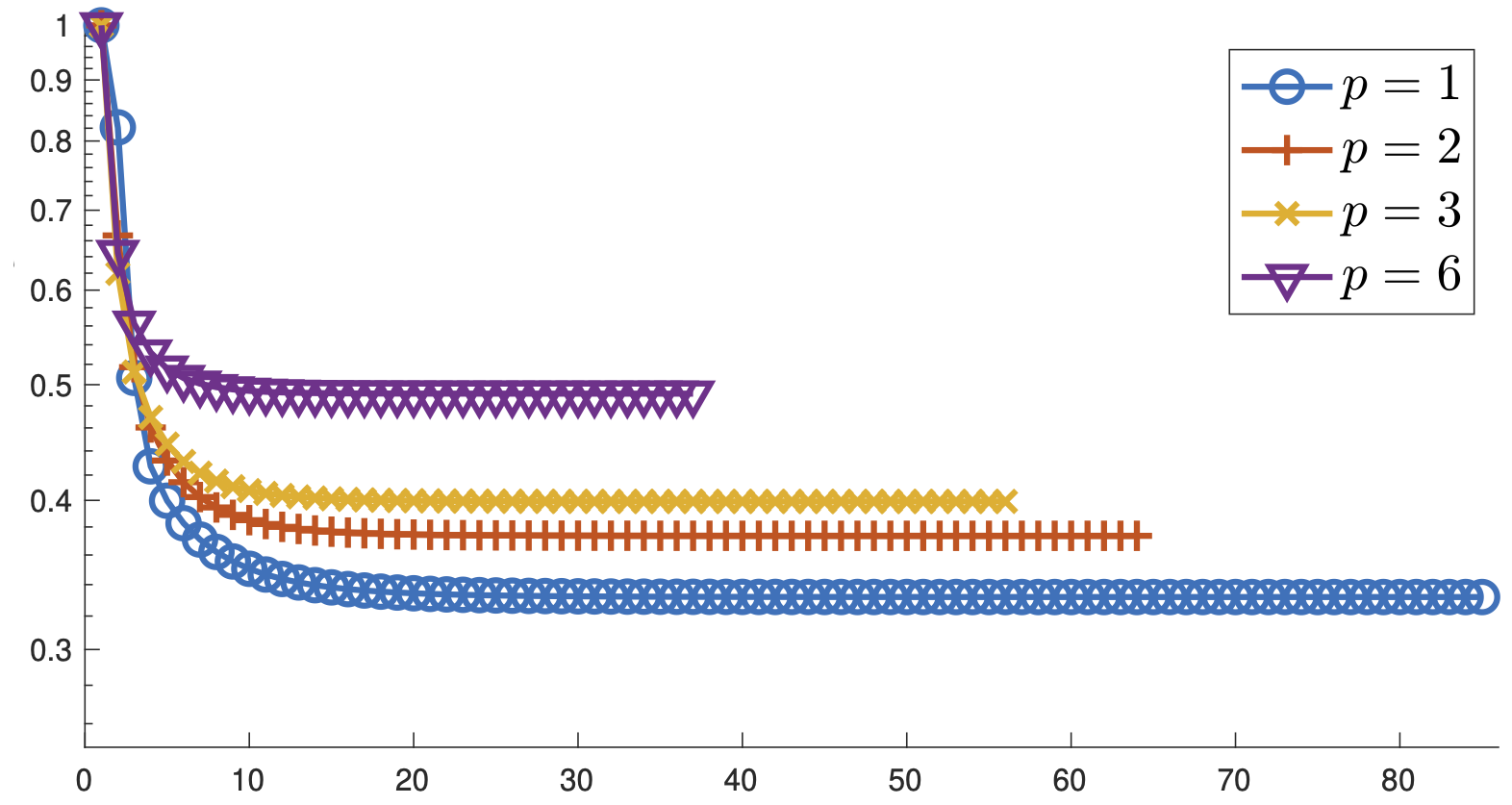}};
	\node[above=of img1, node distance=0cm,yshift=-1cm,font=\color{black}] 
	{MG solver, checkerboard solution, $J=28$};
	\node[below=of img1, node distance=0cm,yshift=1.2cm,font=\color{black}] 
	{\scriptsize iteration};
	\node[left=of img1, node distance=0cm, rotate=90, anchor=center, 
	xshift=0.2cm,yshift=-0.7cm,font=\color{black}]  
	{\scriptsize effectivity index};
	\end{tikzpicture}
\hfil
	\begin{tikzpicture}
	\node (img1)  {\includegraphics[scale=0.12]{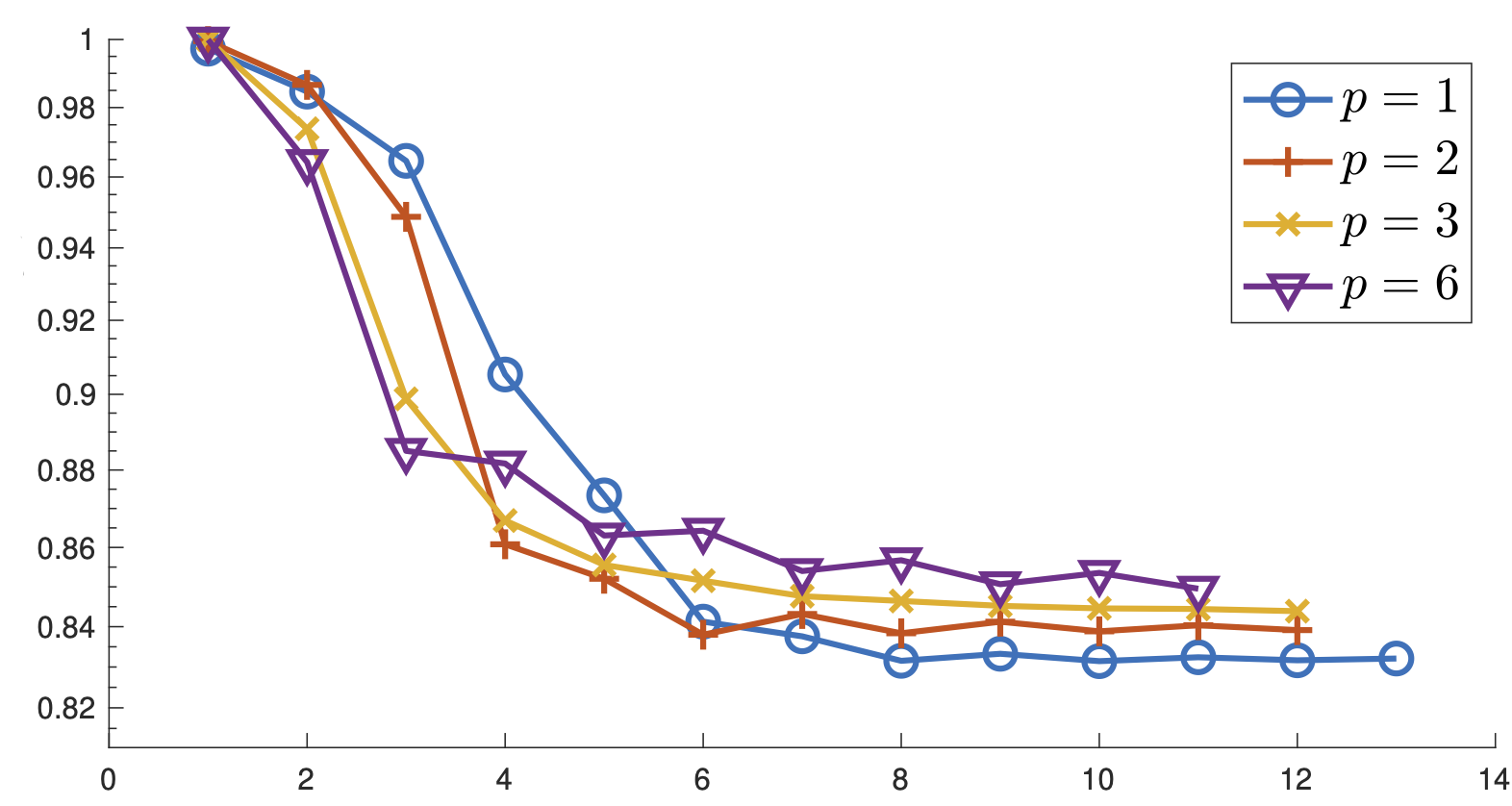}};
	\node[above=of img1, node distance=0cm, yshift=-1cm,font=\color{black}] 
	{DD solver, checkerboard solution, $J=28$};
	\node[below=of img1, node distance=0cm,yshift=1.2cm,font=\color{black}] 
	{\scriptsize iteration};
	\node[left=of img1, node distance=0cm, rotate=90, anchor=center, 
	xshift=0.2cm,yshift=-0.7cm,font=\color{black}] 
	{\scriptsize effectivity index};
	\end{tikzpicture}
\caption{Effectivity index 
$\etalg / \big\|\dif^{-1/2}(\u_{J} - \u_{J}^{i})\big\|$ of the
a~posteriori estimator $\etalg$ for the algebraic error 
$ \big\|\dif^{-1/2}(\u_{J} - \u_{J}^{i})\big\|$. 
Multigrid solver of Algorithm~\ref{Definition_solver} (left) 
and domain decomposition solver of Algorithm~\ref{AS_solver} (right).}
\label{fig_eff}
\end{figure}

Figure~\ref{fig_eff} gives the effectivity indices of the built-in 
a~posteriori estimators $\etalg$ of the algebraic error $\|\dif^{-1/2}(\u_{J} - \u_{J}^{i})\big\|$ for the multigrid solver 
(left) and the domain decomposition solver (right) for the different test cases. 
Recall that $\etalg$ are respectively given by~\eqref{def_eta} 
and~\eqref{def_eta_DD}, whereas the effectivity indices are given by 
$\etalg / \big\|\dif^{-1/2}(\u_{J} - \u_{J}^{i})\big\|$. 
From Theorem~\ref{thm_upper_bound}, we know that the effectivity indices have 
to be smaller than or equal to $1$ but cannot drop below the constant 
$\beta$, 
and this uniformly in the polynomial degree $p$. This $p$-robustness is 
confirmed in Figure~\ref{fig_eff}.

\subsection{Contraction factors of the solvers}

\begin{figure}[htb]
\center
	\begin{tikzpicture}
	\node (img1)  {\includegraphics[scale=0.165]{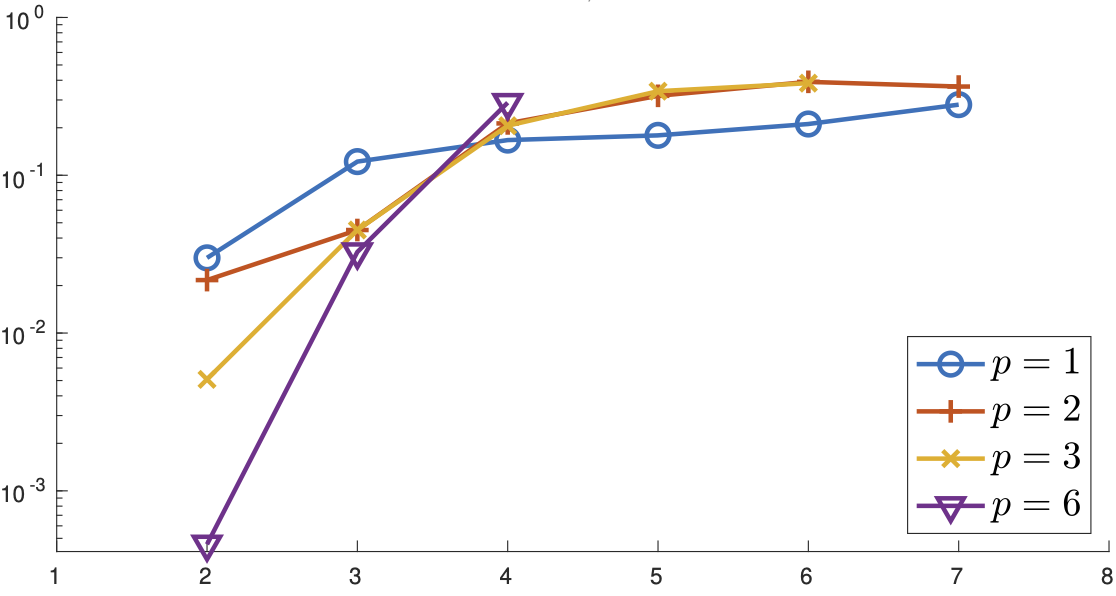}};
	\node[above=of img1, node distance=0cm,yshift=-1cm,font=\color{black}] 
	{MG solver, smooth solution, $J=5$};
	\node[below=of img1, node distance=0cm,yshift=1cm,font=\color{black}] 
	{\scriptsize iteration};
	\node[left=of img1, node distance=0cm, rotate=90, anchor=center, 
	xshift=0.2cm,yshift=-0.7cm,font=\color{black}] 
	{\scriptsize contraction factor};
	\end{tikzpicture}
\hfil
	\begin{tikzpicture}
	\node (img1)  {\includegraphics[scale=0.165]{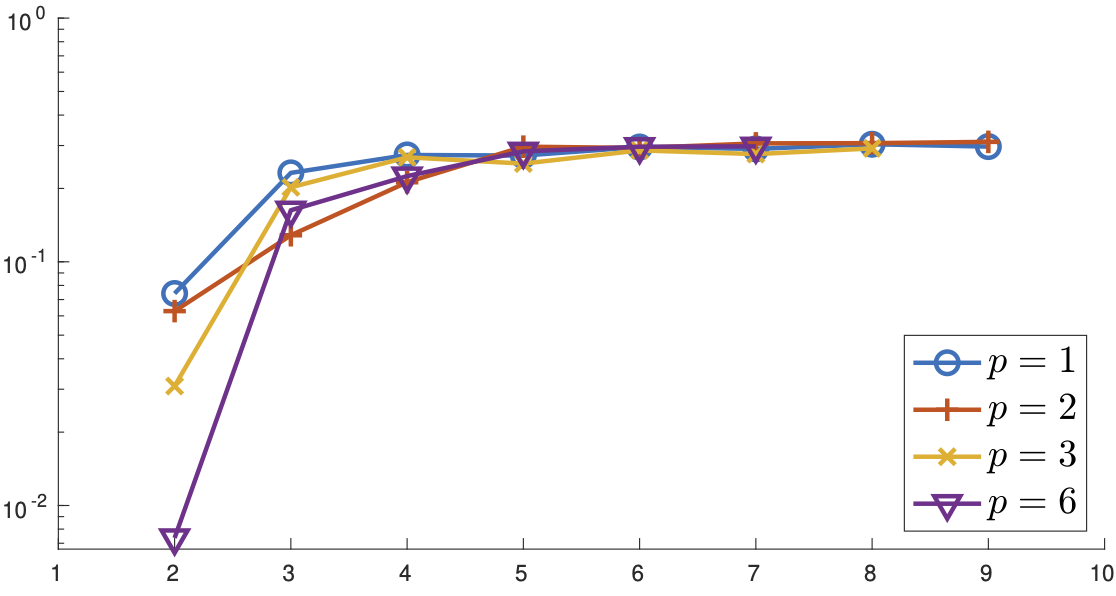}};
	\node[above=of img1, node distance=0cm, yshift=-1cm,font=\color{black}] 
	{DD solver, smooth solution, $J=5$};
	\node[below=of img1, node distance=0cm,yshift=1cm,font=\color{black}] 
	{\scriptsize iteration};
	\node[left=of img1, node distance=0cm, rotate=90, anchor=center, 
	xshift=0.2cm,yshift=-0.7cm,font=\color{black}] 
	{\scriptsize contraction factor};
	\end{tikzpicture}
\\
	\begin{tikzpicture}
	\node (img1)  {\includegraphics[scale=0.16]{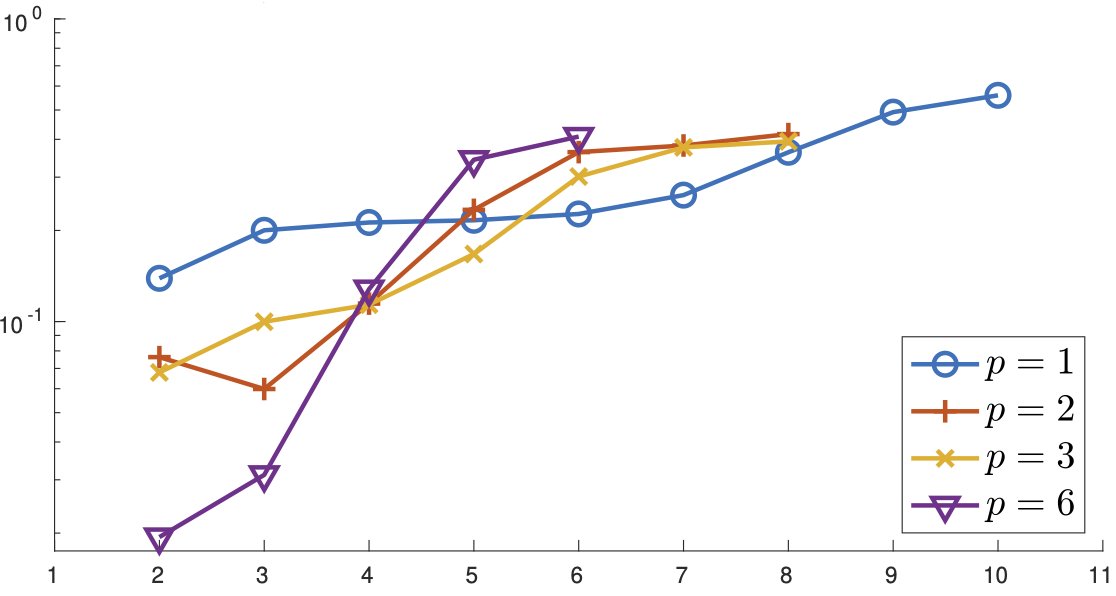}};
	\node[above=of img1, node distance=0cm,yshift=-1cm,font=\color{black}] 
	{MG solver, well wavefront solution, $J=12$};
	\node[below=of img1, node distance=0cm,yshift=1cm,font=\color{black}] 
	{\scriptsize iteration};
	\node[left=of img1, node distance=0cm, rotate=90, anchor=center, 
	xshift=0.2cm,yshift=-0.7cm,font=\color{black}] 
	{\scriptsize contraction factor};
	\end{tikzpicture}
\hfil
	\begin{tikzpicture}
	\node (img1)  {\includegraphics[scale=0.16]{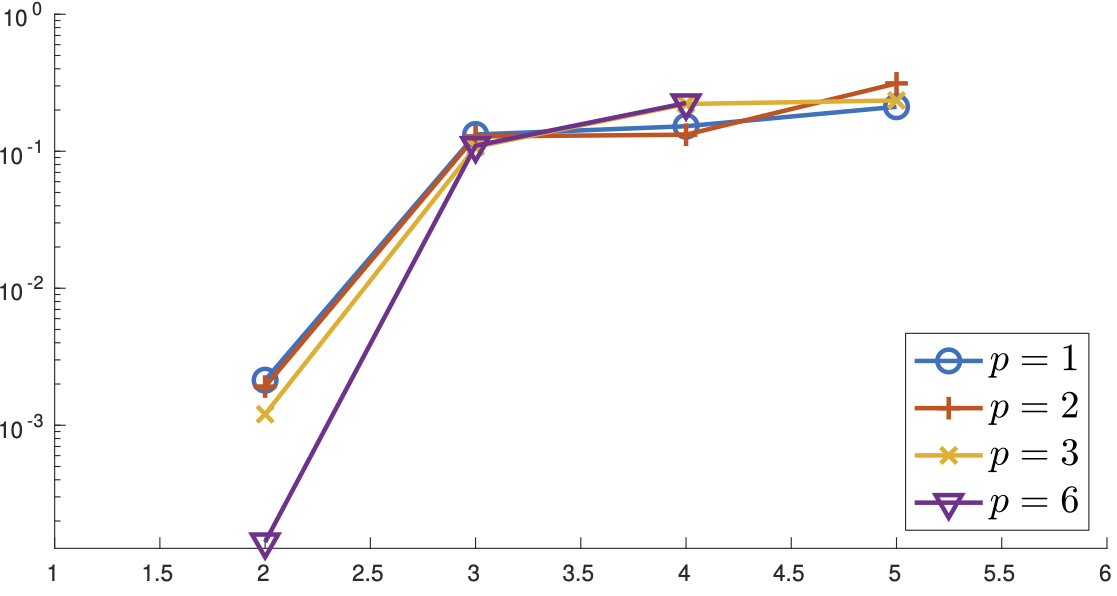}};
	\node[above=of img1, node distance=0cm, yshift=-1cm,font=\color{black}] 
	{DD solver, well wavefront solution, $J=12$};
	\node[below=of img1, node distance=0cm,yshift=1cm,font=\color{black}] 
	{\scriptsize iteration};
	\node[left=of img1, node distance=0cm, rotate=90, anchor=center, 
	xshift=0.2cm,yshift=-0.7cm,font=\color{black}] 
	{\scriptsize contraction factor};
	\end{tikzpicture}
\\
	\begin{tikzpicture}
	\node (img1)  {\includegraphics[scale=0.16]{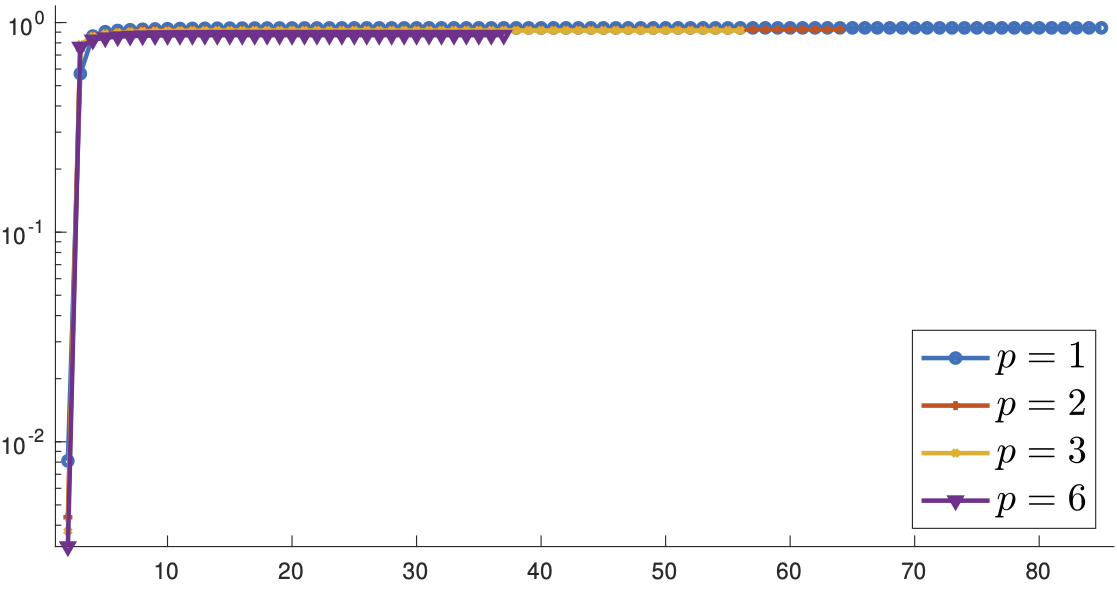}};
	\node[above=of img1, node distance=0cm,yshift=-1cm,font=\color{black}] 
	{MG solver, checkerboard solution, $J=28$};
	\node[below=of img1, node distance=0cm,yshift=1cm,font=\color{black}] 
	{\scriptsize iteration};
	\node[left=of img1, node distance=0cm, rotate=90, anchor=center, 
	xshift=0.2cm,yshift=-0.7cm,font=\color{black}] 
	{\scriptsize contraction factor};
	\end{tikzpicture}
\hfil
	\begin{tikzpicture}
	\node (img1)  {\includegraphics[scale=0.16]{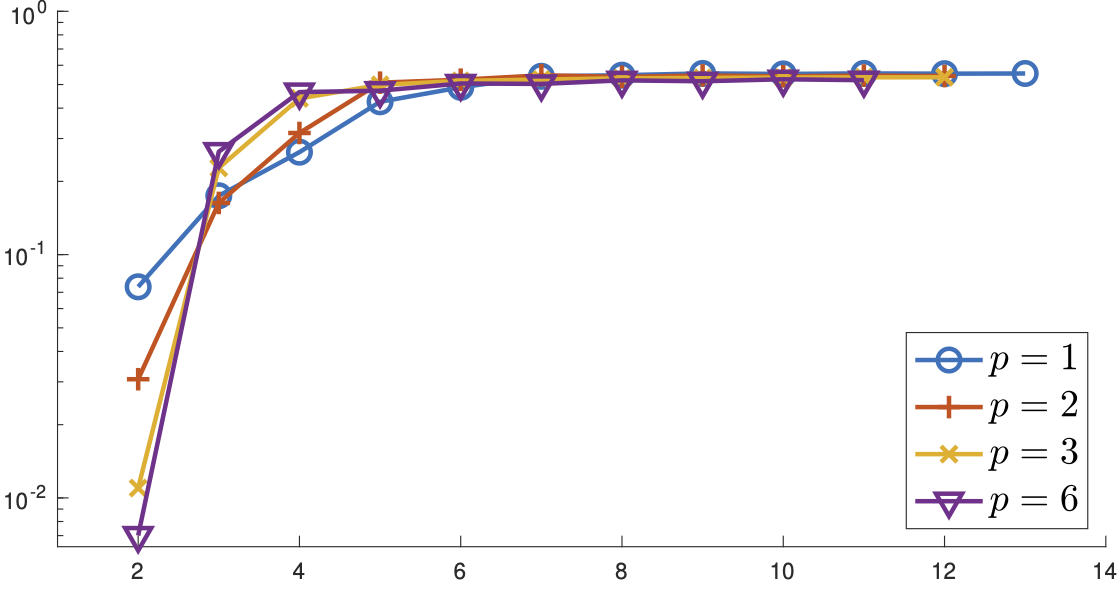}};
	\node[above=of img1, node distance=0cm, yshift=-1cm,font=\color{black}] 
	{DD solver, checkerboard solution, $J=28$};
	\node[below=of img1, node distance=0cm,yshift=1cm,font=\color{black}] 
	{\scriptsize iteration};
	\node[left=of img1, node distance=0cm, rotate=90, anchor=center, 
	xshift=0.2cm,yshift=-0.7cm,font=\color{black}] 
	{\scriptsize contraction factor};
	\end{tikzpicture}
\caption{Contraction factors of the solvers given by 
$\big\|\dif^{-1/2}(\u_{J} - \u_{J}^{i+1})\big\| 
/ \big\|\dif^{-1/2}(\u_{J} - \u_{J}^{i})\big\|$. 
Multigrid solver of Algorithm~\ref{Definition_solver} (left) 
and domain decomposition solver of Algorithm~\ref{AS_solver} (right).}
\label{fig_contr}
\end{figure}

In Figure~\ref{fig_contr}, we report the contraction factors 
$\big\|\dif^{-1/2}(\u_{J} - \u_{J}^{i+1})\big\| 
/ \big\|\dif^{-1/2}(\u_{J} - \u_{J}^{i})\big\|$. 
From Theorem~\ref{thm_converg}, these have to be bounded by the constant 
$\alpha$, in particular uniformly in the polynomial degrees $p$, 
which we indeed observe. More precisely, often, even a stronger 
initial contraction for higher $p$ appears, and then the contraction 
behaves very similarly for different polynomial degrees throughout the rest 
of the iterations.

\subsection{\texorpdfstring{$p$}{p}-robustness in the number of iterations}

Table~\ref{Tab_iter_nonadapt} shows how the $p$-robustness of the solvers 
translates to stable (or even decreasing) iteration numbers necessary to reduce 
the initial a~posteriori algebraic error estimator $\etalg$ by the 
factor $10^{5}$. The performance of the multigrid and the domain decomposition 
solvers appears to be very similar for the smooth solution and 
uniform refinement case, though the multigrid iteration numbers are slightly 
smaller overall.
Robust numerical performance is also seen in the well wavefront case despite 
a rougher analytic solution and a larger number of levels in the graded mesh 
hierarchy.
Though we do not refine uniformly here, we still numerically observe no 
degradation of number of iterations for increasing $p$ with respect to the 
number of mesh levels $J$, probably since the solution is regular enough. It is 
interesting to note that the domain decomposition method requires slightly 
smaller iteration numbers compared to the multigrid. An interpretation of this 
observation can be that as the hierarchy has more levels, the subdomain problems 
can correct the algebraic error much better due to the size of the local 
problems compared to the smaller patches used in the multigrid.
In the checkerboard case exhibiting a singular solution, the iteration numbers 
increase with the growth of the number of mesh levels, in at most a linear way, 
which is in accordance with our main theoretical results; observe, however, that 
the iteration numbers remain stable (or decrease) with respect to $p$. This can 
be seen more easily in the case of the multigrid solver, whereas the increase is 
very mild for the domain decomposition method; therein, however, much larger 
subdomain problems have to be solved.

\begin{table}[tbhp]
\caption{Number of iterations needed for the multigrid solver and domain 
decomposition method to decrease the relative a~posteriori estimator of the 
algebraic error $\etalg$ by $10^5$,
for the test cases~\eqref{Smooth}--\eqref{Checkerboard} and different polynomial 
degrees $p$, number of mesh levels $J$.}
\label{Tab_iter_nonadapt}
{\footnotesize
\begin{center}

\begin{tabular}{c||cc|cc|cc||cc|cc|cc||cc|cc|cc}
\multicolumn{1}{c||}{}
&\multicolumn{6}{c||}{Smooth}                            
&\multicolumn{6}{c||}{Wellwavefront}                     
&\multicolumn{6}{c}{Checkerboard}                         
\\ 
\multirow{2}{*}{$p$}
&\multicolumn{2}{c|}{
	\multirow{2}{*}{\begin{tabular}[c]{@{}c@{}}$ J = 3$\\  MG DD \end{tabular}}}
&\multicolumn{2}{c|}{
	\multirow{2}{*}{\begin{tabular}[c]{@{}c@{}}$ J = 4$\\  MG DD\end{tabular}}}
&\multicolumn{2}{c||}{
	\multirow{2}{*}{\begin{tabular}[c]{@{}c@{}}$ J = 5$\\  MG DD\end{tabular}}}
&\multicolumn{2}{c|}{
	\multirow{2}{*}{\begin{tabular}[c]{@{}c@{}}$ J = 4$\\  MG DD\end{tabular}}}
&\multicolumn{2}{c|}{
	\multirow{2}{*}{\begin{tabular}[c]{@{}c@{}}$ J = 8$\\  MG DD\end{tabular}}}
&\multicolumn{2}{c||}{
	\multirow{2}{*}{\begin{tabular}[c]{@{}c@{}}$ J = 12$\\  MG DD\end{tabular}}}
&\multicolumn{2}{c|}{
	\multirow{2}{*}{\begin{tabular}[c|]{@{}c@{}}$ J = 7$\\  MG DD\end{tabular}}}
&\multicolumn{2}{c|}{
	\multirow{2}{*}{\begin{tabular}[c|]{@{}c@{}}$ J = 14$\\ MG DD\end{tabular}}}
&\multicolumn{2}{c}{
	\multirow{2}{*}{\begin{tabular}[c]{@{}c@{}}$ J = 28$\\  MG DD\end{tabular}}}
\\
&\multicolumn{2}{c|}{}
&\multicolumn{2}{c|}{}
&\multicolumn{2}{c||}{}
&\multicolumn{2}{c|}{}
&\multicolumn{2}{c|}{}
&\multicolumn{2}{c||}{}
&\multicolumn{2}{c|}{}
&\multicolumn{2}{c|}{}
&\multicolumn{2}{c}{}  \\
\hline
1 & 8&9& 8&  9& 7& 9& 12& 9& 10& 5& 10& 5& 19 &10& 37&13& 85& 13 \\
2 & 9&9& 8&  9& 7& 9& 12& 7& 10& 5&  8& 5& 16 &11& 30&12& 64& 12 \\
3 & 8&8& 7&  8& 6& 8& 10& 7&  9& 5&  8& 5& 15 &11& 26&11& 56& 12 \\
6 & 5&7& 4&  7& 4& 7&  9& 4&  8& 4&  6& 4& 11 &10& 18&10& 37& 11 \\
\end{tabular}
\end{center}
}
\end{table}

\subsection{Adaptive number of smoothing steps in the multigrid solver}

In addition to the efficiency of the estimator $\etalg$, 
Theorem~\ref{thm_upper_bound} together with Theorem~\ref{thm_error_contr} 
have yet another important consequence: from~\eqref{loc_eta}, $\etalg$ provides 
a localized estimation of the algebraic error by levels (and by patches). This 
is the starting point for an adaptive extension of the multigrid solver 
following~\cite[Definition~7.1]{Mir_Pap_Voh_lambda}. Here we run an experiment with the solver setting the adaptivity parameter $\theta_{\mathrm{alg}}$ to 0.3.
Instead of employing just one post-smoothing step in 
Algorithm~\ref{Definition_solver} on each level (which, from 
Theorem~\ref{thm_converg}, in contrast to the usual case, is sufficient for 
overall contraction), an adaptive criterion based on the localized 
writing of $\etalg$ 
discerns whether 
additional smoothing per level is beneficial. 
More precisely, we continue smoothing on the given level $j$ if the contribution of the squared error decrease from level $j$, $(\lambda_j^i)^2 \big\|\dif^{-1/2}\brho_j^i\big\|^2 = \lambda_j^i \sum_{\ver \in \V_j} \big\|\dif^{-1/2}\brho_{j,\ver}^i\big\|^2_{\omaj}$, is significant (employing the parameter~$\theta_{\mathrm{alg}}$).
In Table~\ref{Tab_iter_adapt}, 
we present a comparison with the non-adaptive variant for the well wavefront 
test case. The comparison is done in terms of the iteration numbers (\#iter.) 
coinciding with the total number of V-cycles needed to decrease the relative 
a~posteriori estimator of the algebraic error $\etalg$ by the factor $10^5$, 
of global synchronizations (\#sync.) coinciding with the sum over iterations 
of the number of smoothing steps and coarse solves, and we moreover display 
the mean and maximum number of smoothing steps (\#smooth.).
This approach aims to reduce costs by: 1) lowering the overall number of 
V-cycles, since intergrid operators for a hierarchy consisting of many levels 
can become more costly; 2) employing more smoothing steps on lower levels where 
there are less patches and smoothing is cheaper. As we see from 
Table~\ref{Tab_iter_adapt}, the estimator correctly identifies levels 
requiring more smoothing leading to improved numerical performance by adding
one to two more post-smoothing steps only.


\begin{table}[tbhp]
\caption{Comparison of Algorithm~\ref{Definition_solver} and its version with 
adaptive number of smoothing steps for the well wavefront test 
case~\eqref{Wellwavefront}, different polynomial degrees $p$ and number of mesh 
levels $J$.}
\label{Tab_iter_adapt}
\resizebox{\textwidth}{!}{
\begin{tabular}{c||cccccc|cccccc|cccccc}
&\multicolumn{6}{c|}{\multirow{2}{*}{\begin{tabular}[c]{@{}c@{}}
	$ J = 4$\\ \ \ \#iter. \ \ \ \ \ \  
	\#sync.  \ \ \ \ \#smooth. \end{tabular}}}
&\multicolumn{6}{c|}{\multirow{2}{*}{\begin{tabular}[c]{@{}c@{}}
	$ J = 8$\\   \ \ \#iter. \ \ \ \ \ \  
	\#sync.   \ \ \ \ \#smooth.  \end{tabular}}}
&\multicolumn{6}{c}{\multirow{2}{*}{\begin{tabular}[c]{@{}c@{}}
	$ J = 12$\\   \ \ \#iter. \ \ \ \ \ \  
	\#sync.   \ \ \ \ \#smooth.  \end{tabular}}}
\\
\multirow{2}{*}{}
\multirow{3}{*}{$p$}
&\multicolumn{6}{c|}{\multirow{2}{*}{\begin{tabular}[c]{@{}c@{}} 
	\\adapt/MG \  adapt/MG \ mean/max \end{tabular}}}
&\multicolumn{6}{c|}{\multirow{2}{*}{\begin{tabular}[c]{@{}c@{}} 
	\\ adapt/MG  \ adapt/MG  \ mean/max \end{tabular}}}
&\multicolumn{6}{c}{\multirow{2}{*}{\begin{tabular}[c]{@{}c@{}} 
	\\ adapt/MG \ adapt/MG \ mean/max \end{tabular}}}
\\
&\multicolumn{6}{c|}{}
&\multicolumn{6}{c|}{}
&\multicolumn{6}{c}{}  
\\
\hline
1 & &7/12 &&& 55/60  & 1.83/3 & & 7/10 &&& 98/90  & 1.81/2 & 
& 7/10& & & 133/130   &  1.71/3 \ \  \\
2 & &6/12 &&& 47/60  & 1.88/3 & & 6/10 &&& 83/90  & 1.84/2 &
& 		6/8& && 111/104 & 1.71/3 \ \ \\
3 & &6/10 &&&  44/50  & 1.76/2 & & 5/9 \ &&& 66/81  & 1.83/3 &
& 		5/8& && 83/104 & 1.60/3 \ \  \\
6 & &5/9 \  &&& 35/45  & 1.75/2 & & 5/8 \ &&& 63/72 & 1.80/3 &
&    	4/6 &&& 56/78  & 1.44/2 \ \  \\
\end{tabular}
}
\end{table}

\subsection{Localization of the algebraic error on levels and patches}\label{sec_num_loc}

\begin{figure}[ht]
	\includegraphics[width=0.47\textwidth]{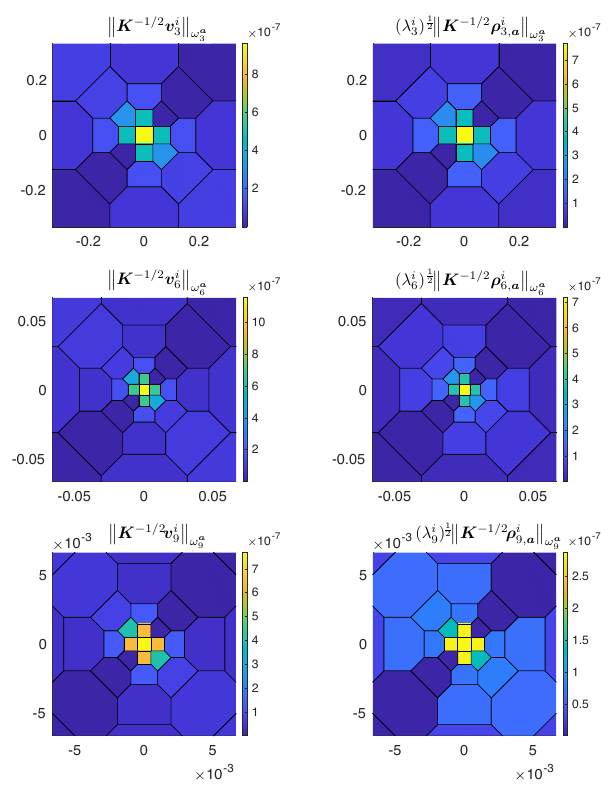}
	\hfill
	\includegraphics[width=0.47\textwidth]{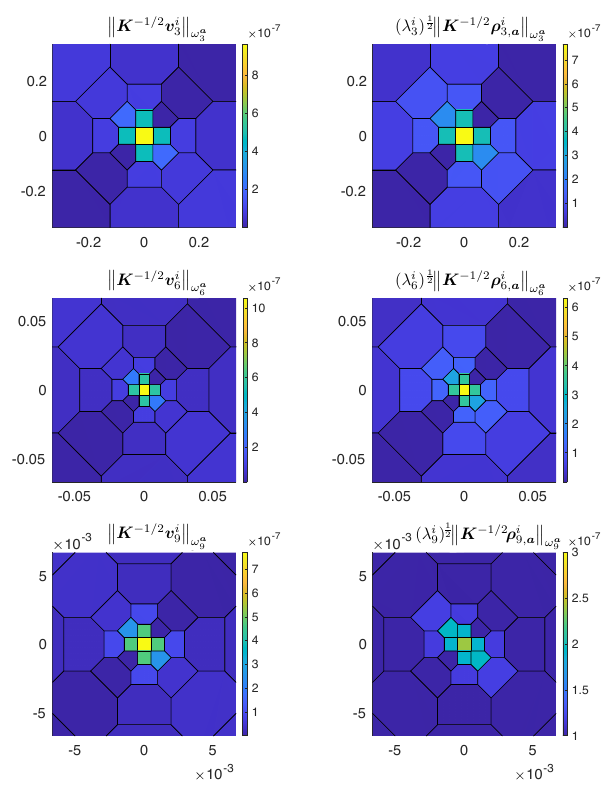}
\caption{Checkerboard problem, $J=9$, zoomed towards the origin. 
Left: patchwise distribution of the contribution of the algebraic error $\big\|\dif^{-1/2}(\u_{J} - \u_{J}^{i})\big\|$ to the given level $j$ on iteration $i=7$ (denoted as $\big\|\dif^{-1/2}\bm{v}_{j}^i\big\|_{\omaj}$)
(first column) compared to the a~posteriori estimator of this algebraic error 
contribution $(\lambda_j^i)^{\frac 1 2} \big\|\dif^{-1/2}\brho_{j,\ver}^i\big\|_{\omaj}$ (second column) for selected levels $j=3,6,9$ in the hierarchy and $p=1$. 
Right: same experiment for $p=3$.}
\label{fig_patch_err_est}
\end{figure}

While in the previous section we focused on 
adaptivity based on the levelwise 
localization of the 
efficient a~posteriori estimator $\etalg$,
\[	\big( \etalg \big)^2 = \sum_{j=0}^J 
	\big(\lambda_j^i \big\|\dif^{-1/2}\brho_j^i\big\|\big)^2,
\]
we now investigate 
numerically the accuracy of the localization both in patches and levels, 
\[
    \big( \etalg \big)^2 
	= \big\|\dif^{-1/2}\brho_0^i\big\|^2  
	+\sum_{j=1}^J \lambda_j^i \sum_{\ver \in \V_j} 
	\big\|\dif^{-1/2}\brho_{j,\ver}^i\big\|^2_{\omaj}.
\]
In Figure~\ref{fig_patch_err_est}, we present the patch contributions for 
different levels of the hierarchy. Since the patches overlap, for ease of 
visualization, we use the Vorono\"i dual graph of the triangular mesh 
connecting the centers of each triangle leading to non-overlapping polygonal 
representations of the patches. The estimator captures correctly the 
distribution and the magnitude of the algebraic error throughout the mesh 
hierarchy. It is interesting to point out that the algebraic error is situated 
in the regions related to mesh refinement. This is in line with a rich 
literature on multilevel methods for graded meshes. From the adaptive meshes 
perspective, smoothing should take place in the patches that were modified in 
the refinement step, see Section~\ref{sec_loc_MG} below.
From the multilevel solver perspective, the local refinement process can also 
be steered by fully adaptive multigrid method, see, e.g. 
R\"ude~\cite{Rud_ful_adpt_MG_93}.

\subsection{Local smoothing in the multigrid solver} \label{sec_loc_MG}

Local solvers as in Chen, Nochetto, and Xu~\cite{Chen_Noch_Xu_MG_loc_ref_12},
Wu and Zheng~\cite{wz2017}, or 
Innerberger et al.~\cite{Inn_Mir_Praet_Strei_hp_MG_24} allow the solver to 
tackle only smaller regions in the smoothing phase on each level: this 
crucially reduces the complexity of the solver, a question that becomes vital 
when using solvers in combination with adaptive mesh refinement algorithms. 
Thus, we present here the \emph{local} multigrid solver as a simple modification 
of our multigrid approach of Algorithm~\ref{Definition_solver}: instead of 
smoothing on all patches on each level, we only employ smoothing to patches 
associated with vertices which are new on the given level with respect to the 
previous one, as well as existing vertices whose patch has shrunk in size in the 
mesh refinement step. Recall that each patchwise problem solved via the V-cycle 
of multigrid provides a local error correction. In comparison, by only smoothing 
locally relative to the mesh-refinement, the V-cycle of the local multigrid will 
only give certain -- the bulk -- of the error corrections of the non-local 
multigrid. This is to say that the larger error components are localized 
relative to the mesh-refinement and precisely these components are 
tackled by local smoothing. 

In Table~\ref{Tab_locMG}, we present a comparison of the original 
Algorithm~\ref{Definition_solver} (MG) with its local modification (locMG). 
We focus on three aspects. 
1) To showcase the localization advantage of locMG 
through the savings in terms of the number of patch-wise solves, we display the 
ratio of the average number of patches where smoothing applies in the locMG to 
all patches (the MG case) (\# patches). 
2) To present the quality of the error correction, 
we rely on~\eqref{error_contr} and~\eqref{loc_eta} to express the computable 
error decrease of multigrid as
\[
	\big\|\dif^{-1/2}(\u_{J}  - \u_{J}^{i+1})\big\|^2
	= \big\|\dif^{-1/2}(\u_{J} - \u_{J}^{i})\big\|^2
	- \big\|\dif^{-1/2}\brho_0^i\big\|^2  
	-\sum_{j=1}^J \lambda_j^i \sum_{\ver \in \V_j} 
	\big\|\dif^{-1/2}\brho_{j,\ver}^i\big\|^2_{\omaj}.
\]
Similarly, we have for the local version of multigrid
\[
	\big\|\dif^{-1/2}(\u_{J}  - \u_{J}^{i+1})\big\|^2
	= \big\|\dif^{-1/2}(\u_{J} - \u_{J}^{i})\big\|^2
	- \big\|\dif^{-1/2}\brho_0^i\big\|^2  
	-\sum_{j=1}^J \lambda_j^i \sum_{\ver \in \V_j^+} 
	\big\|\dif^{-1/2}\brho_{j,\ver}^i\big\|^2_{\omaj},
\]
where now the sum is taken only over the new vertices and immediate neighbors 
$\V_j^+ \subseteq \V_j$. Denoting by $i_\textrm{MG}$ and $i_\textrm{locMG}$ the number of
iterations of the multigrid and local multigrid, respectively, 
the average computable decrease in MG is given by
\begin{align}\label{av_error_corr_MG}
	\frac{1}{i_\textrm{MG}} \sum_{i=1}^{i_\textrm{MG}} 
	\Big( \big\|\dif^{-1/2}\brho_0^i\big\|^2  
	+\sum_{j=1}^J \lambda_j^i \sum_{\ver \in \V_j} 
	\big\|\dif^{-1/2}\brho_{j,\ver}^i\big\|^2_{\omaj}\Big). 
\end{align}
Likewise, the average computable decrease in locMG is given by
\begin{align}\label{av_error_corr_locMG}
	\frac{1}{i_\textrm{locMG}} \sum_{i=1}^{i_\textrm{locMG}} 
	\Big( \big\|\dif^{-1/2}\brho_0^i\big\|^2  
	+\sum_{j=1}^J \lambda_j^i \sum_{\ver \in \V_j^+} 
	\big\|\dif^{-1/2}\brho_{j,\ver}^i\big\|^2_{\omaj}\Big). 
\end{align}
Thus, we present in Table~\ref{Tab_locMG} 
the ratio of~\eqref{av_error_corr_locMG} to~\eqref{av_error_corr_MG} 
to compare the error correction (error corr.).
3)~To compare a global performance, we show the ratio of the overall number of 
iterations of locMG to MG (iter.). One can see that the decrease of the 
algebraic error per iteration is indeed comparable for both methods, while the 
local multigrid employs a considerably smaller number of patchwise 
smoothings/local corrections on each level.

\begin{table}[tbhp]
\caption{Comparison of Algorithm~\ref{Definition_solver} and its version 
with local smoothing for the checkerboard test case~\eqref{Checkerboard}, 
different polynomial degrees $p$ and number of mesh levels $J$.}
\label{Tab_locMG}
\resizebox{\textwidth}{!}{
\begin{tabular}{c||cccccc|cccccc|cccccc}
&\multicolumn{6}{c|}{\multirow{2}{*}{\begin{tabular}[c]{@{}c@{}}
	$ J = 7$\\ \hspace{-1cm} \# patches \hspace{0.4cm}  error corr.  \hspace{0.6cm} iter. \end{tabular}}}
&\multicolumn{6}{c|}{\multirow{2}{*}{\begin{tabular}[c]{@{}c@{}}
	$ J = 14$\\  \hspace{-1cm} \# patches \hspace{0.4cm}  error corr.  \hspace{0.6cm} iter. \end{tabular}}}
&\multicolumn{6}{c}{\multirow{2}{*}{\begin{tabular}[c]{@{}c@{}}
	$ J = 28$\\  \hspace{-1cm} \# patches \hspace{0.4cm}  error corr.  \hspace{0.6cm} iter. \end{tabular}}}
\\
\multirow{2}{*}{}
\multirow{3}{*}{$p$}
&\multicolumn{6}{c|}{\multirow{2}{*}{\begin{tabular}[c]{@{}c@{}} 
	\\locMG/MG \  locMG/MG \ locMG/MG \end{tabular}}}
&\multicolumn{6}{c|}{\multirow{2}{*}{\begin{tabular}[c]{@{}c@{}} 
	\\ locMG/MG  \ locMG/MG  \ locMG/MG \end{tabular}}}
&\multicolumn{6}{c}{\multirow{2}{*}{\begin{tabular}[c]{@{}c@{}} 
	\\ locMG/MG \ locMG/MG \ locMG/MG \end{tabular}}}
\\
&\multicolumn{6}{c|}{}
&\multicolumn{6}{c|}{}
&\multicolumn{6}{c}{}  
\\
\hline
1  & &0.30 &&& \ \ 0.95  & \ \ \ 1.05 & & 0.31  &&& \ \ 0.89  & \  1.05 & 
& 		0.34& & & \ \ 0.94   & \ \ \  1.04 \ \  \\
2  & &0.30 &&&\ \  0.99  &\ \  \ 1 & & 0.30 &&& \ \  0.96   &\ 1.03 &
& 		0.32& && \ \  0.96 & \ \ \ 1.01 \ \ \\
3  & &0.30 &&&  \ \ 0.98  & \ \ \ 1 & & 0.31  &&& \ \ 0.95   & \ 1.04 &
& 		0.32& && \ \ 0.92 & \ \ \ 1.03 \ \  \\
6  & &0.30   &&& \ \ 0.99  & \ \ \ 1 & & 0.21  &&& \ \ 0.94 &  \ 1.06 &
&    	0.25 &&& \ \ 0.95  & \ \ \ 1.03 \ \  \\
\end{tabular}
}
\end{table}

Finally, let us mention that the local smoothing region can also be found 
adaptively in a way disconnected from the local mesh refinement, which is 
then also applicable in uniformly-refined grids, using that our estimators localize the algebraic error on levels and patches as discussed in Section~\ref{sec_num_loc}. Such an approach was designed 
in~\cite{Mir_Pap_Voh_W}. We do not implement and test it here.

\subsection{Non-simply connected domain}~\label{sec_half_whirlpool}
\begin{figure}[ht]
	\center
	\includegraphics[width=0.34\textwidth]{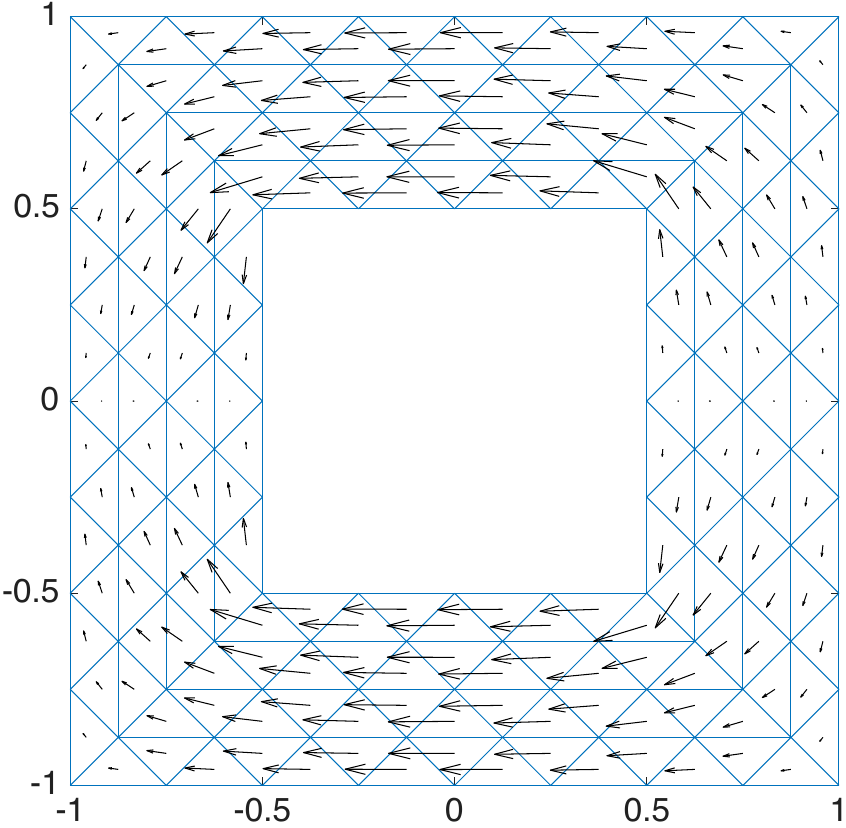} 
\caption{Approximate flow field $\u_{J}$ for the problem on a domain $\Omega$ not simply connected from Section~\ref{sec_half_whirlpool}.}
\label{fig:half_whirlpool}
\end{figure}

Let us now consider an experiment in a domain with a topology not covered by our assumption on $\Omega$ being simply connected. We consider the following test case:

\begin{itemize}
    \item \textbf{Half-whirlpool in a domain with a hole and uniform refinement.}
   Consider $ \Omega  =  (-1,1)^2 \setminus (-1/2,1/2)^2 $ with $\dif = I$ and the source term
    \begin{align}\label{half_whirlpool_f}
f (x, y) 
    = \begin{cases}
    \begin{array}{ll}
        1 \qquad&  x>0, \\
        -1 & x\le 0.
    \end{array}
    \end{cases}
    \end{align}
\end{itemize}
We do not have access to the analytical solution $\u$ in this case, but it is not needed for the following experiment using a uniformly refined mesh. An illustration of the discrete velocity flow field $\u_{J}$ is given in Figure~\ref{fig:half_whirlpool}, where one can see the half-whirlpool phenomenon occurring. 

\begin{figure}[!ht]
    \centering
    \includegraphics[width=0.48\linewidth]{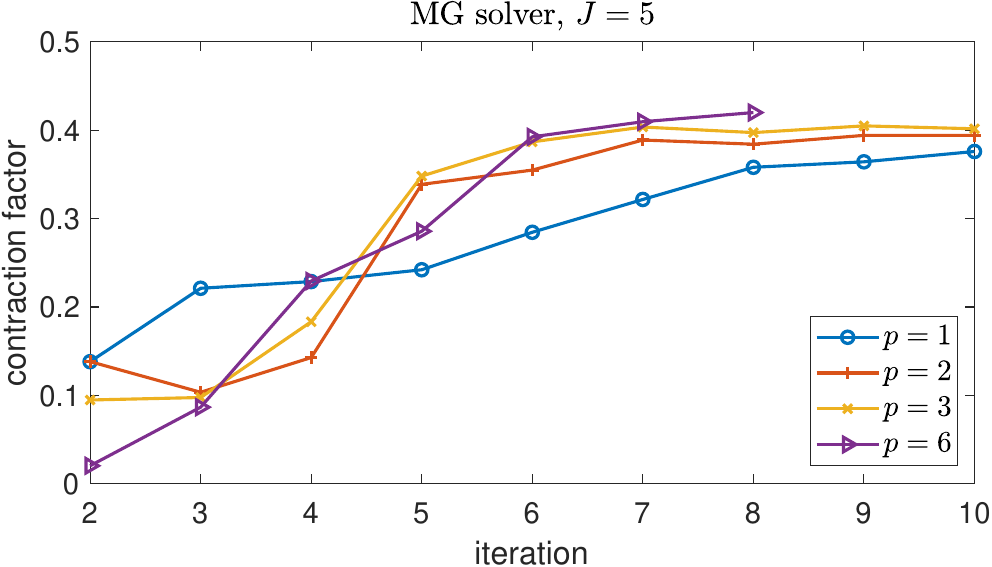}
    \hfill
    \includegraphics[width=0.48\linewidth]{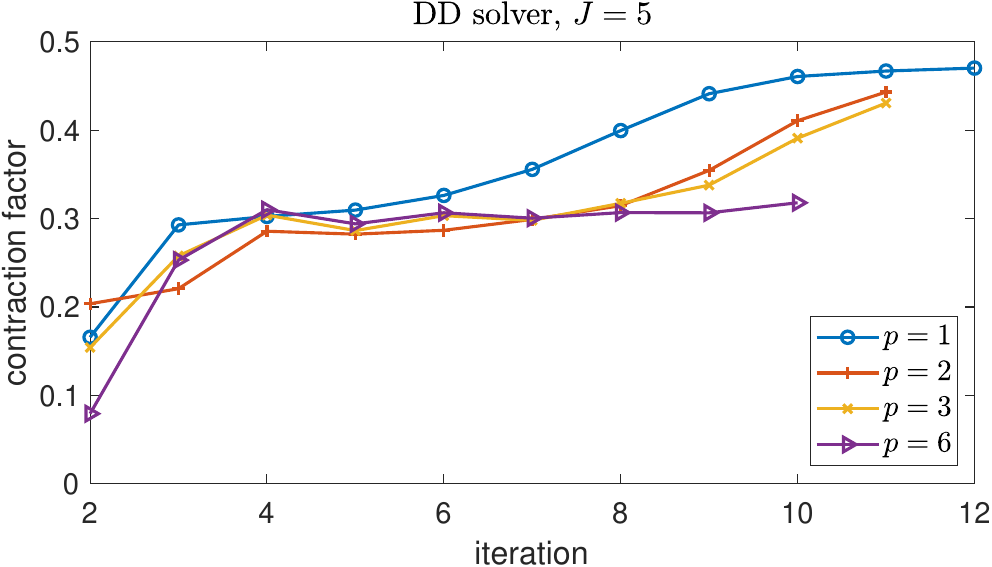}
    \caption{Contraction factors of the solvers given by 
$\big\|\dif^{-1/2}(\u_{J} - \u_{J}^{i+1})\big\| 
/ \big\|\dif^{-1/2}(\u_{J} - \u_{J}^{i})\big\|$. 
Multigrid solver of Algorithm~\ref{Definition_solver} (left) 
and domain decomposition solver of Algorithm~\ref{AS_solver} (right) 
for the problem from Section~\ref{sec_half_whirlpool}}
    \label{fig:half_whirlpool_experiments}
\end{figure}

Since $\Omega$ is not simply connected, this test is not covered by the analysis presented in this manuscript. 
Nonetheless, the results of the multigrid and domain 
decomposition methods can be seen not to be influenced by this change of domain 
topology. Indeed, the overall iteration numbers shown in 
Table~\ref{Tab_iter_half_whirlpool} remain similar to the other tests we have considered. 
This holds also for the contraction factors that are depicted in 
Figure~\ref{fig:half_whirlpool_experiments}.

\begin{table}[ht]
\caption{Number of iterations needed for the multigrid solver and domain 
decomposition method to decrease the relative a~posteriori estimator of the 
algebraic error $\etalg$ by $10^5$,
for the problem from Section~\ref{sec_half_whirlpool} and different polynomial 
degrees $p$, number of mesh levels $J$.}
\label{Tab_iter_half_whirlpool}
{\footnotesize
\begin{center}

\begin{tabular}{c||ccc|ccc|ccc}
\multicolumn{1}{c||}{}
&\multicolumn{9}{c}{Half-whirlpool on a domain with a cavity}                                        
\\ 
%
&\multicolumn{3}{c|}{\multirow{2}{*}{
    \begin{tabular}[c]{@{}c@{}}$ J = 3$\\
    \hspace{-0.2cm}MG \hspace{0.1cm}DD \hspace{0.1cm}\#DoF (div-free) \end{tabular}}}
&\multicolumn{3}{c|}{\multirow{2}{*}{
    \begin{tabular}[c]{@{}c@{}}$ J = 4$\\\hspace{-0.2cm}MG \hspace{0.1cm}DD \hspace{0.1cm}\#DoF (div-free)  \end{tabular}}}
 &\multicolumn{3}{c}{\multirow{2}{*}{
    \begin{tabular}[c]{@{}c@{}}$ J = 5$\\\hspace{-0.2cm}MG \hspace{0.1cm}DD \hspace{0.1cm}\#DoF (div-free)  \end{tabular}}}
\\
$p$&\multicolumn{3}{c|}{}
&\multicolumn{3}{c|}{}
&\multicolumn{3}{c}{}\\
\hline
1 & 11 & 12 & \num{15552} &  11 & 12 & \num{61824}  &  10 & 12 & \num{246528} \\
2 & 11 &  11 & \num{32544} &  10 &  11 & \num{129600} &  10 &  11 & \num{517248} \\
3 & 11 &  10 & \num{55680} &  10 &  11 & \num{221952} &  10 &  11 & \num{886272} \\
6 &  9 &  9 & \num{161952}&   9 &  10 & \num{646464} &   8 &  10 & \num{2583168}\\
\end{tabular}
\end{center}
}
\end{table}

\section{Multilevel stable decomposition for 
di\-ver\-gen\-ce-free Ra\-vi\-art--Thomas piecewise 
polynomials}\label{sec_stable_decompositions}

In our proof of Theorem~\ref{thm_upper_bound} 
(and thus of~Theorem~\ref{thm_converg}), which we present in 
Section~\ref{sec_proofs} below, a crucial ingredient is the existence of a 
suitable stable decomposition of the space $\bV_{J}^0$ defined 
in~\eqref{V_j_0}, i.e., the space of divergence-free Raviart--Thomas 
$\H_{0}({\div}; \Omega)$-conforming piecewise polynomials. The precise result we need is as follows:

\begin{proposition}[Multilevel $p$-robustly stable decomposition of 
$\bV_{J}^0$]\label{SD_assumption}
Let either Assumption~\ref{assumption_refinement_quasiuniformity} or 
Assumption~\ref{graded_grids} hold.
Then for any $\bv_{J} \in \bV_{J}^0$, there exist 
$\bv_0 \in \bV_0^0 \cap {\RT}_0 (\T_{0})$ 
(global coarse-grid lowest-order component) and
$\bv_{j,\ver} \in \bV_{j}^{\ver,0}$,
$ 1 \le j \le J$, $\ver \in \V_j$, 
(local levelwise and patchwise high-order components) 
such that the decomposition
\begin{subequations}\label{SD}
	\begin{align}
	\bv_{J} = \bv_0 + \sum_{j=1}^J \sum_{\ver \in \V_j} \bv_{j,\ver}\label{SD_1}
	\end{align}
is stable as
	\begin{align}
	\big\|\dif^{-1/2}\bv_0\big\|^2 +
	\sum_{j=1}^J \sum_{\ver \in \V_j}
	\big\|\dif^{-1/2}\bv_{j,\ver}\big\|^2_{\omaj}
	\le C_{\rm{SD}}^2\big\|\dif^{-1/2}\bv_{J}\big\|^2,\label{SD_2}
	\end{align}
\end{subequations}
where the positive constant $C_{\rm{SD}}$ only depends on: 
\begin{itemize}
    \setlength\itemsep{0ex}
    \item the mesh shape regularity parameter $\kappa_{\T}$,
    \item the parameters $C_{\rm qu}$ and $C_{\rm ref}$ when Assumption~\ref{assumption_refinement_quasiuniformity} is satisfied or the 
parameters $C^0_{\rm qu}$ and $ C_{\rm loc,qu}$ when Assumption~\ref{graded_grids} is satisfied, 
    \item the diffusion inhomogeneity or anisotropy ratio $\Lambda_{\max}/\Lambda_{\min}$, 
    \item the domain $\Omega$ if $d=3$. 
\end{itemize}
In particular, $C_{\rm{SD}}$ is independent of 
the number of mesh levels~$J$ and of the polynomial degree $p$.
\end{proposition}

\section{Proof of the multilevel stable decomposition of Proposition~\ref{SD_assumption} for $d=2$}\label{sec_MLD_2D}

We establish here a proof of the multilevel stable decomposition of  Proposition~\ref{SD_assumption} in two space dimensions. 
For this purpose, we recall the multilevel $p$-robustly stable decomposition for any function of the two-dimensional stream-function space 
\begin{align}\label{stream_space}
    S_{J} := \{ v_{J} \in H^1_0(\Omega), \  v_{J}|_K \in \PP_{p+1}(K) \ \forall K \in \T_{J} \} = \PP_{p+1}(\T_{J}) \cap H^1_0(\Omega) ,
\end{align}
$p \geq 0$, in the 
$H^1_0(\Omega)$-norm $\| \nabla {\cdot}\|$. This result has been shown in 
\cite[Proposition~7.6]{Mir_Pap_Voh_19} by combining the results of Sch\"oberl et al.~\cite{SchMelPechZag_08} and Xu, Chen, and Nochetto~\cite{Xu_Chen_Noch_opt_MG_loc_ref_09}.

\begin{lemma}[Multilevel $p$-robustly stable decomposition of Lagrange piecewise polynomials for $d=2$]\label{lem_SD_Lag}
Let $d=2$. Under either Assumption~\ref{assumption_refinement_quasiuniformity} 
or Assumption~\ref{graded_grids}, for any $v_{J} \in S_{J}$, there exists a 
decomposition
\begin{subequations}\label{k_rob_SD}
	\begin{align}
	& v_{J} = 
	v_0+ \sum_{j=1}^J \sum_{\ver \in \V_{j}}  v_{j,\ver}, \label{k_rob_SD1}
	\end{align}
with 
$v_0 \in S_{0,1} := \{ v_{0} \in H^1_0(\Omega), \  v_{0}|_K \in \PP_{1}(K) \ \forall K \in \T_{0} \}$
(global coarse-grid lowest-order component) and 
$v_{j,\ver} \in S^{\ver}_{j} := \{ v_{j} \in H^1_0(\omaj), \  v_{j}|_K \in \PP_{p+1}(K) \ \forall K \in \T_{j}^\ver \}$, 
$1 \leq j \leq J$, $\ver \in \V_{j}$, 
(local levelwise and patchwise high-order components) stable as
	\begin{align}
	\| \nabla v_0 \|^2 
	+ \sum_{j=1}^J \sum_{\ver \in \V_{j}} 
	\| \nabla v_{j,\ver} \|^2_{\omaj} 
	\le  C_{\rm SD}^2 \| \nabla v_{J} \|^2, \label{k_rob_SD2}
	\end{align}
\end{subequations}
where the positive constant $C_{\rm{SD}}$ only depends on 
the mesh shape regularity parameter $\kappa_{\T}$ and additionally on
the parameters $ C_{\rm qu}$ and $C_{\rm ref}$ when 
Assumption~\ref{assumption_refinement_quasiuniformity} is satisfied or on the 
parameters $C^0_{\rm qu}$ and $ C_{\rm loc,qu}$ when 
Assumption~\ref{graded_grids} is satisfied.
In particular, $C_{\rm{SD}}$ is independent 
of the number of mesh levels $J$ and of the polynomial degree $p$.
\end{lemma}

With Lemma~\ref{lem_SD_Lag}, we can now prove Proposition~\ref{SD_assumption} in two space dimensions:

\begin{proof}[Proof of Proposition~\ref{SD_assumption} for $d=2$]
The construction of the stable decomposition for $d=2$ uses the property that, 
since $\Omega$ is simply connected, the discrete 
stream-function space of~\eqref{stream_space} satisfies
\begin{align}\label{curl_stream}
	\bV_{J}^0 = \nabla \times \, S_{J},
\end{align}
see e.g., Amrouche, Bernardi, Dauge, and Girault~\cite{Amr_Ber_Dau_Gir_98_potentials3d}
and Arnold~\cite[Section~4.5.5]{Arnold_FEEC_18}, where the two-dimensional curl operator is
\begin{align}\label{curl_1}
	\nabla \times s = \begin{pmatrix} s_y \\ - s_x \end{pmatrix} =
	\begin{pmatrix} 0 & 1 \\ -1 & 0 \end{pmatrix}
	\begin{pmatrix} s_x \\ s_y \end{pmatrix}=
	\begin{pmatrix} 0 & 1 \\ -1 & 0 \end{pmatrix}
	\nabla  s.
\end{align}
Since for all $\bv_{J}, \bw_{J} \in \bV^0_{J}$, 
by~\eqref{curl_stream}, there is $s_{J}, \theta_{J} \in S_{J}$ such that
$\bv_{J} = \nabla \times s_{J}$ and $\bw_{J} = \nabla \times \theta_{J}$, 
we obtain from~\eqref{curl_1} that
\begin{align}\label{curl_grad}
	(\dif^{-1}\bv_{J},\bw_{J}) 
	= (\dif^{-1}\nabla \times s_{J}, \nabla \times \theta_{J}) 
	= ({\mathbfcal{A}} \nabla s_{J}, \nabla \theta_{J})
\end{align}
for 
${\mathbfcal{A}} :=  \begin{pmatrix} 0 &  -1 \\ 1 &  0 \end{pmatrix} \dif^{-1}
\begin{pmatrix} 0 &  1 \\ -1 & 0 \end{pmatrix}$. 
In particular, note that $ (\dif^{-1}\bv_{J},\bw_{J}) =
( \nabla s_{J}, \nabla \theta_{J})$ when $\dif = I $.
Similar properties obviously hold on patches as well.

Let us now show~\eqref{SD} for $d=2$. Let $\bv_{J} \in \bV_{J}^0$. 
By~\eqref{curl_stream}, there is $v_{J} \in S_{J}$ such that
\begin{align}\label{curl_vh}
\bv_{J} = \nabla \times v_{J}.
\end{align}
We can now use the stable decomposition~of Lemma~\ref{lem_SD_Lag} for $v_{J}$: 
there exist $v_0 \in S_{0,1}$ and
$ v_{j,\ver} \in S^{\ver}_{j}$, $1 \le j \le J$, $\ver \in \V_j$, such that
$  v_{J} = v_0+ \sum_{j=1}^J \sum_{\ver \in \V_{j}}  v_{j,\ver}$ 
and~\eqref{k_rob_SD2} holds.
Define
\begin{align}\label{curl_vha}
	\bv_{0} 
	:= \nabla \times v_{0} \quad \text{and} \quad \bv_{j,\ver} 
	:= \nabla \times v_{j,\ver}.
\end{align}
Since $\bV_{0}^{0} \cap {\RT}_{0}(\T_{0}) = \nabla \times S_{0,1}$ and
$\bV_{j}^{\ver,0} = \nabla \times S^{\ver}_{j}$, $1 \le j \le J$, 
$\ver \in \V_j$,
we have $\bv_0 \in \bV_{0}^{0} cap {\RT}_0 (\T_{0})$ and $\bv_{j,\ver} \in \bV_{j}^{\ver,0}.$
Note that by applying $\nabla \times$ on both sides of~\eqref{k_rob_SD1}, 
we have
\begin{align*}
	\bv_{J}  
	= \bv_0 + \sum_{j=1}^J \sum_{\ver \in \V_{j}}  \bv_{j,\ver},  
	\quad \bv_0 \in \bV_{0}^{0} \cap {\RT}_0 (\T_{0}), \ \bv_{j,\ver} \in \bV_{j}^{\ver,0},
\end{align*}
which is the first part~\eqref{SD_1} of the result we want to prove.
Next, note that
\begin{align*}
	\big\| \bv_{0}\big\|^2 
	+	\sum_{j=1}^J \sum_{\ver \in \V_j}\big\| \bv_{j,\ver}\big\|^2_{\omaj}
	& \stackrel{\eqref{curl_vha}}=
	\big\| \nabla \times v_{0}\big\|^2 +
	\sum_{j=1}^J \sum_{\ver \in \V_j} 
	\big\|\nabla \times v_{j,\ver} \big\|^2_{\omaj}\\
	& \stackrel{\eqref{curl_1}}=
	\big\| \nabla v_{0}\big\|^2 +
	\sum_{j=1}^J \sum_{\ver \in \V_j}\big\| \nabla v_{j,\ver} \big\|^2_{\omaj}\\
	& \stackrel{\eqref{k_rob_SD2}}\le  C_{\rm SD}^2 \| \nabla v_{J} \|^2
	\stackrel{\eqref{curl_grad}}= C_{\rm SD}^2 \| \nabla \times v_{J} \|^2
	\stackrel{\eqref{curl_vh}}= C_{\rm SD}^2 \|  \bv_{J} \|^2.
\end{align*}
Finally, the result~\eqref{SD_2} 
follows once we take into account the variations of the 
diffusion coefficient~$\dif$.
\end{proof}

\section{Proof of the multilevel stable decomposition of Proposition~\ref{SD_assumption} for $d=3$}\label{sec_MLD_3D}

Proving Proposition~\ref{SD_assumption} in three space dimensions is harder than in two space dimensions. We proceed by establishing a collection of intermediate results that are of their own interest and presented below in the form of lemmas. We then employ them to prove Proposition~\ref{SD_assumption}.

\subsection{Stable discrete potential for the lowest-order divergence-free 
Raviart--Thomas piecewise polynomials}

In order to present the following result, we introduce the elementwise 
N\'ed\'elec spaces with $d=3$, see N\'ed\'elec~\cite{Ned_mix_R_3_80}, 
by ${\N}_p (K) := [\PP_{p}(K)]^3 + \bx \times \PP_{p}(K)^3$ 
for any element $K \in \T_{J}$. Similarly to~\eqref{eq_spaces}, 
let
\begin{equation}
	\bX \eq \H_{0}(\curl; \Omega) 
	:=  \{ \bv \in \H(\curl; \Omega), \bv \times \n =0 
	\text{ on } \partial \Omega  \}.
\end{equation}
Here, $\bv \times \n = 0$ on $\partial \Omega$ means that $(\nabla {\times} \bv, \bphi) - (\bv, \nabla {\times} \bphi) = 0$ for all $\bphi \in \H^1(\Omega)$.
Moreover, as for~\eqref{V_h}, define
\begin{equation} \label{eq_Ned_0}
	\bX_{J,0} \eq \{\bv_{J} \in \bX, \bv_{J}|_K \in \N_0 (K) \ 
	\forall  K \in \T_{J} \} = \N_0(\T_J) \cap \H_0(\curl; \Omega).
\end{equation}

\begin{lemma}[Stable discrete vector potential]\label{lem_stabpot}
Let $d=3$. For any $\bv_{J} \in \bV_{J}^0 \cap  {\RT}_0 (\T_{J})$, 
there is a discrete vector potential $\bxi_{J} \in \bX_{J,0}$ such that
\begin{align}
	\bv_{J} 
	= \nabla \times \bxi_{J} \quad \text{and} \quad
	\|\bxi_{J}  \| \le C_{\rm pot} \|\bv_{J}\|,\label{vh_curl}
\end{align}
where the constant $C_{\rm pot}$ only depends on the mesh shape regularity 
parameter $\kappa_{\T}$ and the 
domain~$\Omega$.
\end{lemma}

\begin{proof}
This result is called a discrete vector Poincar\'e inequality and can be obtained as in, e.g., Hiptmair and Xu~\cite[Section~4]{Hiptmair_Xu_07}, see Ern et al.~\cite{Ern_Guz_Potu_Voh_Poinc_disc_25} and the references therein for a general view. 
For completeness, we present a direct proof.
Consider the following discrete minimization problem
\begin{align}
	\bxi_{J}
	:=
	\underset{\substack{\bzeta_{J} \in \bX_{J,0} 
	\\ \nabla \times \bzeta_{J} = \bv_{J}}}{\rm arg \ min}
	\| \bzeta_{J} \|.
	\label{xi_argmin}
\end{align}
This problem is well posed since the datum $\bv_{J}$ is a divergence-free 
Raviart--Thomas $\H_{0}({\div}; \Omega)$-conforming piecewise polynomial and 
since $\Omega$ is simply connected, see, e.g., Amrouche, Bernardi, Dauge, and Girault~\cite{Amr_Ber_Dau_Gir_98_potentials3d}
and Arnold~\cite[Section~4.5.5]{Arnold_FEEC_18}.
We also define its counterpart on the 
continuous level,
\begin{align}
	\bxi
	:=
	\underset{\substack{\bzeta \in  \bX   
	\\\ \nabla \times \bzeta = \bv_{J}}}{\rm arg \ min}
	\| \bzeta \|.
	\label{xi_argmin_continuous}
\end{align}
We now employ the commuting projector ${\bm P}_{J}^{0, {\rm curl}}$ of 
Chaumont-Frelet 
and Vohral{\'i}k~\cite[Definition~2]{chaumontfrelet:hal-03817302}. 
Note that ${\bm P}_{J}^{0, {\rm curl}} (\bxi) \in \bX_{J,0}$ with 
$\nabla \times ({\bm P}_{J}^{0, {\rm curl}} (\bxi)) =
{\bm P}_{J}^{0, {\rm div}}(\nabla \times \bxi) = {\bm P}_{J}^{0, {\rm div}}(\bv_{J}) = \bv_{J}$, since $\nabla \times \bxi = \bv_{J} \in \bV_{J}^0 \cap  {\RT}_0 (\T_{J})$ and ${\bm P}_{J}^{0, {\rm div}}$ appearing in~\cite[Theorem~1, equation~(3.15)]{chaumontfrelet:hal-03817302} is a projector. 
Consequently, the minimization property~\eqref{xi_argmin} together 
with~\cite[Theorem~1, equation~(3.18)]{chaumontfrelet:hal-03817302} give
\begin{align}
	\| \bxi_{J} \|
	\le \| {\bm P}_{J}^{0, {\rm curl}} (\bxi) \|
	\le C_{\kappa_{\T}} \| \bxi \|
	\label{xi_disc_continuous_est}
\end{align}
for $ C_{\kappa_{\T}} > 0$ only depending on the shape regularity parameter 
$\kappa_{\T}$, since all discrete functions considered here are of lowest order.

Note now that~\eqref{xi_argmin_continuous} can be equivalently written as
searching for the solution $(\bxi,\bz)$, with $ \bxi \in \bX$ and 
$\bz \in  \bV^0 $, to the problem
\begin{subequations}\label{xi_mixed}
	\begin{alignat}{2}
	(\bxi, \bzeta)+ (\bz, \nabla \times \bzeta) & = 0
	\qquad \qquad & & \forall \bzeta \in  \bX, \label{xi_mixed1}  
	\\
	(\nabla \times \bxi, \bw)   & = (\bv_{J}, \bw)
	& & \forall \bw \in \bV^0.
	\label{xi_mixed2}
	\end{alignat}
\end{subequations}
In particular, note that $\bz \in \bV^0$, and from~\eqref{xi_mixed1} 
it satisfies
\begin{align}
	(\bz, \nabla \times \bzeta) 
	= -(\bxi, \bzeta) \quad \forall \bzeta \in  \bX.
	\label{zeta_Hcurl}
\end{align}
This means that $\bz$ belongs to $\H(\curl; \Omega)$ 
and satisfies $\nabla \times \bz = - \bxi$; altogether, recalling~\eqref{eq_spaces} and~\eqref{eq:Vz}, $\bz \in \H ({\div}; \Omega) \cap \H(\curl; \Omega)$ with $\bz {\cdot} \n = 0$ on $\partial \Omega$ and $\nabla {\cdot} \bz = 0$ in $\Omega$. Thus, by Poincaré--Friedrichs--Weber's inequality, see e.g. Weber~\cite{Weber_80}
or
Boffi, Brezzi, and Fortin~\cite[Section 11.1.2]{Boffi_Brezzi_Fortin_book},
\begin{align}
	\| \bz\| \le C_{\rm PFW} h_\Omega \| \nabla \times \bz \| 
	= C_{\rm PFW} h_\Omega \| \bxi \|,
\label{zeta_weber}
\end{align}
where $C_{\rm PFW}$ is a generic constant of order of unity possibly depending on $\Omega$
and $h_\Omega$ denotes the diameter of the domain $\Omega$.

To prove the stability in~\eqref{vh_curl} it is now sufficient to take 
$\bzeta = \bxi$ and $\bw = \bz$ as test functions in~\eqref{xi_mixed} and 
subtract the two resulting subequations, yielding
\begin{align}
	\| \bxi\|^2  = - (\bv_{J}, \bz)
	\le   \| \bv_{J} \| \| \bz \|
	\stackrel{\eqref{zeta_weber}}\le  
	C_{\rm PFW} h_\Omega \| \bv_{J} \| \| \bxi \|.
	\label{xi_stab}
\end{align}
Finally, combining with~\eqref{xi_disc_continuous_est},
\begin{align*}
	\| \bxi_{J} \|
	\stackrel{\eqref{xi_disc_continuous_est}}\le C_{\kappa_{\T}}   \| \bxi\|
	\stackrel{\eqref{xi_stab}}\le  
	\underbrace{C_{\kappa_{\T}} C_{\rm PFW}  
	h_\Omega}_{C_{\rm pot}} \| \bv_{J} \|.
\end{align*}
\end{proof}

\subsection{One-level stable decomposition for high-order divergence-free
Raviart--Thomas piecewise polynomials}

We now tackle a stable decomposition of $\bV_{J}^0$ from~\eqref{V_j_0} 
on the finest level $\T_{J}$, employing local divergence-free 
spaces $\bV_{J}^{\ver,0}$ from~\eqref{V_ja_0}.

\begin{lemma}[One-level $p$-robust stable decomposition of $\bV_{J}^0$]\label{SD_assumption_d3_OL}
Let $d=3$. 
For any $\bv_{J} \in \bV_{J}^0$, there exist 
$\bdelta_{J,0} \in \bV_{J}^0 \cap {\RT}_0 (\T_{J})$ 
(global lowest-order component) and 
$\bdelta_{J,\ver} \in \bV_{J}^{\ver,0}$, $\ver \in \V_{J}$, 
(local patchwise high-order components) such that the decomposition
\begin{subequations}\label{OL_3dSD}
	\begin{align}
	\bv_{J} = \bdelta_{J,0} 
	+ \sum_{\ver \in \V_{J}} \bdelta_{J,\ver}\label{OL_3dSD_1}
	\end{align}
	is stable as
	\begin{align}
	\big\|\dif^{-1/2}\bdelta_{J,0} \big\|^2 +
	\sum_{\ver \in \V_{J}}\big\|\dif^{-1/2}\bdelta_{J,\ver}\big\|^2_{\omaJ}
	\le C_{\rm{OL}}^2\big\|\dif^{-1/2}\bv_{J}\big\|^2,\label{OL_3dSD_2}
	\end{align}
\end{subequations}
where the constant $C_{\rm{OL}}$ only depends on the
mesh-geometry parameter $\kappa_{\T}$ and the diffusion inhomogeneity or 
anisotropy ratio $\Lambda_{\max}/\Lambda_{\min}$. In particular, $C_{\rm{OL}}$ is independent 
of the number of mesh levels $J$ and of the polynomial degree $p$.
\end{lemma}

\begin{proof}[Proof of Lemma~\ref{SD_assumption_d3_OL} ($p$-robust, relying on Falk and Winther~\cite{Falk_Winther_25})]
Let $\bv_{J} \in \bV_{J}^0$ be given and let 
\[
	\bX_{J} \eq \{\bxi_{J} \in \bX, \bxi_{J}|_K \in \N_p (K) \ 
	\forall  K \in \T_{J} \} = \N_p(\T_J) \cap \H_0(\curl; \Omega).
\]
During the revision of this ma\-nu\-script, there appeared the contribution by Falk and Winther~\cite{Falk_Winther_25} which allows for the following proof of~\eqref{OL_3dSD}. 
First, following~\cite[below equation~(1.1)]{Falk_Winther_25}, there exist operators $W^1, W^2$ and $B^1_\ver, B^2_\ver$ commuting in that $\nabla \times W^1 = W^2 \nabla \times $ and $\nabla \times B^1_\ver = B^2_\ver \nabla \times $. 
Second, since $\Omega$ is simply connected, divergence-free Raviart--Thomas $\H_0(\div; \Omega)$-conforming piecewise polynomials are curls of N\'ed\'elec $\H_0(\curl; \Omega)$-conforming piecewise polynomials. This means that the given $\bv_{J} \in \bV_{J}^0$ satisfies $\bv_{J} = \nabla \times \bxi_{J}$ for some $\bxi_{J} \in \bX_{J}$. 
Third, by~\cite[equation~(1.1)]{Falk_Winther_25}, for this $\bxi_{J}$, there holds $W^1 \bxi_{J} \in \bX_{J} \cap {\N}_0 (\T_{J})$ (global lowest-order component) and $B^1_\ver \bxi_{J} \in \bX_{J}^{\ver} \eq \{ \bzeta_J \in \bX_{J} |_{\omaJ},
	\ \bzeta_J {\times}\n = 0 \text{ on } \partial \omaJ \}$, $\ver \in \V_{J}$, (local patchwise high-order components) and there holds the decomposition
	\begin{align*}
	\bxi_{J} = W^1 \bxi_{J} 
	+ \sum_{\ver \in \V_{J}} B^1_\ver \bxi_{J}.
\end{align*}
Fourth, using the above commutativity, 
	\begin{align*}
	\bv_{J} & = \nabla \times \bxi_{J} = \underbrace{\nabla \times W^1 \bxi_{J}}_{\qe \bdelta_{J,0}}
	+ \sum_{\ver \in \V_{J}} \underbrace{\nabla \times B^1_\ver \bxi_{J}}_{\qe \bdelta_{J,\ver}} \\
    & = W^2 \nabla \times \bxi_{J} 
	+ \sum_{\ver \in \V_{J}} B^2_\ver \nabla \times \bxi_{J} \\
    & = W^2 \bv_{J} 
	+ \sum_{\ver \in \V_{J}} B^2_\ver \bv_{J}, 
\end{align*}
which is the requested decomposition~\eqref{OL_3dSD_1}. Fifth and last, by~\cite[Proposition~8.1]{Falk_Winther_25}, this decomposition is stable as
	\begin{align*}
	\big\|W^2 \bv_{J}\big\|^2 +
	\sum_{\ver \in \V_{J}}\big\|B^2_\ver \bv_{J}\big\|^2_{\omaJ}
	\le C_{\rm{SD}}^2\big\|\bv_{J}\big\|^2,
	\end{align*}
where the constant $C_{\rm{SD}}$ only depends on the
mesh-geometry parameter $\kappa_{\T}$. 
Thus~\eqref{OL_3dSD_2} follows, taking into account the variations of the diffusion coefficient $\dif$.
\end{proof}

\begin{proof}[Proof of Lemma~\ref{SD_assumption_d3_OL} (non $p$-robust, with $C_{\rm{OL}}$ possibly additionally depending on~$p$)]
The proof is presented in six steps and consists in initially introducing 
a lowest-order stable component of $\bv_J$ and constructing afterwards a local 
stable decomposition of the remaining high-order components. 
The latter decomposition relies on techniques following 
Chaumont-Frelet and Vohral{\'i}k~\cite[Theorem~B.1]{Chaum_Voh_Maxwell_equil_23}.

First, for ease of notation, we introduce the Raviart–Thomas 
interpolant of order $q\ge0$ given elementwise as follows: for any 
$ K \in \T_J$ and any $\bw \in [C^1(K)]^3$, construct 
$\bI^{\RT}_{K,q} (\bw) \in {\RT_q(K)}$ such that
\begin{subequations}\label{eq_I_RT}
	\begin{alignat}{2}
	(\bI^{\RT}_{K,q} (\bw){\cdot} \n_K, r_J)_F 
	& =	( \bw{\cdot} \n_K, r_J)_F  
	\qquad & & \forall r_J \in \mathbb{P}_q(F), \ \forall F \in \mathcal{F}_K, 
	\label{eq_interpF} \\
	(\bI^{\RT}_{K,q} (\bw), \br_J)_K           
	& =	( \bw, \br_J)_K 
	& & \forall \br_J \in [\mathbb{P}_{q-1}(K)]^3,   \label{eq_interpK}
	\end{alignat}
\end{subequations}
where $\mathcal{F}_K$ denotes the faces of the simplex $K$.
Moreover, we denote by $\Pi_{K,q}$ the
elementwise $L^2(K)$-orthogonal projection onto $\PP_{q}(K)$ 
for all $K \in \T_{J}$, $q>0$. Similarly, let $\Pi_{q}$ denote the 
$L^2(\Omega)$-orthogonal projection onto
$\PP_{q}(\T_J)$, $q>0$.
Finally, let the elementwise $[L^2(K)]^d$-orthogonal projection 
$ {\bm \Pi}_{K,q} $ be given componentwise by $\Pi_{K,q}$, 
for all $K \in \T_{J}$ and $q>0$.
Using, e.g., \cite[Proposition~2.5.2]{Boffi_Brezzi_Fortin_book}, 
the following elementwise commuting property holds for $\bI^{\RT}_{K,q}$: 
for any $\bw \in [C^1(K)]^3$
\begin{align}\label{eq_commuting}
	\nabla {\cdot} \bI^{\RT}_{K,q} (\bw)
	= \Pi_{K,q} (\nabla {\cdot} \bw).
\end{align}

\textbf{Step 1: Construction of the divergence-free lowest-order component.}
Let $\psi_{J,\ver}$ denote the piecewise affine hat function that takes 
value one in vertex $\ver \in \V_{J}$ and zero in all other vertices of 
$\V_{J}$. Define, for all vertices $\ver \in \V_J$ and all elements 
$K \in \T_J$ sharing $\ver$,
\begin{subequations}
	\begin{alignat}{2}
	\bdelta_{J,0}^{\ver}|_K 
	&:= \bI^{\RT}_{K,0} ((\psi_{J,\ver} \bv_J){|_K}), 
	\label{eq_alpha_loc} \\
	\bdelta_{J,0} &:= \sum_{\ver \in \V_J}  \bdelta_{J,0}^{\ver}.
	\label{eq_alpha_glob}
	\end{alignat}
\end{subequations}
Since $\bv_J \in \bV_J^0$, i.e. $\bv_J \in {\RT}_p(\T_J) 
\cap \H_{0}({\div}; \Omega)$ and $\nabla {\cdot} \bv_J = 0$, there holds
$\bdelta_{J,0}^{\ver} \in {\RT}_0(\T_{J}^{\ver}) \cap \H_{0}({\div};\omaJ)$, 
leading to $\bdelta_{J,0} \in {\RT}_0(\T_{J}) \cap \H_{0}({\div}; \Omega)$. 
Next, we use the commuting property~\eqref{eq_commuting}, the partition of 
unity for the hat functions $\sum_{\ver \in \V_J} \psi_{J,\ver} =1$, and
$\nabla {\cdot} \bv_J = 0$, to obtain
\begin{align}\label{divfree_alpha}
	\begin{split}
	\nabla {\cdot} \bdelta_{J,0}
	& \refvs{\eqref{eq_alpha_glob}}{=}{xx}    
	\sum_{\ver \in \V_J} \nabla {\cdot}  \bdelta_{J,0}^{\ver}
	\refvs{\eqref{eq_alpha_loc}}{=}{xx}
	\sum_{\ver \in \V_J} \sum_{K \in \T_{J}^{\ver}} 
	\nabla {\cdot}  \bI^{\RT}_{K,0} ((\psi_{J,\ver} \bv_J){|_K})
	\refvs{\eqref{eq_commuting}}{=}{xx}
	\sum_{\ver \in \V_J} \sum_{K \in \T_{J}^{\ver}}  
	\Pi_{K,0} (\nabla {\cdot} (\psi_{J,\ver} \bv_J)|_K)\\
	&\refvs{}{=}{xx}
	\sum_{\ver \in \V_J} \sum_{K \in \T_{J}^{\ver}}  
	\Pi_{K,0} (\nabla \psi_{J,\ver} {\cdot} \bv_J 
	+ \psi_{J,\ver} \nabla {\cdot} \bv_J )|_K
	= \Pi_{0} (\nabla 1 {\cdot} \bv_J + \nabla {\cdot} \bv_J ) =0.
	\end{split}
\end{align}
This means that 
$\bdelta_{J,0}  \in \bV_{J}^0 \cap {\RT}_0 (\T_{J})$ as requested.

\textbf{Step 2: Stability of the divergence-free lowest-order component.}
The stability of the lowest-order divergence-free component 
with respect to~$\bv_J$ follows as
\begin{align}\label{eq_locstab_LO}
	\begin{split}
	\| \bdelta_{J,0}  \|_K^2
	&\refvs{\eqref{eq_alpha_glob}}{=}{xx}
	\Big\| \sum_{\ver \in \V_K}  \bdelta_{J,0}^{\ver}   \Big\|_K^2
	\le (d+1)  \sum_{\ver \in \V_K}
	\|  \bdelta_{J,0}^{\ver}   \|_K^2
	\le (d+1)  \sum_{\ver \in \V_K}
	\|  \bdelta_{J,0}^{\ver}   \|_{\omaJ}^2 \\
	&\refvs{\eqref{eq_alpha_loc}}{=}{xx}
	(d+1)  \sum_{\ver \in \V_K} \sum_{K \in \T_{J}^{\ver}}
	\|  \bI^{\RT}_{K,0} ((\psi_{J,\ver} \bv_J){|_K})\|_{K}^2
	\le C_0 (d+1) \sum_{\ver \in \V_K}
	\|  \psi_{J,\ver} \bv_J\|_{\omaJ}^2
	\\
	&   \le C_0 (d+1) \sum_{\ver \in \V_K}
	\|  \psi_{J,\ver} \|^{2}_{\infty,\omaJ} \| \bv_J\|_{\omaJ}^2
	= C_0 (d+1) \sum_{\ver \in \V_K}
	\|  \bv_J\|_{\omaJ}^2,
	\end{split}
\end{align}
where the constant $C_0>0$ arises from the stability of the elementwise 
Raviart--Thomas interpolant defined in~\eqref{eq_I_RT}, see e.g., 
\cite[Chapter~16]{Ern_Guermond_FEM_book_pt1}. Summing over the 
elements yields the desired bound
\begin{align}\label{eq_stab_LO}
	\| \bdelta_{J,0} \|^2
	=
	\sum_{K \in \T_J} \| \bdelta_{J,0}  \|_K^2
	\refvs{\eqref{eq_locstab_LO}}{\le}{xx}
	\sum_{K \in \T_J} C_0 (d+1) \sum_{\ver \in \V_K}
	\|  \bv_J\|_{\omaJ}^2
	\le C_{\kappa_{\T}} C_0 (d+1)
	\|  \bv_J\|^2,
\end{align}
where $C_{\kappa_{\T}} > 0$ only depends on the mesh shape-regularity parameter 
$\kappa_{\T}$.

\textbf{Step 3: 
Construction of the high-order divergence-free local components.} 
In the following steps we pursue the approach of 
Chaumont-Frelet and Vohral{\'i}k~\cite[Theorem~B.1]{Chaum_Voh_Maxwell_equil_23} 
to construct local stable high-order contributions of $\bv_J$. 
Consider the following elementwise construction for all vertices 
$\ver \in \V_J$ and all elements $K \in \T_J$ sharing $\ver$:
\begin{align}\label{delta_ver_constr}
	\bdelta_{J,\ver}|_K := \underset{
	\substack{  
	\bw_{J}  \in {\RT_p(K)}
	\\
	\nabla {\cdot} \bw_{J} = 0
	\\
	\bw_{J} {\cdot} {\bm n} = 
	\big( \bI^{\RT}_{K,p} 
	((\psi_{J,\ver} \bv_J){|_K}) 
	- \bdelta_{J,0}^{\ver} \big) {\cdot} {\bm n} 
	\text{ on } \partial K
	}
	}{\rm arg \ min}
	\big\| \big( \bI^{\RT}_{K,p} ((\psi_{J,\ver} \bv_J){|_K}) 
	- \bdelta_{J,0}^{\ver} \big) - \bw_{J}  \big \|_K.
\end{align}
The well-posedness of such a problem is equivalent to the Neumann 
compatibility condition 
\linebreak$\big( \big( \bI^{\RT}_{K,p} ((\psi_{J,\ver} \bv_J){|_K}) 
- \bdelta_{J,0}^{\ver} \big) {\cdot} {\bm n}, 1 \big)_{\partial K} = 0 $ 
being satisfied,
see e.g. Boffi, Brezzi, and Fortin~\cite{Boffi_Brezzi_Fortin_book}. 
This condition is indeed satisfied, thanks to the definition of
$\bdelta_{J,0}$ and the Raviart--Thomas interpolation using that 
$1 \in \mathbb{P}_0 (\partial K)$
\begin{align*}
	\big( \big( \bI^{\RT}_{K,p} ((\psi_{J,\ver} \bv_J){|_K}) 
	- \bdelta_{J,0}^{\ver} \big) {\cdot} {\bm n}, 1 \big)_{\partial K}
	&\refvs{\eqref{eq_alpha_loc}}{=}{xx}
	\big( \big( \bI^{\RT}_{K,p} ((\psi_{J,\ver} \bv_J){|_K}) 
	-  \bI^{\RT}_{K,0} ((\psi_{J,\ver} \bv_J){|_K})  \big) 
	{\cdot} {\bm n}, 1 \big)_{\partial K}\\
	&
	\refvs{\eqref{eq_interpF}}{=}{xx}
	\big(   (\psi_{J,\ver}  \bv_J )|_K{\cdot} {\bm n}, 1 \big)_{\partial K}
	- \big(  (\psi_{J,\ver}  \bv_J)|_K{\cdot} {\bm n}, 1 \big)_{\partial K} 
	= 0.
\end{align*}
By construction of the normal components and the divergence constraint 
in~\eqref{delta_ver_constr}, the local high-order contributions belong indeed 
to the local spaces $\bV_J^{\ver,0}$ defined by~\eqref{V_ja_0} for all 
$\ver \in \V_J$, i.e.,
\[
\bdelta_{J,\ver} \in {\RT}_p(\T_{J}^{\ver}) \cap \H_{0}({\div}; \omaJ) 
\quad \text{and} \quad
\nabla {\cdot} \bdelta_{J,\ver} = 0.
\]

\textbf{Step 4: Stability of the high-order divergence-free local components.} 
To show stability of the high-order divergence-free local components 
with respect to $\bv_J$, we first introduce, for all $\ver \in \V_J$ 
and all $K \in \T_J^{\ver}$, auxiliary elementwise constructions
\begin{align}\label{delta_ver_aux}
	\widehat{\bdelta}_{J,\ver}|_K :=
	\underset{\substack{
	\bw_J  \in {\RT_p(K)}
	\\
	\nabla {\cdot} \bw_J =  
	\nabla {\cdot} \big( \bI^{\RT}_{K,p} 
	((\psi_{J,\ver} \bv_J){|_K}) - \bdelta_{J,0}^{\ver} \big)
	\\
	\bw_J {\cdot} {\bm n} = 0 \text{ on } \partial K
	}}{\rm arg \ min}
	\|   \bw_J  \|_K.
\end{align}
These problems are also well-posed as the Neumann compatibility condition holds,
\begin{equation} \label{eq_NC} 
	\begin{split}
	\big(  \nabla {\cdot}  \big( \bI^{\RT}_{K,p} 
	((\psi_{J,\ver} \bv_J){|_K}) - \bdelta_{J,0}^{\ver} \big), 1 \big)_K
	&\refvs{\eqref{eq_alpha_loc}}{=}{xxx}
	\big(  \nabla {\cdot} \big( \bI^{\RT}_{K,p} ((\psi_{J,\ver} \bv_J){|_K}) 
	- \bI^{\RT}_{K,0} ((\psi_{J,\ver} \bv_J){|_K})\big), 1 \big)_K
	\\ &
	\refvs{\eqref{eq_commuting}}{=}{xxx}
	\big(   \Pi_{K,p} (\nabla {\cdot}   (\psi_{J,\ver}  \bv_J)|_K) 
	-   \Pi_{K,0} (\nabla {\cdot} (\psi_{J,\ver} \bv_J)|_K) , 1 \big)_K =0,
	\end{split}
\end{equation}
where we have used the definition of the elementwise $L^2(K)$-orthogonal 
projection. Immediately from the
definitions in \eqref{delta_ver_constr} and~\eqref{delta_ver_aux}, we obtain 
\begin{align}\label{delta_hat_rewrite}
	{\widehat	\bdelta}_{J,\ver}|_K  
	= \bI^{\RT}_{K,p} ((\psi_{J,\ver} \bv_J){|_K}) 
	- \bdelta_{J,0}^{\ver}|_K \,
	- \,	\bdelta_{J,\ver}|_K .
\end{align}
Following Ern and Vohral{\'i}k~\cite[Lemma~A.3]{Ern_Voh_p_rob_3D_20}, 
one obtains for $C_{\rm stab} >0$ only depending on mesh shape regularity 
the stability result
\begin{align}\label{cont_pb_stab}
	\underset{\substack{
	\bw_J  \in {\RT_p(K)}
	\\
	\nabla {\cdot} \bw_J =  
	\nabla {\cdot} \big( \bI^{\RT}_{K,p} ((\psi_{J,\ver} \bv_J){|_K}) 
	- \bdelta_{J,0}^{\ver} \big)
	\\
	\bw_J {\cdot} {\bm n} = 0 \text{ on } \partial K
	}}{\rm  min}
	\|   \bw_J  \|_K
	\le C_{\rm stab}
	\underset{\substack{
	\bv  \in \H ({\div}, K)
	\\
	\nabla {\cdot} \bv =  \nabla {\cdot} \big( \bI^{\RT}_{K,p} 
	((\psi_{J,\ver} \bv_J){|_K}) - \bdelta_{J,0}^{\ver} \big)
	\\
	\bv {\cdot} {\bm n} = 0 \text{ on } \partial K
	}}{\rm  min}
	\|   \bv  \|_K.
\end{align}
Thus, we have
\begin{align*}
	\big\| \bdelta_{J,\ver} - \big( \bI^{\RT}_{K,p} 
	((\psi_{J,\ver} \bv_J){|_K}) 
	- \bdelta_{J,0}^{\ver} \big)\big\|_K
	\refvs{\eqref{delta_hat_rewrite}}{=}{xxx}
	\|	{\widehat	\bdelta}_{J,\ver} \|_K
	&\refvs{\eqref{delta_ver_aux}}{=}{xxx}
	\underset{\substack{
	\bw_J  \in {\RT_p(K)}
	\\
	\nabla {\cdot} \bw_J 
	=  \nabla {\cdot} \big( \bI^{\RT}_{K,p} ((\psi_{J,\ver} \bv_J){|_K}) 
	- \bdelta_{J,0}^{\ver} \big)
	\\
	\bw_J {\cdot} {\bm n} = 0 \text{ on } \partial K
	}}{\rm  min}
	\|   \bw_J  \|_K
	\\ &
	\refvs{\eqref{cont_pb_stab}}{\le}{xxx}
	C_{\rm stab}
	\underset{\substack{
	\bv  \in \H ({\div}, K)
	\\
	\nabla {\cdot} \bv =  \nabla {\cdot} \big( \bI^{\RT}_{K,p} 
	((\psi_{J,\ver} \bv_J){|_K}) - \bdelta_{J,0}^{\ver} \big)
	\\
	\bv {\cdot} {\bm n} = 0 \text{ on } \partial K
	}}{\rm  min}
	\|   \bv  \|_K \\
	& =  C_{\rm stab} \| \nabla \zeta_K \|_K,
\end{align*}
where $\zeta_K \in H^1_\star (K) := \{ v \in H^1(K), (v,1)_K = 0 \}$ 
solves the equivalent primal problem
\begin{align}\label{primal_aux}
	(\nabla \zeta_K, \nabla v)_K 
	= (\nabla {\cdot} \big( \bI^{\RT}_{K,p} ((\psi_{J,\ver} \bv_J){|_K}) 
	- \bdelta_{J,0}^{\ver} \big) ,v)_{K} \quad \forall v \in H^1_\star (K).
\end{align}

Next, similarly as done in~\eqref{eq_NC}, we use 
\eqref{eq_alpha_loc} and \eqref{eq_commuting} to obtain
\begin{equation}\label{eq_div}
	\begin{split}
	\nabla {\cdot} \big( \bI^{\RT}_{K,p} ((\psi_{J,\ver} \bv_J){|_K}) 
	- \bdelta_{J,0}^{\ver} \big) 
	& = \Pi_{K,p} (\nabla {\cdot}   (\psi_{J,\ver}  \bv_J)|_K) 
	-   \Pi_{K,0} (\nabla {\cdot} (\psi_{J,\ver} \bv_J)|_K) \\
	& = \Pi_{K,p} ((\nabla \psi_{J,\ver} {\cdot} \bv_J +  0 )|_K) 
	- \Pi_{K,0} ((\nabla \psi_{J,\ver} {\cdot} \bv_J +  0 )|_K)\\
	& = \nabla \psi_{J,\ver}|_K {\cdot}   
	(\bv_J|_K  - {\bm \Pi}_{K,0} (\bv_J|_K)),
	\end{split}
\end{equation}
since $\bv_{J}$ is divergence-free and thus actually $p$-degree polynomial in 
each component, see, e.g., \cite[Corollary~2.3.1]{Boffi_Brezzi_Fortin_book}. 
Consequently, we can pursue our main estimate using the Cauchy--Schwarz and the 
Poincaré inequality
\begin{align}\label{shifted_stab}
	\begin{split}
	\big\| \bdelta_{J,\ver} - \big( \bI^{\RT}_{K,p} 
	((\psi_{J,\ver} \bv_J){|_K}) &- \bdelta_{J,0}^{\ver} \big)\big\|_K
	\lesssim
	\| \nabla \zeta_K \|_K
	= \underset{\substack{
	v \in H^1_\star (K)
	\\
	\| \nabla v \|_K = 1
	}}{\rm  max} (\nabla \zeta_K, \nabla v)_K
	\\&
	\refvs{\substack{\eqref{primal_aux}\\\eqref{eq_div}}}{=}{xx}
	\underset{\substack{
	v \in H^1_\star (K)
	\\
	\| \nabla v \|_K = 1
	}}{\rm  max} (\nabla \psi_{J,\ver} {\cdot}   
	(\bv_J  - {\bm \Pi}_{K,0} (\bv_J|_K)),  v)_K
	\\ &
	\refvs{}{\le}{xx}
	\underset{\substack{
	v \in H^1_\star (K)
	\\
	\| \nabla v \|_K = 1
	}}{\rm  max} \Big\{ \|\nabla \psi_{J,\ver}\|_{K,\infty} \|\bv_J  
	- {\bm \Pi}_{K,0} (\bv_J|_K)\|_K \frac{h_K}{\pi}  \| \nabla v \|_K \Big\}
	\\
	&
	\refvs{}{=}{xx} \|\nabla \psi_{J,\ver}\|_{K,\infty} \|\bv_J  
	- {\bm \Pi}_{K,0} (\bv_J|_K)\|_K \frac{h_K}{\pi}\\
	& \refvs{}{\lesssim}{xx} \|\bv_J  - {\bm \Pi}_{K,0} (\bv_J|_K)\|_K
	\leq \|\bv_J\|_K,
	\end{split}
\end{align}
where the last but one estimate follows as 
$\|\nabla \psi_{J,\ver}\|_{K,\infty} \lesssim h_K^{-1}$, 
where the hidden constant only depends on the mesh shape regularity, 
and the last estimate from the stability of the $[L^2(K)]^d$-orthogonal 
projection.
This finally leads to the desired stability
\begin{align}\label{delta_Ja_stab}
	\begin{split}
	\| \bdelta_{J,\ver} \|_K
	&\le
	\big\| \bdelta_{J,\ver} - \big( \bI^{\RT}_{K,p} 
	((\psi_{J,\ver} \bv_J){|_K}) - \bdelta_{J,0}^{\ver} \big)\big\|_K
	+  \big\|  \bI^{\RT}_{K,p} ((\psi_{J,\ver} \bv_J){|_K}) 
	- \bdelta_{J,0}^{\ver} \big\|_K
	\\ &
	\refvs{\substack{\eqref{shifted_stab}\\\eqref{eq_alpha_loc}}}{\lesssim}{xxx}
	\|\bv_J\|_K + \big\|  \bI^{\RT}_{K,p} 
	((\psi_{J,\ver} \bv_J){|_K}) \big\|_K  
	+ \| \bI^{\RT}_{K,0} ((\psi_{J,\ver} \bv_J){|_K}) \|_K
	\lesssim
	3 \|\bv_J\|_K,
	\end{split}
\end{align}
where we have also used the fact that $\|\psi_{J,\ver}\|_{K,\infty} = 1$ and 
the stability of the elementwise Raviart--Thomas interpolant defined 
in~\eqref{eq_I_RT}, which brings the dependence of the hidden constant on the 
polynomial degree~$p$.

\textbf{Step 5: Decomposition property~\eqref{OL_3dSD_1}.} It remains to show 
that the sum of the contributions we introduced indeed gives the original 
$\bv_J \in \bV_J^0$, i.e., the equality  
\begin{equation}\label{eq_dec}
	\sum_{\ver \in \V_J} \bdelta_{J, \ver} = \bv_J - \bdelta_{J,0},
\end{equation}
where $\bdelta_{J,0}$ is the global lowest-order 
contribution~\eqref{eq_alpha_glob}
and $\bdelta_{J, \ver} $ for $\ver \in \V_J$ are the local high-order 
contributions~\eqref{delta_ver_constr}.
First, note that for all mesh elements $K \in \T_J$, there holds
\begin{align}\label{decomp_a}
	\begin{split}
	\Big( \sum_{\ver \in \V_J} \bdelta_{J, \ver} \Big) 
	|_{\partial K} {\cdot} {\bm n}
	& \refvs{\eqref{delta_ver_constr}}{=}{xx}
	\sum_{\ver \in \V_J} \big( \bI^{\RT}_{K,p} ((\psi_{J,\ver} 
	\bv_J){|_{ \partial K}}) 
	- \bdelta_{J,0}^{\ver} \big)|_{\partial K} {\cdot} {\bm n} \\
	& \refvs{\eqref{eq_alpha_glob}}{=}{xx}
	\big(\bI^{\RT}_{K,p} (\bv_J|_K) - \bdelta_{J,0} 
	\big)|_{\partial K} {\cdot} {\bm n} 
	= ( \bv_J - \bdelta_{J,0} )|_{\partial K} {\cdot} {\bm n},
	\end{split}
\end{align}
also using that $\bI^{\RT}_{K,p}$ from~\eqref{eq_I_RT} is linear and 
projector. Thus, the functions on both sides of~\eqref{eq_dec} 
have the same normal components on all mesh faces.
Next, using that $\bv_J$ is divergence-free, we obtain
\begin{align}\label{decomp_b}
	\begin{split}
	\nabla {\cdot} 	\Big( \sum_{\ver \in \V_J} \bdelta_{J, \ver} \Big)
	=	 	
	\sum_{\ver \in \V_J} \nabla {\cdot} \bdelta_{J, \ver}
	\refvs{\eqref{delta_ver_constr}}{=}{xx} 0
	\refvs{\eqref{divfree_alpha}}{=}{xx}
	\nabla {\cdot} (\bv_J - \bdelta_{J,0} ).
	\end{split}
\end{align}
Thus, the functions on both sides of~\eqref{eq_dec} have the same divergence.

Thanks to~\eqref{delta_ver_constr}, for all $K \in \T_J$, there holds
\[
\big( \bdelta_{J,\ver} 
-  \big(\bI^{\RT}_{K,p} ((\psi_{J,\ver} \bv_J){|_K}) 
- \bdelta_{J,0}^{\ver} \big), \bw_J \big)_{K} =0
\quad \forall \bw_{J}  \in {\RT_p(K)},
\nabla {\cdot} \bw_{J} = 0, \bw_{J} {\cdot} {\bm n} = 0.
\]
Summing over all vertices, together with the patchwise 
construction~\eqref{eq_alpha_glob} and the fact that the hat functions form a 
partition of unity, i.e., $\sum_{\ver \in \V_J} \psi_{J,\ver}  = 1$, leads to
\begin{align}\label{decomp_c}
	\begin{split}
	0 & = \Big( \sum_{\ver \in \V_J} \big( \bdelta_{J,\ver} 
	-  \big(\bI^{\RT}_{K,p} ((\psi_{J,\ver} \bv_J){|_K}) 
	- \bdelta_{J,0}^{\ver} \big)\big), \bw_J \Big)_{K}\\
	&= \Big( \sum_{\ver \in \V_J} \bdelta_{J,\ver} 
	-  \big(\bI^{\RT}_{K,p} ( \bv_J|_K) - \bdelta_{J,0} \big), 
	\bw_J \Big)_{K} \\
	&= \Big( \sum_{\ver \in \V_J} \bdelta_{J,\ver} 
	-  \big( \bv_J - \bdelta_{J,0} \big), \bw_J \Big)_{K}
	\quad \forall \bw_{J}  \in {\RT_p(K)},
	\nabla {\cdot} \bw_{J} = 0, \bw_{J} {\cdot} {\bm n} = 0.
	\end{split}
\end{align}
Thus, the difference of the functions on both sides of~\eqref{eq_dec} is 
elementwise orthogonal to normal-component-free and divergence-free functions. 
Consequently, \eqref{decomp_a}--\eqref{decomp_c} together yield the 
desired decomposition~\eqref{OL_3dSD_1}.

\textbf{Step 6: Stability property~\eqref{OL_3dSD_2}.}
Using the results of steps~2 and~4, one readily obtains
\[
\| \bdelta_{J,0} \|^2 + \sum_{\ver \in \V_J} 
\| \bdelta_{J, \ver} \|_{\omaJ}^2
\refvs{\eqref{eq_stab_LO}}{\lesssim}{xx}
\| \bv_{J} \|^2 + \sum_{\ver \in \V_J}  \sum_{K \in \T_J} 
\| \bdelta_{J, \ver} \|_{K}^2
\refvs{\eqref{delta_Ja_stab}}{\lesssim}{xx}
\| \bv_{J} \|^2.
\]
Therefrom, the desired stability property~\eqref{OL_3dSD_2} 
follows taking into account the variations of the diffusion coefficient $\dif$.

\end{proof}

\subsection{Multilevel stable decomposition for lowest-order 
divergence-free Raviart--Thomas piecewise polynomials}

After having obtained a one-level high-order stable decomposition in 
Lemma~\ref{SD_assumption_d3_OL}, we decompose further its global 
lowest-order contribution in a stable multilevel way. Recall the notations 
$\bV^0_j$ from~\eqref{V_j_0} and $\bVajz$ from~\eqref{V_ja_0}.

\begin{lemma}[Multilevel lowest-order stable decomposition of $\bV_{J}^0$]\label{SD_assumption_d3_ML}
Let $d=3$.
Let 
either Assumption~\ref{assumption_refinement_quasiuniformity} or 
Assumption~\ref{graded_grids} hold. 
Let the polynomial degree $p=0$. Then, for any $\bv_{J} \in \bV_{J}^0$, 
there exist $\bv_0 \in \bV_0^0$ (global coarse-grid component) and 
$\bv_{j,\ver} \in \bVajz$, $\ver \in \V_j$, 
(local levelwise and patchwise components) such that
\begin{subequations}\label{ML_3dSD}
\begin{align}
	\bv_{J} 
	= \bv_0 + \sum_{j = 1}^J \sum_{\ver \in \V_j} \bv_{j,\ver}\label{ML_3dSD_1}
\end{align}
stable as
\begin{align}
	\big\|\dif^{-1/2}\bv_0\big\|^2 +
	\sum_{j = 1}^J \sum_{\ver \in \V_j} 
	\big\|\dif^{-1/2}\bv_{j,\ver}\big\|^2_{\omaj}
	\le C_{\rm{ML}}^2\big\|\dif^{-1/2}\bv_{J}\big\|^2,\label{ML_3dSD_2}
\end{align}
\end{subequations}
where the constant $C_{\rm{ML}}$ only depends on
mesh-geometry parameters $\kappa_{\T}$, $ C_{\rm qu}$ and $C_{\rm ref}$ or 
$C^0_{\rm qu}$ and $C_{\rm loc,qu}$, the diffusion inhomogeneity or 
anisotropy ratio $\Lambda_{\max}/\Lambda_{\min}$, and 
the domain $\Omega$. 
In particular, $C_{\rm{ML}}$ is independent of the number of mesh levels $J$.
\end{lemma}

\begin{proof}
First, using that $\Omega$ is simply connected, we introduce the stable 
discrete vector potential $\bxi_{J} \in \bX_{J,0}$ (recall~\eqref{eq_Ned_0}) of 
$\bv_{J}$ using Lemma~\ref{lem_stabpot}. Then, we decompose $\bxi_{J}$ 
following 
Hiptmair, Wu, and Zheng~\cite[Theorem~5.2]{Hipt_Wu_Zheng_cvg_adpt_MG_12}, 
which will lead to a stable decomposition of the original $\bv_{J}$.

We first introduce some notation. Let $\mathcal{E}_j^+$ denote the set of new 
edges arising in the $j$-th mesh refinement step and their immediate neighboring 
edges, i.e., edges whose edge-patch has shrunk in the refinement step. 
Let $\bphi_{j,e} $ be the lowest-order basis function of the 
N\'ed\'elec space $\bX_{j,0}$ associated to the edge $e \in \mathcal{E}_j^+$ and 
let ${\omega_j^e} := {\rm supp}(\bphi_{j,e})$. Similarly, let $\V_j^+$ denote 
the set of new vertices arising in the $j$-th mesh refinement step and their 
immediate neighboring vertices. Denote by $\psi_{j,a}$ the piecewise linear hat 
function which takes value one at the vertex $\ver \in \V_j$ and vanishes in all 
other vertices of level $j$. Then, 
from~\cite[Theorem~5.2]{Hipt_Wu_Zheng_cvg_adpt_MG_12}, for any 
$\bxi_{J} \in \bX_{J,0}$, there exist $\bxi_0 \in \bX_{0,0}$, 
$\bxi_{j,e} \in {\rm span} \{ \bphi_{j,e} \}$ with $e \in \mathcal{E}_j^+$, 
$\xi_{j,a} \in {\rm span} \{ \psi_{j,a} \}$ with $\ver \in \V_j^+$, 
$1 \le j \le J$, such that the decomposition
\begin{subequations}
	\begin{align}
	\bxi_{J} = \bxi_0 + \sum_{j=1}^J \sum_{e \in \mathcal{E}_j^+} \bxi_{j,e}
	+ \sum_{j=1}^J \sum_{\ver \in \V_j^+} \nabla \xi_{j,a} \label{hcurl_decomp}
	\end{align}
is stable as
	\begin{align}
	\| \bxi_0 \|^2_{\H(\curl; \Omega)}
	+ \sum_{j=1}^J \sum_{e \in \mathcal{E}_j^+} 
	\| \bxi_{j,e} \|^2_{\H(\curl; {\omega_j^e})}
	+ \sum_{j=1}^J \sum_{\ver \in \V_j^+} \| \nabla \xi_{j,a} \|^2_{\omaj}
	\le C_{\rm 3D}^2\|\bxi_{J}  \|^2_{\H(\curl; \Omega)},\label{hcurl_stab}
	\end{align}
\end{subequations}
where 
$ \| \bzeta \|^2_{\H(\curl; \omega)} 
:=  \| \bzeta \|^2_{\omega} + \| \nabla \times \bzeta \|^2_{\omega} $ 
for any $\bzeta \in \H(\curl; \omega) $, $\omega \subset \mathbb{R}^d$.

Let the edges in $\mathcal{E}_j^+$ be re-indexed by associating them to 
patches such that any edge is only counted once and such that each edge basis 
function has tangential trace zero on the boundary of the patch it is 
associated to. We write 
$\mathcal{E}_j^+ = \cup_{\ver \in \V_j} \mathcal{E}_j^\ver$, which gives 
$\sum_{\ver \in \V_j} \sum_{e \in \mathcal{E}_j^\ver}  
= \sum_{e \in \mathcal{E}_j^+} $.
Then, we apply the decomposition~\eqref{hcurl_decomp} to~\eqref{vh_curl}, 
leading to
\begin{align*}
	\bv_{J} = \nabla \times \bxi_{J} = \nabla \times \bxi_0
	+ \sum_{j=1}^J \sum_{e \in \mathcal{E}_j^+} \nabla \times  \bxi_{j,e}
	+ 0
	= \bv_0
	+ \sum_{j=1}^J \sum_{\ver \in \V_j} \bv_{j,a},
\end{align*}
with $\bv_0 := \nabla \times \bxi_0 \in \bV_{0}^{0}$ and 
$\bv_{j,\ver} \eq \sum_{e \in \mathcal{E}_j^\ver} \nabla \times \bxi_{j,e} 
\in \bV_{j}^{\ver,0}$.
Finally, this decomposition is stable as
\begin{align*}
	\| \bv_0 \|^2
	+   \sum_{j=1}^J \sum_{\ver \in \V_j}  \|\bv_{j,\ver}\|^2_{\omaj}
	& =
	\| \nabla \times \bxi_0 \|^2
	+  \sum_{j=1}^J \sum_{\ver \in \V_j}
	\big\|  \sum_{e \in \mathcal{E}_j^\ver}   
	\nabla \times \bxi_{j,e}\big\|^2_{\omaj} \\
	& \lesssim
	\| \nabla \times \bxi_0 \|^2
	+  \sum_{j=1}^J
	\sum_{e \in \mathcal{E}_j^+}
	\|\nabla \times \bxi_{j,e}\|^2_{\omega_j^e}\\
	& \stackrel{\eqref{hcurl_stab}}\le C_{\rm 3D}^2 \big(\| \bxi_{J} \|^2 
	+ \| \nabla \times \bxi_{J} \|^2 \big)
	\stackrel{\eqref{vh_curl}}\le
	C_{\rm 3D}^2 (C_{\rm pot}^2 +1) \|\bv_{J} \|^2.
\end{align*}
The desired result~\eqref{ML_3dSD_2} is then obtained after taking into 
consideration the variations of the diffusion coefficient $\dif$.
\end{proof}

\subsection{Proof of Proposition~\ref{SD_assumption} for $d=3$}

We finally prove Proposition~\ref{SD_assumption} for~$d=3$. Recall that the proof for $d=2$ was given in Section~\ref{sec_MLD_2D}.

\begin{proof}[Proof of Proposition~\ref{SD_assumption} for $d=3$]
Let $\bv_{J} \in \bV_{J}^0$. We first apply the one-level 
high-order decomposition of Lemma~\ref{SD_assumption_d3_OL}. 
This gives us $\balpha_{J} \in \bV_{J}^0 \cap {\RT}_0 (\T_{J}) $ and
$\bdelta_{J,\ver} \in \bV_{J}^{\ver,0}$, $\ver \in \V_{J}$, such that
$$\bv_{J} = \balpha_{J} +
\sum_{\ver \in \V_{J}} \bdelta_{J,\ver}$$
with the stability bound~\eqref{OL_3dSD_2}.
Now, we decompose $\balpha_{J}$ thanks to Lemma~\ref{SD_assumption_d3_ML}. 
This gives $\balpha_0 \in \bV_0^0 \cap {\RT}_0 (\T_0) $ and
$\bdelta_{j,0}^{\ver} \in \bV_{j}^{\ver,0} \cap{\RT}_0 (\T_{\ver})$, 
$\ver \in \V_j$, such that
$$\balpha_{J} = \balpha_0 +
\sum_{j = 1}^J \sum_{\ver \in \V_j} \bdelta_{j,0}^{\ver}$$
and the stability bound~\eqref{ML_3dSD_2}.
Define $\bv_0 := \balpha_0 \in \bV_0^0 \cap \ {\RT}_0 (\T_0)$, 
$\bv_{j,\ver} := \bdelta_{j,0}^{\ver} \in \bV_{j}^{\ver,0} 
\cap{\RT}_0 (\T_{\ver})$ for $1\le j \le J-1$, and 
$\bv_{J,\ver} := \bdelta_{J,0}^{\ver}  + \bdelta_{J,\ver} \in \bV_{J}^{\ver,0}$.
This yields the decomposition~\eqref{SD_1}. 
The stability~\eqref{SD_2} follows by
\begin{align*}
	\big\|\dif^{-1/2}\bv_0\big\|^2
	& +
	\sum_{j = 1}^J 
	\sum_{\ver \in \V_j}\big\|\dif^{-1/2}\bv_{j,\ver}\big\|^2_{\omaj}
	\\
	& = \big\|\dif^{-1/2}\balpha_0\big\|^2 +
	\sum_{j = 1}^{J-1} 
	\sum_{\ver \in \V_j}
	\big\|\dif^{-1/2}\bdelta_{j,0}^{\ver}\big\|^2_{\omaj}
	+\sum_{\ver \in \V_J} 
	\big\|\dif^{-1/2}(\bdelta_{J,0}^{\ver} 
	+ \bdelta_{J,\ver}) \big\|^2_{\omega_J^{\ver}}
	\\
	& \stackrel{\eqref{ML_3dSD_2}}\le 
	2 C_{\rm{ML}}^2\big\|\dif^{-1/2}\balpha_{J}\big\|^2 
	+ 2 \sum_{\ver \in \V_{J}} 
	\big\|\dif^{-1/2}\bdelta_{J,\ver} \big\|^2_{\omega_{J}^{\ver}}
	\stackrel{\eqref{OL_3dSD_2}}\le 2 C_{\rm{OL}}^2
	C_{\rm{ML}}^2\big\|\dif^{-1/2}\bv_{J}\big\|^2.
\end{align*}
\end{proof}

	\begin{remark}
		Note that in \cite{Ewing_Wang_92,Ewing_Wang_94} a 
		stable decomposition of the form~\eqref{SD} was presented in 
		two space dimensions and for the case of quadrilateral and 
		triangular elements, respectively. In \cite{Cai_etal}, 
		the stable decomposition was extended to three space dimensions 
		for the case of tetrahedral elements and lowest-order polynomial 
		degree. 
		Thus, the result of Proposition~\ref{SD_assumption} is more 
		general as it holds for $d \in \{ 2,3\}$ and polynomial degree $p\ge0$.
        Most importantly, thanks to~\cite{SchMelPechZag_08} used in Lemma~\ref{lem_SD_Lag} and~\cite{Falk_Winther_25} used in the proof of Lemma~\ref{SD_assumption_d3_OL}, there is no dependence on the polynomial degree $p$ in~\eqref{SD_2}. 
	\end{remark}

\section{Proof of the reliability and efficiency of the algebraic error estimators of Theorem~\ref{thm_upper_bound}}\label{sec_proofs}

In this section, we present a proof of Theorem~\ref{thm_upper_bound}. 
We rely for this purpose on the stable decomposition result of 
Proposition~\ref{SD_assumption} established above. First, we present a remark:

\begin{remark}[Lower bound on the optimal step-sizes]
In the multilevel setting, as in \emph{\cite[Lem\-ma~10.1]{Mir_Pap_Voh_lambda}}, 
\eqref{rho_j} together with a patch overlap argument yields
\begin{equation}\label{overlapping_rhoa}
	\big\|\dif^{-1/2}\brho_j^i\big\|^2
	\le (d+1) \sum_{\ver \in \V_j} 
	\big\|\dif^{-1/2} \brho_{j,\ver}^i\big\|^2_{\omaj}.
\end{equation}
A direct consequence of~\eqref{overlapping_rhoa} and \eqref{brhoaj_lam}, 
together with the definition $\lambda_j^i = 1$ when $j=0$ or $\brho_j^i = 0$, 
is
\begin{equation}\label{lam_LB}
	\lambda_j^i \ge \frac{1}{d+1} \quad 0 \le j \le J.
\end{equation}
Similarly, in the domain decomposition setting, 
\eqref{rho_h_i_a} together with a patch overlap on the coarse mesh $\T_H$ 
yields
\begin{equation}\label{overlapping_rhoa_DD}
	\big\|\dif^{-1/2}\brho_{h}^i\big\|^2
	\le (d+1) \sum_{\ver \in \V_H} \big\|\dif^{-1/2} \rhoa \big\|^2_{\omaH}.
\end{equation}
Thus, from~\eqref{overlapping_rhoa_DD} and \eqref{brhoaj_lam_DD}, together with 
the definition $\lambda_h^i = 1$ when $\brho_{h}^i = 0$,
\begin{equation}\label{lam_LB_DD}
	\lambda_h^i \ge \frac{1}{d+1}.
\end{equation}

\end{remark}

\subsection{Multilevel setting}\label{ML_proofs}

We can now present the proof of Theorem~\ref{thm_upper_bound} 
for the estimator built in Algorithm~\ref{Definition_solver}.

\begin{proof}[Proof of Theorem~\ref{thm_upper_bound} (multilevel setting)]
The proof follows closely the proof of \cite[Theorem~6.6]{Mir_Pap_Voh_lambda}, 
thus only the main steps are presented here.  
First, note that~\eqref{lower_bound} follows directly from the 
Pythagorean formula of Theorem~\ref{thm_error_contr} by using that
 $\big\|\dif^{-1/2}(\u_{J} - \u_{J}^{i+1})\big\| \ge 0$. To 
 show~\eqref{upper_bound}, we use the 
stable decomposition \eqref{SD_1} applied to $\u_{J} - \u_{J}^i$,
\begin{equation} \label{eq_SD_err}
	\u_{J} - \u_{J}^i 
	= \bv_0 + \sum_{j=1}^J \sum_{\ver \in \V_{j}}  \bv_{j,\ver}.
\end{equation}
Since $\u_{J} - \u_{J}^i \in \bV_{J}^0$ and using~\eqref{mfe_div_free} 
as well as the construction of the multilevel solver, one obtains
\begin{align*}
	\big\|\dif^{-1/2}(\u_{J} - \u_{J}^i)\big\|^2
	\stackrel{\eqref{mfe_div_free}}= {} & 
	- (\dif^{-1} \u_{J}^i, \u_{J} - \u_{J}^i)
	\stackrel{\eqref{eq_SD_err}}= 
	- \Big(\dif^{-1} \u_{J}^i, \bv_0 
	+ \sum_{j=1}^J \sum_{\ver \in \V_{j}}  \bv_{j,\ver} \Big) \\
	\stackrel{\eqref{rho_0}}= {} &
	(\dif^{-1}\brho_0^i,\bv_0)
	- \sum_{j=1}^J \sum_{\ver \in \V_{j}} (\dif^{-1}  
	\u_{J}^i, \bv_{j,\ver})_{\omaj}\\
	\stackrel{\eqref{rho_j_compact}}= {} &
	(\dif^{-1}\brho_0^i,\bv_0)
	+ \sum_{j=1}^J
	\sum_{\ver \in \V_{j}}
	(\dif^{-1} \brho_{j,\ver}^i, \bv_{j,\ver})_{\omaj} \\
	{} &
	+ \sum_{j=1}^J
	\sum_{m=0}^{j-1} 
	\sum_{\ver \in \V_{j}} \lambda_m^i 
	(\dif^{-1} \brho_{m}^i,  \bv_{j,\ver})_{\omaj}.
\end{align*}
In the same spirit 
of~\cite[Proof of Theorem 6.6, page S138]{Mir_Pap_Voh_lambda}, 
one can use Young's inequalities, the lower bound on the step-sizes
\eqref{lam_LB}, and finite overlap of patches as 
$ \Big|\sum\limits_{n=1}^{d+1} a_n \Big|^2 
\le (d+1) \sum\limits_{n=1}^{d+1} |a_n |^2$,  to estimate
\begin{equation}\label{eq_MG_Y}\begin{split}
	\big\|\dif^{-1/2}(\u_{J} - \u_{J}^i)\big\|^2
	&
   \le
	 C_{\rm{SD}}^2  \big\|\dif^{-1/2}\brho_0^i\big\|^2 
     + \frac{1}{4 C_{\rm SD}^2} \big\|\dif^{-1/2}\bv_0\big\|^2\\
     & \quad
    + \sum_{j=1}^J
	\sum_{\ver \in \V_{j}} 
    	\big( C_{\rm{SD}}^2 \big\|\dif^{-1/2} \brho_{j,\ver}^i\big\|_{\omaj}^2 
        +\frac{1}{4 C_{\rm SD}^2} 
		 \big\|\dif^{-1/2} \bv_{j,\ver}\big\|_{\omaj}^2 \big) \\
	 &\quad + \sum_{j=1}^J 
	\sum_{\ver \in \V_{j}} \Big(
    C_{\rm{SD}}^2 \big\|\dif^{-1/2} \big(\sum_{m=0}^{j-1}  
	\lambda_m^i \brho_{m}^i \big)\big\|_{\omaj}^2
	+\frac{1}{4 C_{\rm SD}^2} \big\|\dif^{-1/2} \bv_{j,\ver}\big\|_{\omaj}^2
    \Big)\\
	 &  \! \! \! \stackrel{\eqref{lam_LB}}\le C_{\rm{SD}}^2 
	 \Big( \big\|\dif^{-1/2}\brho_0^i\big\|^2
     + (d+1)  \lambda_j^i 
     \sum_{j=1}^J 
	\sum_{\ver \in \V_{j}} 
        \big\|\dif^{-1/2} \brho_{j,\ver}^i\big\|_{\omaj}^2
    \\
    & \quad 
       + J (d+1)
        \big\|\dif^{-1/2} \big(\sum_{m=0}^{J}  
		\lambda_m^i \brho_{m}^i \big)\big\|^2
        \Big)
        \\
    & \quad  +\frac{1}{2 C_{\rm SD}^2}  
    \Big(
    \big\|\dif^{-1/2}\bv_0\big\|^2
	+ \sum_{j=1}^J 
	\sum_{\ver \in \V_{j}} \big\|\dif^{-1/2} \bv_{j,\ver}\big\|_{\omaj}^2
    \Big).
\end{split}\end{equation}
Recalling the local writing \eqref{loc_eta} of the estimator $\etalg$, 
its property \eqref{brhoaj_lam}, and the stability of the 
decomposition~\eqref{SD_2}, leads to
\begin{align*}
	\big\|\dif^{-1/2}(\u_{J} - \u_{J}^i)\big\|^2
	&
   \le 
   C_{\rm SD}^2(d+1) \big( (\etalg)^2 
   + J^2 \sum_{m=0}^{J}   (\lambda_m^i 
   \big\|\dif^{-1/2}  \brho_{m}^i \big\|)^2 \big) 
   + \frac{1}{2}\big\|\dif^{-1/2}(\u_{J} - \u_{J}^i)\big\|^2  \\
   &
   \stackrel{\eqref{def_eta}}\le 
   C_{\rm SD}^2(d+1)(J^2+1)(\etalg)^2 
   + \frac{1}{2}\big\|\dif^{-1/2}(\u_{J} - \u_{J}^i)\big\|^2 . 
\end{align*}
Finally, a rewriting allows us to obtain~\eqref{upper_bound}
\begin{align}\label{SD_bound}
	\begin{split}
	\big\|\dif^{-1/2}(\u_{J} - \u_{J}^i)\big\|^2
	\le C^2 (\etalg)^2. 
	\end{split}
\end{align}
where $ C := (2 C_{\rm SD}^2(d+1)(J^2+1))^{1/2} $ only depends on 
$C_{\rm{SD}}$ of~\eqref{SD_2}, the space dimension $d \leq 3$, and at most linearly 
on the number of mesh levels $J$.
\end{proof}

\subsection{Two-level domain decomposition setting}

We now proceed to the proof of Theorem~\ref{thm_upper_bound} for the estimator 
built in Algorithm~\ref{AS_solver}. First, we present a few preparatory steps 
that will allow us to re-use the results we presented in the multilevel setting.
Recall from~Section~\ref{sec_DD_setting} that the fine mesh $\T_h$ is 
assumed to be obtained from the coarse mesh $\T_H$ after $J$ refinement 
steps.

Our goal is to write a multilevel presentation of the algebraic residual lifting
$\rhoa \in \Vaz$ computed by solving the subdomain 
problem \eqref{patch_problem_1}
on $\omaH$. Recall that $\omaH$ is the open subdomain
corresponding to the coarse-grid patch $\T_{0}^\ver$.
For the purpose of the analysis only, define the local MFE spaces on 
$\omaH$ associated with the intermediate mesh levels $\T_j$, $1 \le j \le J$, as
\begin{align}
	\bV_{j,H}^{\ver} & := \{ \bv_j \in \bV_j |_{\omaH},
	\ \bv_j {\cdot} \n = 0 \text{ on } \partial \omaH \}. \label{V_ja_MLDD}
\end{align}
In contrast to the spaces $\bVaj$ from Section~\ref{sec_ml_solver_setting},
which are defined on the $j$-level (small) patches $\omaj$,
we stress that the spaces $\bV_{j,H}^{\ver}$ are defined on the subdomains 
(large patches) $\omaH$ and their fine meshes $\T_j$; there are much fewer 
spaces $\bV_{j,H}^{\ver}$ than there are spaces $\bVaj$ but the spaces 
$\bV_{j,H}^{\ver}$ have much higher dimension than $\bVaj$.
Finally, define the divergence-free subspace of $\bV_{j,H}^{\ver}$,
\begin{align}
	\bV_{j,H}^{\ver,0} = \{\bvja \in \bV_{j,H}^{\ver},\ 
	\nabla {\cdot} \bvja = 0\}. 
	\label{V_ja_0_MLDD}
\end{align}

Let a coarse mesh vertex $\ver \in \V_{H}$ be fixed. We now consider an 
orthogonal
multilevel decomposition of $\rhoa \in \Vaz$ from~\eqref{patch_problem_1} 
on the subdomain
$\omaH$.
For $1 \le j \le J$, define for (analysis purposes, not constructed in practice) 
$\brho_{j,\ver}^i \in \bV_{j,H}^{\ver,0}$ as the solutions to
\begin{equation}\label{rho_j_compact_AS}
	(\dif^{-1}\brho_{j,\ver}^i,\bv_{j,\ver})_{\omaH} = 
	- (\dif^{-1}\u^i_{H},\bv_{j,\ver})_{\omaH} -
	\sum_{m=1}^{j-1}(\dif^{-1}\brho_{m,\ver}^i,\bv_{j,\ver})_{\omaH}  \quad
	\forall \, \bvja  \in  \bV_{j,H}^{\ver,0} ,\quad \forall 1 \le j \le J.
\end{equation}
Taking above $j=J$ and noting that $\bV_{J,H}^{\ver,0} = \Vaz$ in our notation, 
it follows that
$$
\sum_{j=1}^J(\dif^{-1}\brho_{j,\ver}^i,\va)_{\omaH} = 
- (\dif^{-1}\u^i_{H},\va)_{\omaH}
\quad \forall \va \in \Vaz,
$$
which, together with~\eqref{patch_problem_1}, implies that
\begin{equation}\label{rho_ML}
	\rhoa = \sum_{j=1}^J \brho_{j,\ver}^i,
\end{equation}
so that $\brho_{j,\ver}^i$ decompose $\rhoa$. 
Moreover, this decomposition is indeed orthogonal:
\begin{equation}\label{orthog_patch_1}
	(\dif^{-1}\brho_{j,\ver}^i,\brho_{m,\ver}^i)_{\omaH} = 0, 
	\ \ 1 \le j,m \le J, \, j \ne m,
\end{equation}
and
\begin{equation}\label{orthog_patch_2}
	\big\|\dif^{-1/2}\rhoa\big\|_{\omaH}^2 
	= \sum_{j=1}^J \big\|\dif^{-1/2}\brho_{j,\ver}^i\big\|_{\omaH}^2.
\end{equation}

With the above developments, the convergence analysis in the 
domain decomposition setting is a modification of the multigrid analysis 
presented in the previous Section~\ref{ML_proofs}. 
The key point is that the multilevel representation \eqref{rho_ML} allows us to 
use the stable decomposition from Proposition~\ref{SD_assumption}:

\begin{proof}[Proof of Theorem~\ref{thm_upper_bound} 
	(two-level domain decomposition setting)]
As in the multigrid case, note that the lower bound~\eqref{lower_bound} follows directly from the Pythagorean formula of Theorem~\ref{thm_error_contr} 
by using that $\big\|\dif^{-1/2}(\u_{J} - \u_{J}^{i+1})\big\| \ge 0$. 
It remains to show~\eqref{upper_bound}.
The proof here is more similar to \cite[Lemma~7.7]{Mir_Pap_Voh_19}, 
since the patches we are considering are bigger than the ones used in the 
stable decomposition of Proposition~\ref{SD_assumption}.
Thus, we first define, for every $j \in \{1, \ldots,  J \}$ and 
$\ver \in \V_{H}$, the set $\I_{\ver,
j} \subset \V_j$ containing vertices in $\T_j$ of the interior of the 
patch $ \omaH$ such that $\{\I_{\ver, j} \}_{\ver \in \V_{H}}$ cover $\V_j$ 
and are mutually disjoint. Let $\bv_{j,\ver} \in \bV_{j}^{\ver,0}$. 
This allows us to write
\begin{equation} \label{eq_sum}
	\sum_{\ver \in \V_j} \bv_{j,\ver} 
	= \sum_{\ver \in \V_{H}} \sum_{\tb \in \I_{\ver, j}} \bv_{j, \tb}.
\end{equation}
Moreover, since the vertices of $\I_{\ver, j}$ are localized in the interior of 
the patch $\omaH$, we have
\begin{equation} \label{eq_sum_incl}
	\sum_{\tb \in \I_{\ver, j}} \bv_{j, \tb} \in \bV_{j,H}^{\ver,0}.
\end{equation}
Once again, we use the stable decomposition \eqref{SD_1} applied to the error 
\begin{equation}\label{eq_SD_err_DD}
	\u_{J} - \u_{h}^i 
	= \bv_0 + \sum_{j=1}^J \sum_{\ver \in \V_{j}}  \bv_{j,\ver}
\end{equation}
and the construction of the multilevel solver
\begin{align*}
	\big\|\dif^{-1/2}(\u_{J} - \u_{h}^i)\big\|^2
	\stackrel{\eqref{mfe_div_free}}= {} & 
	- (\dif^{-1} \u_{h}^i, \u_{J} - \u_{h}^i ) 
	\stackrel{\eqref{eq_SD_err_DD}}= 
	- \Big(\dif^{-1} \u_{h}^i, \bv_0 
	+ \sum_{j=1}^J \sum_{\ver \in \V_{j}}  \bv_{j,\ver} \Big) \\
	\stackrel{\eqref{rho_0_DD}}= {} &
	(\dif^{-1}\brho_{H}^i,\bv_0)
	- \sum_{j=1}^J \sum_{\ver \in \V_{j}}(\dif^{-1} \u_{h}^i, 
	\bv_{j,\ver})_{\omaj}
	\\
	\stackrel{\eqref{u_0_i_DD}}= {} &
	(\dif^{-1}\brho_{H}^i,\bv_0)
	-  \sum_{j=1}^J \sum_{\ver \in \V_{j}}
	(\dif^{-1}  (\u_{H}^i - \brho_{H}^i), \bv_{j,\ver})_{\omaj}
	\\
	\stackrel{\eqref{eq_sum}}= {} &
	(\dif^{-1}\brho_{H}^i,\bv_0)
	-  \sum_{j=1}^J \sum_{\ver \in \V_{H}}
	\Big(\dif^{-1}  (\u_{H}^i - \brho_{H}^i), \sum_{\tb \in \I_{\ver,j}} 
	\bv_{j,\tb}\Big)_{\omaH}
	\\
	\stackrel{\substack{\eqref{eq_sum_incl}\\\eqref{rho_j_compact_AS}}}= {} &
	(\dif^{-1}\brho_{H}^i,\bv_0)
	+
	\sum_{j=1}^J \sum_{\ver \in \V_{H}}
	\Big\{
	\Big(\dif^{-1} \brho_{j,\ver}^i,   
	\sum_{\tb \in \I_{\ver,j}} \bv_{j,\tb} \Big)_{\omaH}\\
	{} &	
	+   \sum_{m=1}^{j-1} \Big(\dif^{-1} \brho_{m,\ver}^i,  
	\sum_{\tb \in \I_{\ver,j}} \bv_{j,\tb}\Big)_{\omaH}
	+ \Big(\dif^{-1} \brho_{H}^i,  \sum_{\tb \in \I_{\ver,j}} 
	\bv_{j,\tb}\Big)_{\omaH}
	\Big\} .
\end{align*}
Using Young's inequality and the stability of the decomposition~\eqref{SD_2}, 
one can follow the approach of~\cite[Lemma 7.7]{Mir_Pap_Voh_19}, cf.~\eqref{eq_MG_Y}--\eqref{SD_bound}, to show  
that there holds
\begin{equation}\label{SD_bound_AS}
	\big\|\dif^{-1/2} (\u_{J} - \u_{h}^i) \big\|^2 
	\le
	\widetilde{C}^2 \Big(\big\|\dif^{-1/2}\brho^i_{H}\big\|^2
	+  \sum_{j=1}^J 
	\sum_{\ver \in \V_{H}}
	\big\|\dif^{-1/2}\brho_{j,\ver}^i\big\|_{\omaH}^2 \Big),
\end{equation}
where $ \widetilde{C}^2$ only depends on $C_{\rm{SD}}$ of~\eqref{SD_2} 
and at most linearly on the number of mesh levels $J$.
Finally, we obtain the result
\begin{align*}
	\big\|\dif^{-1/2}(\u_{J} - \u_{h}^i)\big\|^2
	& \stackrel{\eqref{SD_bound_AS}}\le   \widetilde{C}^2 \Big(
	\big\|\dif^{-1/2}\brho_{H}^i\big\|^2 
	+   \sum_{j=1}^J \sum_{\ver \in \V_{H}} 
	\big\|\dif^{-1/2}\brho_{j,\ver}^i\big\|_{\omaH}^2 \Big)
	\\
	&
	\stackrel{\eqref{orthog_patch_2}}= \widetilde{C}^2 \Big(
	\big\|\dif^{-1/2}\brho^i_{H}\big\|^2 
	+  \sum_{\ver \in \V_{H}}  \big\|\dif^{-1/2}\rhoa \big\|_{\omaH}^2 \Big)
	\\
	& \stackrel{\eqref{lam_LB_DD}}\le 
	\widetilde{C}^2 (d+1) \Big( \big\|\dif^{-1/2}\brho^i_{H}\big\|^2
	+ \lambda_{h}^i \sum_{\ver \in \V_{H}}  
	\big\|\dif^{-1/2}\rhoa \big\|_{\omaH}^2
	\Big)\\
	& \stackrel{\eqref{loc_eta_DD}}= \widetilde{C}^2 (d+1) (\etalg)^2.
\end{align*}
\end{proof}

\section{Conclusions}\label{sec_conclusions}
In this work, we have presented an a-posteriori-steered multigrid solver and an 
a-posteriori-steered two-level domain decomposition method to solve iteratively 
an algebraic system originating from a saddle-point mixed finite element 
discretization of a second-order elliptic problem. 
The update in the solvers is constructed in such way that the norm 
of its components  
provides an a posteriori estimate of the algebraic error.
We proved that the a~posteriori estimators are 
efficient independently of the polynomial degree $p$ used in the discretization. 
This efficiency of the a~posteriori estimators 
leads to the equivalent result of the associated solvers 
contracting the algebraic error at each iteration, in particular independently 
of $p$. 
Though not the topic here, one may also treat local smoothing and dependency with respect the 
number of mesh levels $J$, see, e.g., 
Innerberger et al.~\cite{Inn_Mir_Praet_Strei_hp_MG_24} and the references 
therein. The idea consists in only solving the local problems associated to 
degrees of freedom that are new with respect to the previous mesh, which is 
pertinent in hierarchies leading to graded meshes.

\bibliographystyle{siam}
\bibliography{biblio}

\providecommand{\noopsort}[1]{}\def\polhk#1{\setbox0=\hbox{#1}{\ooalign{\hidewidth \lower1.5ex\hbox{`}\hidewidth\crcr\unhbox0}}} \def\cprime{$'$}
\begin{thebibliography}{10}

\bibitem{AlonsRodr_Cam_San_div_free_18}
{\sc A.~Alonso~Rodr\'{\i}guez, J.~Cama\~{n}o, E.~De~Los~Santos, and F.~Rapetti}, {\em A graph approach for the construction of high order divergence-free {R}aviart-{T}homas finite elements}, Calcolo, 55 (2018).
\newblock Paper No. 42, 28.

\bibitem{Amr_Ber_Dau_Gir_98_potentials3d}
{\sc C.~Amrouche, C.~Bernardi, M.~Dauge, and V.~Girault}, {\em Vector potentials in three-dimensional non-smooth domains}, Math. Methods Appl. Sci., 21 (1998), pp.~823--864.

\bibitem{Arnold_FEEC_18}
{\sc D.~N. Arnold}, {\em Finite Element Exterior Calculus}, SIAM, Philadelphia, PA, 2018.

\bibitem{Arn_Falk_Wint_MG_H_div_H_curl_00}
{\sc D.~N. Arnold, R.~S. Falk, and R.~Winther}, {\em Multigrid in {$H({\rm div})$} and {$H({\rm curl})$}}, Numer. Math., 85 (2000), pp.~197--217.

\bibitem{Bast_Beni_Voh_Yot_DD_MFE_25}
{\sc M.~Bastidas~Olivares, A.~Beni~Hamad, M.~Vohral{\'{\i}}k, and I.~Yotov}, {\em A posteriori algebraic error estimates and nonoverlapping domain decomposition in mixed formulations: energy coarse grid balancing, local mass conservation on each step, and line search}, Comput. Methods Appl. Mech. Engrg., 444 (2025), p.~118090.

\bibitem{Ben_Gol_Lie_05}
{\sc M.~Benzi, G.~H. Golub, and J.~Liesen}, {\em Numerical solution of saddle point problems}, Acta Numer., 14 (2005), pp.~1--137.

\bibitem{Boffi_Brezzi_Fortin_book}
{\sc D.~Boffi, F.~Brezzi, and M.~Fortin}, {\em Mixed finite element methods and applications}, vol.~44 of Springer Series in Computational Mathematics, Springer, Heidelberg, 2013.

\bibitem{Bren_MG_MFE_92}
{\sc S.~C. Brenner}, {\em A multigrid algorithm for the lowest-order {R}aviart-{T}homas mixed triangular finite element method}, SIAM J. Numer. Anal., 29 (1992), pp.~647--678.

\bibitem{Brenner_MFE_Solvers_09}
\leavevmode\vrule height 2pt depth -1.6pt width 23pt, {\em Fast solvers for mixed finite element methods}, in Mixed finite element technologies, vol.~509 of CISM Courses and Lect., Wien: Springer, 2009, pp.~57--88.

\bibitem{BrennerOhSung18_MGDarcy}
{\sc S.~C. Brenner, D.-S. Oh, and L.-Y. Sung}, {\em Multigrid methods for saddle point problems: {D}arcy systems}, Numer. Math., 138 (2018), pp.~437--471.

\bibitem{BDM_86}
{\sc F.~Brezzi, J.~Douglas, Jr., and L.~D. Marini}, {\em Recent results on mixed finite element methods for second order elliptic problems}, in Vistas in applied mathematics, Transl. Ser. Math. Engrg., Optimization Software, New York, 1986, pp.~25--43.

\bibitem{Brubeck_Farrell_hoFEM_21}
{\sc P.~D. Brubeck and P.~E. Farrell}, {\em A scalable and robust vertex-star relaxation for high-order {FEM}}, SIAM J. Sci. Comput., 44 (2022), pp.~A2991--A3017.

\bibitem{Cai_etal}
{\sc Z.~Cai, R.~R. Parashkevov, T.~F. Russell, J.~D. Wilson, and X.~Ye}, {\em Domain decomposition for a mixed finite element method in three dimensions}, SIAM J. Numer. Anal., 41 (2003), pp.~181--194.

\bibitem{Cant_DeTurck_Gluck_calc_top_3D_02}
{\sc J.~Cantarella, D.~DeTurck, and H.~Gluck}, {\em Vector calculus and the topology of domains in 3-space}, Amer. Math. Monthly, 109 (2002), pp.~409--442.

\bibitem{Chaum_Voh_Maxwell_equil_23}
{\sc T.~Chaumont-Frelet and M.~Vohral{\'{\i}}k}, {\em $p$-robust equilibrated flux reconstruction in {${\boldsymbol H}(\mathrm{curl})$} based on local minimizations. {A}pplication to a posteriori analysis of the curl--curl problem}, SIAM J. Numer. Anal., 61 (2023), pp.~1783--1818.

\bibitem{chaumontfrelet:hal-03817302}
\leavevmode\vrule height 2pt depth -1.6pt width 23pt, {\em A stable local commuting projector and optimal {$hp$} approximation estimates in {${\boldsymbol H}(\mathrm{curl})$}}, Numer. Math., 156 (2024), pp.~2293--2342.

\bibitem{Chav_Coh_Jaf_Dup_Rib_MFE2D_84}
{\sc G.~Chavent, G.~Cohen, J.~Jaffré, M.~Dupuy, and I.~Ribera}, {\em {Simulation of two-dimensional waterflooding by using mixed finite elements}}, Society of Petroleum Engineers Journal, 24 (1984), pp.~382--390.

\bibitem{Chen_Noch_Xu_MG_loc_ref_12}
{\sc L.~Chen, R.~H. Nochetto, and J.~Xu}, {\em Optimal multilevel methods for graded bisection grids}, Numer. Math., 120 (2012), pp.~1--34.

\bibitem{Cow_Man_Whee_BDD_MFE_95}
{\sc L.~C. Cowsar, J.~Mandel, and M.~F. Wheeler}, {\em Balancing domain decomposition for mixed finite elements}, Math. Comp., 64 (1995), pp.~989--1015.

\bibitem{dgs2023}
{\sc L.~Diening, L.~Gehring, and J.~Storn}, {\em Adaptive mesh refinement for arbitrary initial triangulations}, Found. Comput. Math.,  (2025).
\newblock DOI~10.1007/s10208-025-09698-7.

\bibitem{Dorfler_marking_96}
{\sc W.~D{\"o}rfler}, {\em A convergent adaptive algorithm for {P}oisson's equation}, SIAM J. Numer. Anal., 33 (1996), pp.~1106--1124.

\bibitem{Ern_Guermond_FEM_book_pt1}
{\sc A.~Ern and J.-L. Guermond}, {\em Finite elements {I}---{A}pproximation and interpolation}, vol.~72 of Texts in Applied Mathematics, Springer, Cham, 2021.

\bibitem{Ern_Guz_Potu_Voh_Poinc_disc_25}
{\sc A.~Ern, J.~Guzm{\'a}n, P.~Potu, and M.~Vohral{\'{\i}}k}, {\em Discrete {P}oincar\'e inequalities: a review on proofs, equivalent formulations, and behavior of constants}, IMA J. Numer. Anal.,  (2025).
\newblock DOI~10.1093/imanum/draf089.

\bibitem{Ern_Voh_p_rob_3D_20}
{\sc A.~Ern and M.~Vohral\'{\i}k}, {\em Stable broken {$H^1$} and {$H({\rm div})$} polynomial extensions for polynomial-degree-robust potential and flux reconstruction in three space dimensions}, Math. Comp., 89 (2020), pp.~551--594.

\bibitem{Ewing_Wang_92}
{\sc R.~E. Ewing and J.~Wang}, {\em Analysis of the {S}chwarz algorithm for mixed finite elements methods}, RAIRO Mod\'{e}l. Math. Anal. Num\'{e}r., 26 (1992), pp.~739--756.

\bibitem{Ewing_Wang_94}
\leavevmode\vrule height 2pt depth -1.6pt width 23pt, {\em Analysis of multilevel decomposition iterative methods for mixed finite element methods}, RAIRO Mod\'{e}l. Math. Anal. Num\'{e}r., 28 (1994), pp.~377--398.

\bibitem{Falk_Winther_25}
{\sc R.~S. Falk and R.~Winther}, {\em Local space-preserving decompositions for the bubble transform}, Found. Comput. Math.,  (2025).
\newblock DOI~10.1007/s10208-025-09700-2.

\bibitem{Gir_Rav_NS_86}
{\sc V.~Girault and P.-A. Raviart}, {\em Finite element methods for {N}avier-{S}tokes equations. {Theory} and algorithms}, vol.~5 of Springer Series in Computational Mathematics, Springer-Verlag, Berlin, 1986.

\bibitem{Glow_Whe_MFE_DD_88}
{\sc R.~Glowinski and M.~F. Wheeler}, {\em Domain decomposition and mixed finite element methods for elliptic problems}, in First {I}nternational {S}ymposium on {D}omain {D}ecomposition {M}ethods for {P}artial {D}ifferential {E}quations ({P}aris, 1987), SIAM, Philadelphia, 1988, pp.~144--172.

\bibitem{Hecht84}
{\sc F.~Hecht}, {\em Construction of a basis at free divergence in finite element and application to the {N}avier-{S}tokes equations}, in Numerical solutions of nonlinear problems ({R}ocquencourt, 1983), INRIA, Rocquencourt, 1984, pp.~284--297.

\bibitem{Heinrichs_88_line_relaxation_MG}
{\sc W.~Heinrichs}, {\em Line relaxation for spectral multigrid methods}, J. Comput. Phys., 77 (1988), pp.~166--182.

\bibitem{Hipt_Hopp_MFE3D_99}
{\sc R.~Hiptmair and R.~H.~W. Hoppe}, {\em Multilevel methods for mixed finite elements in three dimensions}, Numer. Math., 82 (1999), pp.~253--279.

\bibitem{Hipt_Wu_Zheng_cvg_adpt_MG_12}
{\sc R.~Hiptmair, H.~Wu, and W.~Zheng}, {\em Uniform convergence of adaptive multigrid methods for elliptic problems and {M}axwell's equations}, Numer. Math. Theory Methods Appl., 5 (2012), pp.~297--332.

\bibitem{Hiptmair_Xu_07}
{\sc R.~Hiptmair and J.~Xu}, {\em Nodal auxiliary space preconditioning in {${\bf H}({\bf curl})$} and {${\bf H}({\rm div})$} spaces}, SIAM J. Numer. Anal., 45 (2007), pp.~2483--2509.

\bibitem{Inn_Mir_Praet_Strei_hp_MG_24}
{\sc M.~Innerberger, A.~Mira\c{c}i, D.~Praetorius, and J.~Streitberger}, {\em {$hp$}-robust multigrid solver on locally refined meshes for {FEM} discretizations of symmetric elliptic {PDE}s}, ESAIM Math. Model. Numer. Anal., 58 (2024), pp.~247--272.

\bibitem{JayadharanKhattatovYotov21}
{\sc M.~Jayadharan, E.~Khattatov, and I.~Yotov}, {\em Domain decomposition and partitioning methods for mixed finite element discretizations of the {B}iot system of poroelasticity}, Comput. Geosci., 25 (2021), pp.~1919--1938.

\bibitem{kpp2013}
{\sc M.~Karkulik, D.~Pavlicek, and D.~Praetorius}, {\em On 2{D} newest vertex bisection: optimality of mesh-closure and {$H^1$}-stability of {$L_2$}-projection}, Constr. Approx., 38 (2013), pp.~213--234.

\bibitem{Kellogg_75}
{\sc R.~B. Kellogg}, {\em On the {P}oisson equation with intersecting interfaces}, Appl. Anal., 4 (1975), pp.~101--129.

\bibitem{Kirby_FE_de_Rham_14}
{\sc R.~C. Kirby}, {\em Low-complexity finite element algorithms for the de {R}ham complex on simplices}, SIAM J. Sci. Comput., 36 (2014), pp.~A846--A868.

\bibitem{Mathew_93}
{\sc T.~P. Mathew}, {\em Schwarz alternating and iterative refinement methods for mixed formulations of elliptic problems. {II}. {C}onvergence theory}, Numer. Math., 65 (1993), pp.~469--492.

\bibitem{Mir_Pap_Voh_19}
{\sc A.~Mira{\c c}i, J.~Pape{\v z}, and M.~Vohral{\'{\i}}k}, {\em A multilevel algebraic error estimator and the corresponding iterative solver with $p$-robust behavior}, SIAM J. Numer. Anal., 58 (2020), pp.~2856--2884.

\bibitem{Mir_Pap_Voh_lambda}
\leavevmode\vrule height 2pt depth -1.6pt width 23pt, {\em A-posteriori-steered $p$-robust multigrid with optimal step-sizes and adaptive number of smoothing steps}, SIAM J. Sci. Comput., 43 (2021), pp.~S117--S145.

\bibitem{Mir_Pap_Voh_W}
\leavevmode\vrule height 2pt depth -1.6pt width 23pt, {\em Contractive local adaptive smoothing based on {D}\"orfler's marking in a-posteriori-steered $p$-robust multigrid solvers}, Comput. Methods Appl. Math., 21 (2021), pp.~445--468.

\bibitem{Mitchell_91}
{\sc W.~F. Mitchell}, {\em Adaptive refinement for arbitrary finite-element spaces with hierarchical bases}, J. Comput. Appl. Math., 36 (1991), pp.~65--78.

\bibitem{Mitchell_10}
{\sc W.~F. Mitchell}, {\em The {$hp$}-multigrid method applied to {$hp$}-adaptive refinement of triangular grids}, Numer. Linear Algebra Appl., 17 (2010), pp.~211--228.

\bibitem{Ned_mix_R_3_80}
{\sc J.-C. N{\'e}d{\'e}lec}, {\em Mixed finite elements in {${\mathbb{R}}\sp{3}$}}, Numer. Math., 35 (1980), pp.~315--341.

\bibitem{PV22}
{\sc J.~Pape{\v z} and M.~Vohral{\'i}k}, {\em Inexpensive guaranteed and efficient upper bounds on the algebraic error in finite element discretizations}, Numer. Algorithms, 89 (2022), pp.~371--407.

\bibitem{Ra_Tho_MFE_77}
{\sc P.-A. Raviart and J.-M. Thomas}, {\em A mixed finite element method for 2nd order elliptic problems}, in Mathematical aspects of finite element methods (Proc. Conf., Consiglio Naz. delle Ricerche (C.N.R.), Rome, 1975), Springer, Berlin, 1977, pp.~292--315. Lecture Notes in Math., Vol. 606.

\bibitem{Rud_ful_adpt_MG_93}
{\sc U.~R{\"u}de}, {\em Fully adaptive multigrid methods}, SIAM J. Numer. Anal., 30 (1993), pp.~230--248.

\bibitem{Scheichl_PhD_00}
{\sc R.~Scheichl}, {\em Iterative solution of saddle point problems using divergence-free finite elements with applications to groundwater flow}, ProQuest LLC, Ann Arbor, MI, 2000.
\newblock Thesis (Ph.D.)--University of Bath (United Kingdom).

\bibitem{Schei_MFE_dec_03}
\leavevmode\vrule height 2pt depth -1.6pt width 23pt, {\em Decoupling three-dimensional mixed problems using divergence-free finite elements}, SIAM J. Sci. Comput., 23 (2002), pp.~1752--1776.

\bibitem{SchMelPechZag_08}
{\sc J.~Sch{\"o}berl, J.~M. Melenk, C.~Pechstein, and S.~Zaglmayr}, {\em Additive {S}chwarz preconditioning for {$p$}-version triangular and tetrahedral finite elements}, IMA J. Numer. Anal., 28 (2008), pp.~1--24.

\bibitem{Schoeberl_Zulehner_03}
{\sc J.~Sch\"{o}berl and W.~Zulehner}, {\em On {S}chwarz-type smoothers for saddle point problems}, Numer. Math., 95 (2003), pp.~377--399.

\bibitem{stevenson2008}
{\sc R.~Stevenson}, {\em The completion of locally refined simplicial partitions created by bisection}, Math. Comp., 77 (2008), pp.~227--241.

\bibitem{TakacsZulehner13_MGall}
{\sc S.~Takacs and W.~Zulehner}, {\em Convergence analysis of all-at-once multigrid methods for elliptic control problems under partial elliptic regularity}, SIAM J. Numer. Anal., 51 (2013), pp.~1853--1874.

\bibitem{Thomasset81_book}
{\sc F.~Thomasset}, {\em Implementation of finite element methods for {N}avier-{S}tokes equations}, Springer Series in Computational Physics, Springer-Verlag, New York-Berlin, 1981.

\bibitem{trax_97}
{\sc C.~T. Traxler}, {\em An algorithm for adaptive mesh refinement in {$n$} dimensions}, Computing, 59 (1997), pp.~115--137.

\bibitem{Weber_80}
{\sc C.~Weber}, {\em A local compactness theorem for {M}axwell's equations}, Math. Methods Appl. Sci., 2 (1980), pp.~12--25.

\bibitem{Wheeler_Yotov_00}
{\sc M.~F. Wheeler and I.~Yotov}, {\em Multigrid on the interface for mortar mixed finite element methods for elliptic problems}, Comput. Methods Appl. Mech. Engrg., 184 (2000), pp.~287--302.

\bibitem{wz2017}
{\sc J.~Wu and H.~Zheng}, {\em Uniform convergence of multigrid methods for adaptive meshes}, Appl. Numer. Math., 113 (2017), pp.~109--123.

\bibitem{Xu_Chen_Noch_opt_MG_loc_ref_09}
{\sc J.~Xu, L.~Chen, and R.~H. Nochetto}, {\em Optimal multilevel methods for {$H({\rm grad})$}, {$H({\rm curl})$}, and {$H({\rm div})$} systems on graded and unstructured grids}, in Multiscale, nonlinear and adaptive approximation, Springer, Berlin, 2009, pp.~599--659.

\end{thebibliography}
\end{document}